\documentclass[mnsc,nonblindrev]{informs3} 

\OneAndAHalfSpacedXI

\usepackage{natbib}
 \bibpunct[, ]{(}{)}{,}{a}{}{,}%
 \def\bibfont{\small}%

\defcitealias{hmis_data_dict}{U.S. HUD 2021}
\defcitealias{hud_pdr}{U.S. HUD 2007}
\defcitealias{lahsa_pit_ind}{LAHSA 2024a}
\defcitealias{lahsa_hic}{LAHSA 2024b}
\defcitealias{vi-spdat_change}{OrgCode and Community Solutions 2015}
\defcitealias{lahsa_hmis}{LAHSA 2015}
\defcitealias{lahsa2018report}{LAHSA 2018b}
\defcitealias{lahsa_intro_ces}{LAHSA 2018a}
\defcitealias{LACoC}{LAHSA 2023}
\defcitealias{lahsa-vi-spdat}{LAHSA 2017}

\usepackage{multirow}
\usepackage{booktabs}
\usepackage{siunitx} 
\usepackage{graphicx}

\usepackage{makecell}
\usepackage{enumerate}
\usepackage{enumitem}
\usepackage{bbm}
\usepackage{amsmath}
\usepackage{mathtools}
\usepackage{soul}
\usepackage[section]{placeins}
\usepackage{comment}
\usepackage{float}

\usepackage[table, dvipsnames]{xcolor}
\definecolor{tablecell}{RGB}{180,198,231}
\newcolumntype{s}{>{\columncolor{tablecell}} c}
\definecolor{blue}{rgb}{0,0,0}

\usepackage{pgfplots}
\usetikzlibrary{arrows.meta}

\pgfplotsset{compat=newest,
    width=6cm,
    height=3cm,
    scale only axis=true,
    max space between ticks=25pt,
    try min ticks=5,
    every axis/.style={
        axis y line=left,
        axis x line=bottom,
        axis line style={thick,->,>=latex, shorten >=-.4cm}
    },
    every axis plot/.append style={thick},
    tick style={black, thick}
}
\tikzset{
    semithick/.style={line width=0.8pt},
}
\usepgfplotslibrary{groupplots}
\usepgfplotslibrary{dateplot}

\usepackage{url}

\TheoremsNumberedThrough     
\ECRepeatTheorems

\EquationsNumberedThrough    

\MANUSCRIPTNO{MS-0001-1922.65}

\newcommand{\delck}[1]{{\color{magenta} [{\sffamily\bfseries DELETED CK:} #1]}}
\newcommand{\noteck}[1]{{\color{magenta} [{\sffamily\bfseries NOTE CK  :} {\em #1}]}}
\newcommand{\newck}[1]{{\color{magenta} #1}}

\newcommand{\newpv}[1]{{\color{teal} #1}}

\newcommand{\notebt}[1]{{\color{red} [{\sffamily\bfseries NOTE BT  :} {\em #1}]}}
\newcommand{\newbt}[1]{{\color{red} #1}}

\newcommand{\indep}{\perp \!\!\! \perp}

\begin{document}

\RUNAUTHOR{Tang et.\ al.}

\RUNTITLE{Learning Policies for Online Allocation of Scarce Resources}


\TITLE{{Learning Optimal and Fair Policies for Online \\ Allocation of Scarce Societal Resources \\ from Data Collected in Deployment}}

\ARTICLEAUTHORS{%
\AUTHOR{Bill Tang}
\AFF{USC CAIS Center for Artificial Intelligence in Society, \EMAIL{yongpeng@usc.edu}} 
\AUTHOR{Çağıl Koçyiğit}
\AFF{Luxembourg Centre for Logistics and Supply Chain Management, University of Luxembourg, \EMAIL{cagil.kocyigit@uni.lu}}
\AUTHOR{Eric Rice}
\AFF{USC CAIS Center for Artificial Intelligence in Society, \EMAIL{ericr@usc.edu}} 
\AUTHOR{Phebe Vayanos}
\AFF{USC CAIS Center for Artificial Intelligence in Society, \EMAIL{phebe.vayanos@usc.edu}} 
} 

\ABSTRACT{%

We study the problem of allocating scarce societal resources of different types (e.g., permanent housing, deceased donor kidneys for transplantation, ventilators) to heterogeneous allocatees on a waitlist (e.g., people experiencing homelessness, individuals suffering from end-stage renal disease, Covid-19 patients) based on their observed covariates. We leverage administrative data collected in deployment to design an online policy that maximizes expected outcomes while satisfying budget constraints, in the long run. Our proposed policy waitlists each individual for the resource maximizing the difference between their estimated mean treatment outcome and the estimated resource dual-price or, roughly, the opportunity cost of using the resource. Resources are then allocated as they arrive, in a first-come first-serve fashion. We demonstrate that our data-driven policy almost surely asymptotically achieves the expected outcome of the optimal out-of-sample policy under mild technical assumptions. We extend our framework to incorporate various fairness constraints. We evaluate the performance of our approach on the problem of designing policies for allocating scarce housing resources to people experiencing homelessness in Los Angeles based on data from the homeless management information system. In particular, we show that using our policies improves rates of exit from homelessness by~$5.16\%$ and that policies that are fair in either allocation or outcomes by race come at a very low price of fairness.

}%




\maketitle

%

\section{Introduction}\label{sec: intro}

We study the problem of allocating basic resources of different types to heterogeneous individuals within high-stakes social settings subject to budget and fairness constraints. We are particularly motivated by housing allocation for individuals experiencing homelessness in Los Angeles (LA) County. According to the Los Angeles Homeless Services Authority (LAHSA), more than {\color{blue}75,000 people are currently unsheltered in LA County as of June 2024~\citepalias{lahsa_pit_ind}}. By comparison, the latest Housing Inventory Count in {\color{blue}June 2024 from LAHSA reported only around 28,600 permanent housing unit interventions~\citepalias{lahsa_hic}, many of which are already occupied}. While permanent housing interventions are a key aspect of addressing homelessness ~\citepalias{hud_pdr}, the number of individuals experiencing homelessness far exceeds the capacity of available resources. Due to the scarcity of resources available relative to demand for them, there is a need to allocate housing in an effective way as individuals and housing resources arrive. 

Additionally, since housing resources are a public good with significant impact for individuals experiencing homelessness, it is important that an allocation policy is, in some sense, fair. A lack of fairness in the system may erode the trust and effectiveness of any policy. Indeed, in our motivating example, LAHSA policymakers and community members are particularly concerned with issues of fairness in allocation and outcomes within homelessness services. A 2018 report from LAHSA's Ad Hoc Committee on Black People Experiencing Homelessness found that while Black individuals represented~$9\%$ of the overall population in LA
County, they made up~$40\%$ of the population experiencing homelessness~\citepalias{lahsa2018report}. Furthermore, while housing resources are allocated at equal rates across racial groups, Black individuals receiving Permanent Supportive Housing (PSH) experienced higher rates of returns to homelessness relative to other racial groups. These findings have led to further investigations into causes of racial inequities in outcomes within the current system~\citep{cplPSHinequity} and a policy focus on racial equity to address existing injustices. Therefore, an allocation system needs to be not just effective, but also fair to be adopted and trusted by policymakers and community members.

Currently in the U.S., many local communities pool resources into a centralized planning system known as Continuum of Care (CoC), which coordinates housing and services funding for people experiencing homelessness. Within CoCs, Coordinated Entry Systems (CES), a network of homeless service providers and funders, use a standardized process to manage housing and supportive services and connect individuals to appropriate interventions. Individuals seeking housing in the LA CoC have various entry points into the system, such as access centers, street outreach, or emergency shelters. Following entry, they may take a self-reported survey to be assessed for eligibility and vulnerability. The assessment survey, known as Vulnerability Index--Service Prioritization Decision Assistance Tool (VI-SPDAT), consists of a series of questions to measure vulnerability, such as prior housing history and disabilities, and weights responses to output a score between 0 and 17, where a higher score indicates greater vulnerability~\citep{vi-spdat}. Generally, the score is used in determining and prioritizing the most vulnerable individuals to receive more supportive resources. Many communities, based on OrgCode recommendations, have used score cutoffs for prioritization or guidance in the past~\citepalias{lahsa_intro_ces}. Individuals scoring 8-17 are considered ``high risk'' and are prioritized for Permanent Supportive Housing ({PSH}), the most supportive and intensive type of housing resource. Individuals scoring 4-7 are recommended for Rapid-Rehousing ({RRH}), or short-term rental subsidies, while those scoring less than 4 are only eligible for services, which we refer to as Service Only ({SO})~\citepalias{vi-spdat_change}. Based on these scores and the VI-SPDAT assessments, case managers, who are social service workers assisting individuals experiencing homelessness, will determine the appropriate resource to match an individual to. 

While the above prioritization policy, based on just score cuts, is interpretable, both the thresholds and the weights on responses to construct the score are not tied to any data on intervention outcomes nor to arrival rates of individuals and resources. This leads to three issues. First, it is not clear that assigning the most vulnerable individuals {PSH}, or even {RRH}, would lead to the best overall outcomes, such as reducing returns to homelessness after intervention. Second, there may not be enough housing available to serve all individuals prioritized for that resource because the thresholds do not account for arrival rates. Indeed, even OrgCode, the organization that developed the VI-SPDAT, has called for alternative approaches beyond the VI-SPDAT in 2020~\citep{orgcode_change}. Finally, a threshold prioritization tool cannot address stakeholder concerns around racial disparities in housing outcomes. An allocation system that is equitable in outcomes would need to account for the heterogeneity of treatment effects by protected groups such as race, gender, or age~\citep{Jo2022}. In contrast to the existing tool, we seek to incorporate treatment outcomes and effects into the design of a policy for maximizing desirable societal outcomes and ensuring certain notions of fairness in allocation or outcomes.

Utilizing treatment outcomes and effects in policy design, however, is challenging since they are \textit{counterfactual}, or ``what if'' quantities. This means we need to know what \textit{would happen} to an individual if they were assigned to different resources. While running randomized control trials is the gold standard for inference of causal treatment effects, in high-stakes domains such as public health and social services, experimentation may be unethical or impractical. For example, it would be unethical to deny housing from individuals experiencing dangerous and vulnerable situations. A proper experiment may also take years before conclusive results, leading to negative life consequences for many individuals during that time span. Therefore, we aim to leverage observational data for learning treatment outcomes. The historical data that LAHSA and other communities have access to was collected in deployment from administered VI-SPDAT assessments, which contains an individual's covariates, and the LA County Homeless Management Information System (HMIS) database, which contains an individual's trajectory of interactions with the system. These interactions include events such as engagements with street outreach workers, emergency shelter access, and housing allocations, if any. From these trajectories, we construct a return to homelessness definition that is used as the outcome of interest.   

Motivated by these issues with the existing housing allocation system, our goal is to design an efficient policy for allocating scarce housing resources to maximize positive expected outcomes, subject to budget and fairness constraints. Using historical observational data, we learn treatment outcomes based on an arriving individual's characteristics, such as prior housing history, to appropriately match individuals to housing resources. In particular, {\color{blue} in keeping with the existing allocation system, the assignment policy must make decisions about which resources to allocate to individuals as the resources arrive in the system, rather than waiting for a sufficient number of individuals and housing resources to aggregate in the system for offline matching. Stakeholders, including people experiencing homelessness, policy-makers, and the public at large,} are dissatisfied with scarce housing sitting empty for long periods of time {\color{blue} while thousands of people who need them sleep on the streets or in shelters each night. Furthermore, offline matching over fixed time intervals results in undesirable behavior where individuals with the same characteristics can be recommended different resources due to variations between different batches. When two individuals with the same characteristics receive different resources, this allocation does not seem fair nor interpretable to policymakers and individuals seeking assistance. This concern is captured by the notion of \emph{individual fairness} introduced by \cite{dwork2012fairness}, which calls for treating similar individuals similarly. A policy that assigns the same resources based on individual characteristics ensures consistency and interpretability in allocation recommendations. Additionally, implementing dynamic policies that adjust allocations based on the current state of the system is infeasible in this context, as the information necessary, such as the state of an individual in the system at any point in time, is not systematically tracked.}

While we primarily focus on housing, our approach is applicable in general to allocation of scarce treatments in other public domains where system outcome efficiency and fairness are important. One example is the problem of kidney transplantations in which patients arrive to the system and join a waiting pool for a donor organ to be procured~\citep{Trichakis2013, Dickerson2015}. Similar to our case, there is a significant organ shortage relative to the number of waiting patients and observational data exists from the United Network for Organ Sharing. Another example is the allocation of scarce healthcare equipment such as ventilators and critical care beds to patients arriving to the hospital with Covid-19 during spikes in hospitalizations. Observational data from electronic health records for hospitalized patients with Covid-19 can be used to learn policies for allocating equipment as patients arrive~\citep{johnston2020, Radovanovic2021}.


\subsection{Contributions of the Paper}\label{sec: contributions}
 
Our contributions in this paper are:
 \begin{itemize}
 \item {\color{blue}We first consider policy design under budget constraints only. {We propose a sample-based dual-price queuing policy by solving an offline policy learning problem using only historical (observational) data and propose an online queuing-based implementation of the derived policy.} This policy establishes separate queues for each treatment, and each individual is waitlisted for the treatment that maximizes the difference between their estimated expected outcome, given their covariates, and the estimated resource dual price, which approximates the opportunity cost of assigning that treatment. Resources are allocated as they arrive in the system to the first individual {in the queue for that treatment}. The proposed policy ensures that capacity constraints are always met, meaning resources are assigned only when they are available in the system. {We prove that, under mild technical assumptions, our policy attains the same long‑run expected average outcome as the best offline policy that already knows the true distribution of covariates and (counterfactual) treatment outcomes.}
 We establish a connection between our technical assumptions on the convergence of estimators and the standard assumptions commonly made in observational data studies, arguing that estimators satisfying our assumptions exist when the standard assumptions for observational data hold. Our framework unifies policy learning from observational data and dual-price allocation methods.}
 \item We extend our framework to incorporate various fairness constraints, including statistical parity in allocations and outcomes between groups of a protected characteristic (e.g., race, gender, or age). We introduce fairness-constrained extensions of the sample-based dual-price queuing policy{\color{blue}, which adjust the estimated opportunity costs for different groups using the estimated dual variables associated with the fairness constraints}. While these extensions are designed to meet the respective fairness constraints in-sample, we demonstrate through numerical experiments that they can also improve fairness out-of-sample. 
 \item We demonstrate through computational experiments on synthetic and real data that our method's guarantees of asymptotic optimality lead to improvements in out-of-sample performance, even in finite samples. Specifically, we evaluate the policy's performance in allocating scarce housing resources to people experiencing homelessness in Los Angeles, based on data from the HMIS database. {\color{blue}We show that our policies improve historical rates of exits from homelessness by~$5.16\%$ (which amounts to approximately 500 additional individuals exiting homelessness per year), and policies that are fair in either allocation or outcomes by race come at a very low price of fairness.}
 \end{itemize}

%

\subsection{Related Works}\label{sec: relatedworks}

Our work is closely related to the literature on data-driven allocation of public resources and fairness, policy learning from observational data, network revenue management, and online matching. 

\textbf{Data-Driven Allocation of Public Resources \color{blue}{and Fairness}.} A closely related strand of literature focuses on data-driven allocation in scarce public resource settings such as kidney allocation{\color{blue}~\citep{Trichakis2013, Dickerson2015, li2019incorporating}}, hospital ICUs~\citep{grand2022robustness}, {\color{blue}refugee resettlement~\citep{bansak2018improving, ahani2021placement, ahani2023dynamic, freund2023group, bansak2024outcome}}, or homeless services~\citep{Javad2018, Kube2019, Rahmattalabi2022, kaya2022improving, kube2023fair}. {\color{blue}Our work shares similarities with several papers in this line of literature (e.g., \cite{bansak2018improving, ahani2021placement, kube2023fair, freund2023group}) in that we use historical data to learn treatment outcomes and incorporate these estimates into resource allocation optimization problems. However, in contrast to this stream of literature, we prove the asymptotic optimality of our proposed sample-based policy, explicitly formulating the necessary properties that the estimators must satisfy and discussing their relationship to standard assumptions in observational data studies
Additionally, some of the related works, including \cite{bansak2018improving} and \cite{kube2023fair}, formulate the corresponding resource allocation problems as mixed-integer linear optimization problems for offline matching, whereas our approach only requires solving a finite linear program to derive an online implementation. Furthermore, we consider various fairness notions, which have not been studied by the majority of these works.} 

{\color{blue}The works most closely related to ours are~\cite{Trichakis2013}, \cite{li2019incorporating}, \cite{ahani2023dynamic}, and \cite{freund2023group}, as they, like our framework, leverage dual variables in constructing resource allocation decisions.}
\cite{Trichakis2013} propose learning a scoring policy for kidney allocation that satisfies fairness constraints, but assumes knowledge of the treatment outcomes. Their scoring rule is a linear approximation of treatment outcomes minus the dual multipliers of various eligibility and fairness constraints. {\color{blue}\cite{li2019incorporating} and \cite{ahani2023dynamic} propose allocations to maximize dual variable adjusted outcomes, where they estimate the dual variables based on simulated future arrival trajectories. In comparison, we only solve a single finite linear program to estimate the dual variables, provide theoretical guarantees for our proposed method, and incorporate fairness constraints. Concurrent to our work, \cite{freund2023group} tackle group fairness constraints in dynamic refugee settlement and propose dual-price assignment decisions, relying on dual variables associated with capacity and fairness constraints. They demonstrate that the regret of their proposed approach vanishes as the number of refugee cases grows. While they assume that estimates of treatment outcomes (i.e., employment probabilities) are directly available to decision-makers, our asymptotic optimality result explicitly incorporates the estimation process from data, establishing asymptotic optimality as the number of training samples increases. {Furthermore, while their fairness notion generalizes our notion of statistical parity in outcomes in Section~\ref{sec: sp in out}, our notion of statistical parity in allocation in Section~\ref{sec: sp in alloc} addresses a distinct but complementary fairness notion.} Another important difference is that we search for a policy that is a function of an individual's covariates—therefore treating individuals with the same covariates equally—whereas they propose a matching that allows refugees with the same covariates to be treated differently.} 

In another related work,~\cite{BHATTACHARYA2012} focus on maximizing welfare under a budget constraint for social programs with a binary treatment by learning asymptotically consistent treatment effects and propose an allocation policy based on data from a randomized experiment. In comparison to their binary treatment setting, we allow for multiple scarce treatment types, are able to impose a variety of fairness constraints, and formulate an online implementation in which an individual is matched to an appropriate resource as they arrive in the system. 

{\color{blue}Our work is also broadly related to fairness in decision-making within high-stakes domains. Many such problems incorporate fairness constraints that limit some notion of disparity between protected groups~\citep{Trichakis2013, Corbett-Davies2017, Javad2018, Nguyen2021, Rahmattalabi2022}. Recent work has shown that different notions of fairness studied in high-stakes domains may be impossible to satisfy simultaneously~\citep{Mashiat2022, Jo2022}. Other fairness concepts, e.g., based on social welfare, ordinal efficiency and envy-freeness, have been explored in algorithmic division~\citep{Rahmattalabi2021, Nguyen2021, Freeman2020, manshadi2021fair} and mechanism design~\citep{athanassoglou2011house}.}  Our work follows in the vein of imposing fairness constraints such as parity in resource allocation rates between different racial groups, as this is how policymakers within LA homeless services evaluate fairness~\citepalias{lahsa2018report}.

\textbf{Policy Learning from Observational Data.} Another strand of literature focuses on policy evaluation and learning from observational data by estimating counterfactuals or correcting for selection bias in historical treatment assignments. Various works focus on constructing unbiased estimators for policy evaluation from observational data under an unconfoundedness assumption, which can then be used for policy learning. These include the direct method of using plug-in estimates of counterfactual treatment outcomes~\citep{Qian2011}, propensity score based reweighting strategies~\citep{Jswaminathan15a}, and doubly robust combinations of outcome estimates and propensity-based reweighting~\citep{dr_dudik}.~\cite{Kitagawa2018} and~\cite{Athey2021} study regret bounds for policy learning under propensity-weighted or doubly robust objective functions relative to some well-specified policy class with structure while~\cite{Zhou2022} extend the analysis to the multiple treatment case. Other works emphasize interpretable policy classes such as prescriptive trees~\citep{kallus2017, Bertsimas2019, optprestrees, Zhou2022}. While these works derive important analytical results for policy classes with suitable structure, they either use a specified class like trees in practice, or suggest a general policy learning algorithm without explicitly considering how to optimize over the abstract policy class that may have budget, fairness, or functional form constraints. 
In contrast, we consider the general policy class of measurable functions and derive an explicit assignment policy parameterized by sample-based dual multipliers that is asymptotically optimal and can flexibly handle capacity and fairness constraints.
{\color{blue}In other words, by considering a specific policy design problem and connecting it to the dual-price controls framework, we are able to characterize an asymptotically optimal, interpretable, and computationally cheap policy across all (unrestricted) policy classes.}

\textbf{Network Revenue Management.} Another related stream of literature is network revenue management (NRM), which focuses on the sale of products comprised of scarce resources to customers who request a specific product with an offer price~\citep{talluri2004}. In domains such as airlines and hotels, the system must decide to accept or reject an arriving customer's offer and the decision depends on the stochastic demand and revenue for products relative to resource capacities. In particular, our approach is inspired by work on bid-price controls, a class of policies where an offer is accepted only if the revenue generated is greater than the sum of optimal dual variables (bid-prices) associated with the resource capacity constraints~\citep{talluri1998}. Early work by~\cite{talluri1998} showed that such a control policy is asymptotically optimal as the number of requests and resources grow as long as the relative relationship of supply and demand stays constant. Since NRM problems are often formulated as dynamic programs, recent works have focused on settings where large state spaces cannot be solved tractably and computing optimal bid-prices is challenging. Tractable solutions are generally heuristics, approximate dynamic programming, or Lagrangian relaxation methods~\citep{Adelman2007, kunnumkal2016, Zhang2017, Li2022}. While our approach is inspired by bid-price controls within NRM, our setting differs in that the system must choose amongst a set of resources to match an arrival to\, rather than an accept or reject decision for a specific requested resource.


{\color{blue}
\textbf{Online Matching and Resource-Constrained Dynamic Problems.} Our work also relates to the class of online matching problems \citep{karp1990optimal, mehta2007adwords, aggarwal2011online, mehta2013online, devanur2013randomized, vera2019bayesian}, where one must decide which resource, among a set of capacitated resources, to immediately and irrevocably assign to each arriving individual, which results in a reward but also uses a resource. The goal is to match individuals and resources in an online fashion to maximize rewards while respecting resource capacities. Similar to our approach, online matching algorithms make matching decisions by trading off between the potential match reward and the opportunity cost of using up a scarce resource as in the case of the AdWord problem \citep{mehta2007adwords}. The derivation of the tradeoff function often results from analyzing the offline or full knowledge primal-dual linear programs where the dual variables are used to parameterize the opportunity cost of using a resource. While in many online matching settings it is possible, and even desirable, for individuals with the same features to be matched to different resources, this behavior is undesirable in our setting, which is why we adopt a policy design framework rather than a matching one. When two individuals with the same features are assigned different resources, this allocation does not seem fair nor interpretable to policymakers and individuals seeking assistance. Therefore, it is important in our motivating application that two individuals with the same features receive the same type of resource.}

{\color{blue}Additionally, methodologically, our work is related to more general resource-constrained dynamic problems~\citep{balseiro2023survey}, which capture some of the online matching settings and network revenue management problems as special cases. Specifically, our work is closely related to the fluid certainty-equivalent control results~\citep{balseiro2023survey}, which leverage dual-based optimality conditions. However, our framework is based on offline policy learning translated into an online/dynamic implementation, and we work with historical observational data generated under an unknown historical policy, whereas they assume that the underlying distribution of contexts and their relation to rewards are known.}

\subsection{Organization and Notation}
The remainder of this paper is structured as follows. Section~\ref{sec: problem set up}  introduces the policy design problem and outlines the technical assumptions required for our optimality results. 
{\color{blue}Section~\ref{sec: sample-based} studies the proposed sample-based dual-price queuing policy, along with its optimality guarantee.}
Section~\ref{sec: fairness extensions} studies the fairness-constrained extensions of the dual-price queuing policy. Finally, Section~\ref{sec: emp results} presents numerical results based on synthetic and real data. All proofs are relegated to the online appendix.


{\it Notation.} All random variables in this paper are defined as measurable functions on an abstract probability space~$(\Omega, \mathcal F, \mathbb P)$ and are capitalized ($e.g.$,~$X$) while their realizations are denoted by the same symbols in lowercase letters ($e.g.$,~$x$) unless stated otherwise. Throughout the paper, (in)equalities containing random variables should be interpreted to hold~$\mathbb P$ almost surely. We let~$\mathcal L_\infty(\mathcal X, \mathcal C)$ denote the set of all bounded Borel-measurable functions from a Borel set~$\mathcal X$ to a Borel set~$\mathcal C$. For any set~$\mathcal{A}$, we use~$\vert \mathcal{A} \vert$ to denote the cardinality of the set. We use~$\mathbbm{1}[\cdot]$ to denote the indicator function of a logical expression~$E$, where~$\mathbbm{1}[E] = 1$ if~$E$ is true and~$\mathbbm{1}[E] = 0$ otherwise.

\section{Problem Statement, Problem Formulation, and Assumptions} \label{sec: problem set up}

{\color{blue}

\subsection{Problem Statement and Problem Formulation}

We now formalize our problem statement. To model long-run policy performance, we consider an infinite time horizon  wherein individuals indexed in the set $\mathcal K$ arrive over time. Each individual $k\in \mathcal K$ is characterized by their \emph{covariates} $X_k \in \mathcal X \subseteq \mathbb R^d$ and their \emph{potential outcomes} $\{Y^t_k\}_{t \in \mathcal T}$, where $\mathcal T = \{0,1,\ldots, m\}$ collects the treatment options, with treatment 0 representing no-treatment (or the control), and $Y^t_k \in \mathbb R$ corresponds to the individual's outcome under treatment $t$, see e.g.,~\cite{hernan2023causal} for an overview of the potential outcomes framework in causal inference. For example, in the context of the housing allocation problem that motivates us, the vector of covariates $X_k$ may collect the VI-SPDAT and HMIS information of homeless client $k$; the set of treatments $\mathcal T$ may collect the various housing options such as PSH, RRH, or no resource; and the potential outcome $Y_k^t$ represents, for instance, whether client $k$ would exit homelessness or not if assigned resource $t$. In particular, we regard both $X_k$ and the $Y_k^t$s, observed at individual $k$'s \emph{time of arrival}, as fixed quantities over time. This modeling choice is realisitic in our motivating application if the time from when individuals arrive to the system and when they receive a resource is not too long. Crucially, we emphasize that an individual's covariates and potential outcomes \emph{can change over a long period of time}, which motivates our online implementation proposed in Section~\ref{sec: sample-based}. We numerically observe that our policy has low wait times with minimal impact on individual outcomes in Sections~\ref{sec: policy_queue_waittimes} and~\ref{sec: waittime simulation app}. Without much loss of generality we assume that $\{(X_k,\{Y_k^t\}_{t\in \mathcal T} )\}_{k\in \mathcal K}$ are independent and identically distributed (i.i.d.) across individuals and therefore may drop the index $k$ from now on.

We do not know the joint distribution of $(X,\{Y^t\}_{t\in \mathcal T} )$. However, we have access to a historical dataset of i.i.d.\ observations which (with a slight abuse of notation) we denote by $\{(X_i, T_i, Y_i)\}_{i\in \mathcal N}$ and index in the set $\mathcal N$. Here $X_i$ is the covariate vector of the $i$th observation, $T_i$ is the historical treatment assigned to it, and $Y_i = Y_i^{T_i}$ is the \emph{observed outcome}, i.e., the outcome under the treatment received.  Critically, we do not know and cannot control the historical treatment assignment policy and the outcomes $Y_i^t$ for $t \neq T_i$ remain unobserved.

We aim to use this historical observational dataset to design a policy~$\pi \in \mathcal L_\infty(\mathcal X, [0,1]^{m+1})$ for assigning treatments to different individuals conditional on their covariates. Specifically,~$\pi^t(x)$ denotes the probability of assigning treatment~$t$ to an individual with covariates~$x$, so that~$\sum_{t \in \mathcal{T}} \pi^t(x) = 1$ for all~$x \in \mathcal{X}$. We focus on policies that are functions of the covariates only (being e.g., independent of the state of the system) for explainability and fairness purposes. Indeed, aligned with the notion of individual fairness \citep{dwork2012fairness}, such policies ensure that all individuals that are similarly situated are treated similarly and make it easier to justify why a particular allocation decision was made.

The objective of the policy design problem is to maximize the asymptotic expected average outcome of assigning treatments to arriving individuals. We can formulate the expected average outcome under a policy~$\pi$ as
\begin{equation*}
    \begin{aligned}
     \lim_{\vert \mathcal{K} \vert \to \infty} \mathbb E \left[ \frac{1}{\vert \mathcal{K} \vert} \sum_{k \in \mathcal K} \sum_{t \in \mathcal T} \pi^t(X_k) Y_k^t  \right] = \lim_{\vert \mathcal{K} \vert \to \infty} \mathbb E \left[ \frac{1}{\vert \mathcal{K} \vert} \sum_{k \in \mathcal K} \sum_{t \in \mathcal T}\pi^t(X_k) \mathbb E [Y_k^t \vert X_k]  \right]
     = \mathbb E \left[\sum_{t \in \mathcal T} \pi^t(X) m^t(X) \right], 
    \end{aligned}
\end{equation*}
where $\mathbb E$ denotes the expectation operator with respect to the (unknown) joint distribution $\mathbb P$ of $(X, T, \{ Y^t \}_{t \in \mathcal T})$ and~$m^t \in \mathcal L_\infty(\mathcal X, \mathbb R$) is the conditional average treatment outcome under treatment~$t \in \mathcal T$ for covariate vector~$x \in \mathcal X$, i.e., $m^t(x) = \mathbb E [Y^t \vert X=x]$. The last equality follows from the definition of~$m^t$ and the assumption that individuals arrive i.i.d.

In the scare resource settings that motivate us, each treatment type~$t \in \mathcal T \setminus \{0\}$ has limited availability, i.e., the arrival rate of treatments of type~$t \in \mathcal T \setminus \{0\}$ is less than the arrival rate of individuals. Specifically, for treatment~$t \in \mathcal T$, we let~$b^t \in [0,1]$ denote the asymptotic capacity of treatment~$t$ per individual. For example, if~$b^t = 0.3$, then we can allocate treatment~$t$ to at most~$30\%$ of all individuals. Since treatment~$0$ is the {no-treatment} case, we set~$b^0 = 1$. We then require policy~$\pi$ to satisfy the capacity constraints 
\begin{equation*}\label{eq: almost-surely capacity constraint}
    \begin{aligned}
        \lim_{\vert \mathcal{K} \vert \to \infty} \frac{1}{\vert \mathcal{K} \vert} \sum_{k \in \mathcal K} \pi^t(X_k) \leq b^t \quad \forall t \in \mathcal T.
    \end{aligned}
\end{equation*}
As individuals arrive i.i.d. 
and~$\pi$ is bounded, by the strong law of large numbers we have
\begin{equation*}
    \begin{aligned}
        \lim_{\vert \mathcal{K} \vert \to \infty} \frac{1}{\vert \mathcal{K} \vert} \sum_{k \in \mathcal K} \pi^t(X_k) = \mathbb{E} [\pi^t(X)] \quad \forall t \in \mathcal T.
    \end{aligned}
\end{equation*}
The policy design problem is thus given by 
\begin{equation}\label{eq:UpperBoundTrueP}\tag{$\mathcal{P}$}
    \begin{aligned}
    z^\star = &\max_{\pi \in \Pi} &&\mathbb E \left[\sum_{t \in \mathcal T} \pi^t(X) m^t(X) \right]\\
    &\;\;\,\text{s.t.} && \mathbb E\left[ \pi^t(X) \right] \leq b^t \quad\forall t \in \mathcal T ,
    \end{aligned}
\end{equation}
where~$$\Pi = \left\{\pi \in L_\infty(\mathcal X, {[0,1]^{m+1}}) \,:\, \sum_{t \in \mathcal T} \pi^t(x) = 1 \quad\forall x \in \mathcal X \right\}.$$
Without loss of generality, we take~$m^t$ to be non-negative for all~$t \in \mathcal T$. In fact, if~$m^t$ did take on negative values, one could shift the values of all~$m^t$ by adding a sufficiently large constant to all of them. This does not change the optimal solution to \eqref{eq:UpperBoundTrueP} and simply shifts the optimal objective value by a constant.

Note that we cannot solve problem \eqref{eq:UpperBoundTrueP} because the functions $m^t$ and the distribution $\mathbb P$ are unknown, and we do not have access to samples from this distribution due to the observational nature of the data. Moreover, policy $\pi$ in Problem~\eqref{eq:UpperBoundTrueP} is an \emph{offline} policy in the sense that it only depends on the personal characteristic of individuals when they arrive in the system. However, it must be implemented \emph{online}, since in our setting, both individuals and resources arrive over time. In Section~\ref{sec: sample-based}, we will propose a method for learning an asymptotically optimal policy for Problem~\eqref{eq:UpperBoundTrueP} from the observational data. We will also introduce an online implementation of this policy that maintains the same performance. In the remainder of this section, we introduce the assumptions underlying our approach, highlighting why they are reasonable in our setting.
 
}

\subsection{Assumptions}\label{sec: assumptions}

\subsubsection{Observational Data and Identifiability Assumptions}\label{subsec: causal assumption}

A fundamental challenge of policy learning from observational data is that, for any sample~$i$, we observe the outcome~$Y_i = Y_i^{T_i}$ associated with the treatment~$T_i$ assigned by the historical policy but we do not observe the outcome that would have been obtained if any other treatment~$t \neq T_i$ had been assigned.
Therefore, even evaluating the performance of a \emph{counterfactual} policy~$\pi$ necessitates appropriate \emph{identifiability} conditions to make causal quantities such as $\mathbb{E}[Y^t]$ and $m^t$ possible to estimate from the available data. We make the following assumptions that are standard within the causal inference literature \citep{hernan2023causal}.
%
\begin{assumption}\label{ass: causal}
The following conditions hold.
\begin{itemize}
    \item[(i)] Conditional Exchangeability:~$Y^t \indep T \mid X \quad \forall t \in \mathcal{T}$
    
    \item[(ii)] Positivity:~$\mathbb{P}(T = t \mid X=x) > 0 \quad \forall t \in \mathcal{T}$ and~$\forall x \in \mathcal X$ 

    \item[(iii)] Consistency:~$Y^{T} = Y$ 
    
\end{itemize}
\end{assumption}
Assumption~\ref{ass: causal} (i) says that the outcome~$Y^t$ of assigning treatment~$t \in \mathcal T$ and the allocated treatment~$T$ are conditionally independent given the covariates~$X = x$. Equivalently, the joint distribution of potential outcomes is the same for all individuals in each treatment group~$t \in \mathcal{T}$ given they have the same value of~$x$. This assumption is satisfied if we measure and condition on all covariates that are common causes of treatment assignment and outcomes. This implies that we can consider all individuals with the same~$x$ in each treatment group as exchangeable and we can calculate~$m^t(x)$ as~$\mathbb{E}[Y \mid T=t, X=x]$. In context of our motivating example, this means the VI-SPDAT assessment responses and HMIS database contain all the relevant information for individuals experiencing homelessness that case managers use to make their final decisions on housing assignments. In our case, this assumption is likely reasonable since the VI-SPDAT is used as a prioritization tool by case managers and contains a comprehensive overview of an individual's vulnerability factors.

Assumption~\ref{ass: causal} (ii) requires that the conditional probability of treatment assignment given features~$X = x$ is positive for each treatment. Intuitively, this means that each individual receives any particular treatment with positive probability under the historical policy. Assumption~\ref{ass: causal} (iii) states that the observed outcome under an observed treatment~$T$ is actually the potential outcome. Practically, this means that our treatments are well specified or else if there were multiple versions of a particular treatment~$t$, then~$Y^t$ is not well defined. 

\subsubsection{Conditional Average Treatment Outcomes Assumptions}\label{sec: CATE} 
{\color{blue}Our asymptotic optimality results 
rely on the following technical assumption on conditional average treatment outcomes, which ensures that a dual-price queuing policy—constructed under full distributional knowledge and \emph{not} from data—is optimal in~\eqref{eq:UpperBoundTrueP}.}
\begin{assumption} \label{ass: meanoutcomes}
For any~$t \in \mathcal T$, the following conditions hold.
\begin{enumerate}[label=(\roman*)] 
  \item \label{ass: meanoutcomes i} The random variable~$m^t(X)$ is continuously distributed and has a bounded support, i.e., there exists a~$C \in \mathbb{R}_+$ such that~$|m^t(x)| \leq C$ for all ~$x \in \mathcal{X}$.
  
  \item \label{ass: meanoutcomes ii} For any~$t' \in \mathcal{T}$ and~$t' \neq t$,~$m^t(X) - m^{t'}(X)$, the difference in conditional treatment effects between~$t$ and~$t'$, is continuously distributed.
\end{enumerate}
\end{assumption}
Assumption~\ref{ass: meanoutcomes} requires that the conditional average treatment outcomes~$m^t(X)$ for all treatments~$t \in \mathcal{T}$, and the treatment effect differences~$m^t(X) - m^{t'}(X)$ between any two treatments~$t$,~$t'$ are bounded and continuously distributed. {\color{blue}Note that boundedness of $m^t(X)$ for all~$t \in \mathcal{T}$ implies boundedness of~$m^t(X) - m^{t'}(X)$ for any~$t$,~$t'$.} Intuitively, this assumption holds if treatment outcomes and effects take on continuous values. 

The following example illustrates a simple case where the conditional average treatment outcomes are linear functions of a continuous random variable and satisfy Assumption~\ref{ass: meanoutcomes}. 
\begin{example} [Assumptions~\ref{ass: meanoutcomes}~\ref{ass: meanoutcomes i} and ~\ref{ass: meanoutcomes ii}  hold] 
Suppose that there exists a single covariate~$X$ that is continuously distributed on the support~$[0,1]$. There are two potential treatments indexed by~$0$ and~$1$, and their respective conditional average treatment outcome functions are given by~$m^0(x) = x$ and~$m^1(x) = 2x$. In this case, the random variables~$m^0(X)$,~$m^1(X)$,~$m^1(X) - m^0(X)$ and~$m^0(X) - m^1(X)$ are all bounded and continuously distributed.
\end{example}

We also show in the following two examples wherein conditions (i) and (ii) in Assumption~\ref{ass: meanoutcomes} are both needed, thereby demonstrating that one does not imply the other. These also give intuition as to when Assumption~\ref{ass: meanoutcomes} does not hold.

\begin{example} [Assumption~\ref{ass: meanoutcomes}~\ref{ass: meanoutcomes i} does not imply Assumption~\ref{ass: meanoutcomes}~\ref{ass: meanoutcomes ii}]
Suppose again that there exists a single covariate~$X$ that is continuously distributed on the support~$[0,1]$. There are two potential treatments indexed by~$0$ and~$1$, and their respective conditional average treatment outcome functions are given by~$m^0(x) = x$ and~$m^1(x) = x + 1$. The random variables~$m^0(X)$ and~$m^1(X)$ are therefore bounded and continuously distributed, which implies that Assumption~\ref{ass: meanoutcomes}~\ref{ass: meanoutcomes i} holds. On the other hand, ~$m^1(x) - m^0(x) = 1$, and therefore~$m^1(X) - m^0(X)$ is not continuously distributed. Thus, Assumption~\ref{ass: meanoutcomes}~\ref{ass: meanoutcomes ii} does not hold.
\end{example}

\begin{example} [Assumption~\ref{ass: meanoutcomes}~\ref{ass: meanoutcomes ii} does not imply Assumption~\ref{ass: meanoutcomes}~\ref{ass: meanoutcomes i}] 
Suppose again that there exists a single covariate~$X$ that is continuously distributed on the support~$[0,1]$. There are two potential treatments indexed by~$0$ and~$1$, and their respective conditional average treatment outcome functions are given by~$m^0(x) = 1$ and~$m^1(x) = x$. Then, we have that~$m^1(x) - m^0(x) = x - 1$ and~$m^0(x) - m^1(x) = 1 - x$, and so both~$m^1(X) - m^0(X)$ and~$m^0(X) - m^1(X)$ are bounded and continuously distributed. However, ~$m^0(X)$ is constant valued in this case and not continuously distributed. Therefore, Assumption~\ref{ass: meanoutcomes}~(ii) does not imply Assumption~\ref{ass: meanoutcomes}~(i).
\end{example}

\subsubsection{Sample-Based Estimates Assumptions}\label{subsec: sample-based-assumption}
In Section~\ref{sec: sample-based}, we propose a sample-based dual-price queuing policy. In order to learn this policy from data, we first need to estimate the conditional treatment outcomes. 
Let~$\hat{m}^t(\cdot \vert \{(X_i, T_i, Y_i)\}_{i \in \mathcal{N}}) \in \mathcal L_{\infty}(\mathcal{X}, \mathbb{R})$ 
denote an estimator of the conditional mean treatment response~$m^t$ from the data. 
For ease of notation, we suppress the dependency of~$\hat{m}^t$ on the 
historical samples and simply write~$\hat{m}_n^t$ to stress that we use a sample of size~$n$. 
The following assumption formally states properties we need our estimators to have and how ``good'' they need to be asymptotically. 
\begin{assumption}\label{ass: estimator of m^t}
For any~$t \in \mathcal T$, the following condition holds\newck{.} 
\begin{itemize}
    \item[(i)] The estimate~$\hat{m}_n^t$ is measurable and has bounded support, i.e., there exists a~$\hat{C} \in \mathbb{R}_+$ such that~$\vert \hat{m}_n^t(x) \vert \leq \hat{C}$ for all~$x \in \mathcal{X}$ and~$n \in \mathbb{N}$.  
    \item[(ii)] The estimate~$\hat{m}_n^t$ converges uniformly and almost surely to~$m^t$, i.e.,~$\lim_{n \rightarrow \infty} \sup_{x \in \mathcal{X}} |\hat{m}_n^t(x) - m^t(x)| = 0$.
\end{itemize}
\end{assumption}
Assumption~\ref{ass: estimator of m^t} (i) is a technical assumption needed for our theoretical results and is standard within other policy learning works as well. Assumption~\ref{ass: estimator of m^t} (ii), the uniform convergence of~$\hat{m}_n^t$ to~$m^t$ for each treatment, is necessary for the asymptotic optimality of the sample-based dual-price policy; see Theorem~\ref{theorem: asymptotic property of sample-based dual-price queuing policy}. This is because we do not make any functional assumptions about the conditional expected outcomes~$m^t$,~$t \in \mathcal T$.  In fact, if one has knowledge of the form of~$m^t$, for example that it is a linear function with parameters~$\theta$, then Assumption~\ref{ass: estimator of m^t} can be weakened to uniform convergence of the estimator parameters~$\theta$. In our general setting, the asymptotic convergence guarantees we require in Assumption~\ref{ass: estimator of m^t} (ii) are satisfied by non-parametric estimators like nearest-neighbors, kernel regression, and sieve estimators~\citep{mack1982weak, cheng1984strong, Liero1989, hansen2004uniform, chen2015optimal, xie2023uniform} if standard assumptions on observational data also hold. {\color{blue} Specifically, under Assumption~\ref{ass: causal}, the conditional average treatment outcomes~$m^t$ for~$t \in \mathcal{T}$, which are \textit{causal} quantities, can be expressed and estimated as \textit{statistical} quantities that are solely a function of the observed data~\citep{hernan2023causal}. When Assumption~\ref{ass: causal} holds and we can estimate a causal quantity as a statistical quantity from observed data, then the aforementioned non-parametric estimators will satisfy the uniform convergence property of Assumption~\ref{ass: estimator of m^t} (ii) ~\citep{Liero1989}. 
}

\section{Dual-Price Queuing Policy: Solution Approach and Online Implementation}
\label{sec: sample-based}
{\color{blue} In this section, we characterize an online policy, which we call the sample-based dual-price queuing policy, and show that it asymptotically attains the optimal value~$z^\star$ of problem \eqref{eq:UpperBoundTrueP} as the number of data samples tends to infinity. To that end, we need the conditional average treatment outcome functions~$m^t$ for each~$t \in \mathcal{T}$, which we estimate from our observational dataset. There exist a variety of methods within the causal inference literature for learning conditional average treatment effects and outcomes such as metalearning, double machine learning, and doubly robust learning~\citep{chernozhukov2018double, Kunzel2019, Foster2019, Kennedy2020}. Since these methods allow for arbitrary estimators to be used for modeling outcomes and treatment assignments, which are components used to adjust for observational selection bias in the estimation procedure, we can use nonparametric estimators that satisfy the theoretical asymptotic assumptions we make in Assumption~\ref{ass: estimator of m^t}.

To showcase the approach that yields the characterization of our policy, consider the following sample approximation of problem~\eqref{eq:UpperBoundTrueP}, where we replace the expectations with sample averages and~$m^t$'s with their estimators~$\hat{m}_n^t$'s.
\begin{equation}\label{eq: sample-based UpperBoundTrueP}\tag{$\mathcal{\hat{P}}$}
    \begin{aligned}
    \hat{z}_n^\star = &\max_{{\pi} \in \Pi} &&\frac{1}{n} \sum_{i=1}^n \sum_{t \in \mathcal T} {\pi}^t(x_i) \hat{m}_n^t(x_i)\\
    &\;\;\,\text{s.t.} && \frac{1}{n} \sum_{i=1}^n {\pi}^t(x_i) \leq b^t \quad\forall t \in \mathcal T
    \end{aligned}
\end{equation}}
{\color{blue}The Lagrangian dual of this problem is given by
\begin{equation}\label{eq:sample-based DualOfTrueP}
    \begin{aligned}
        \hat{\nu}_n^\star &= \min_{\mu \in \mathbb R_+^{m+1}} \;\max_{\pi \in \Pi}\;  \frac{1}{n} \sum_{i=1}^n \sum_{t \in \mathcal T} {\pi}^t(x_i) \hat{m}_n^t(x_i) + \sum_{t \in \mathcal T} \mu^t \left( b^t - \frac{1}{n} \sum_{i=1}^n {\pi}^t(x_i) \right)\\
        &= \min_{\mu \in \mathbb R_+^{m+1}} \;\max_{\pi \in \Pi}\; \frac{1}{n} \sum_{i=1}^n  \sum_{t\in \mathcal T} \pi^t(x_i) (\hat{m}_n^t(x_i) - \mu^t) + \sum_{t\in \mathcal T} \mu^t  b^t,
    \end{aligned}
\end{equation}
where~$\mu$ collects the dual variables of the capacity constraints. By weak duality, we have~$\hat{\nu}_n^\star \geq \hat{z}_n^\star$. For any~$\mu \in \mathbb R_+^{m+1}$, the inner maximization problem is solved by the policy
\begin{equation}\label{eq: sample-based policy from Lagrangian relaxation}
    \begin{aligned}
    \pi^t(x) = \begin{cases}
    1 &\text{if}\quad t = \min \arg \max_{t' \in \mathcal T}\; \left( \hat{m}_n^{t'}(x) - \mu^{t'} \right)\\
    0 &\text{otherwise}
    \end{cases}
    \quad\forall t \in \mathcal T, \,\forall x \in \mathcal X.
    \end{aligned}
\end{equation}
We suppress the dependence of policy~$\pi$ on~$\mu$ notationally in order to avoid clutter.
We use a lexicographic tie-breaker in characterizing \eqref{eq: sample-based policy from Lagrangian relaxation}. We emphasize that our results do not rely on this particular tie-breaking rule. By substituting policy \eqref{eq: sample-based policy from Lagrangian relaxation} into the dual problem \eqref{eq:sample-based DualOfTrueP}, we obtain}
\begin{equation}\label{eq: sample-based assignment simplified}\tag{$\mathcal{\hat{D}}$}
\begin{aligned}
\hat{\nu}^\star_{n} = \min_{{\mu} \in \mathbb{R}_+^{m+1}} \; \frac{1}{n}  \sum_{i=1}^n \max_{t \in \mathcal T}\left( \hat{m}_n^t(x_i) - \mu^t \right) + \sum_{t \in \mathcal T} \mu^t  b^t\newpv{.}
\end{aligned}
\end{equation}
Problem~\eqref{eq: sample-based assignment simplified} is a convex optimization problem and can be reformulated as a finite linear program that can easily be solved by off-the-shelf solvers.  Denote by~$\hat{\mu}_{n}^\star$ an optimal solution to problem~\eqref{eq: sample-based assignment simplified}, which exists under Assumption~\ref{ass: estimator of m^t} (i) because we can restrict the feasible set of \eqref{eq: sample-based assignment simplified} to a compact set without loss of generality; see Lemma~\ref{lemma: sample compact solution set alt}. 
{\color{blue}We define a candidate policy~$\hat{\pi}_{n}^{\star}$ through}
\begin{equation}\label{eq: sample policy from dual}
    \begin{aligned}
    \hat{\pi}_{n}^{\star,t}(x) = \begin{cases}
    1 &\text{if}\quad t = \min \arg \max_{t' \in \mathcal T}\; (\hat{m}_n^{t'}(x) - \hat{\mu}_{n}^{\star, t'})\\
    0 &\text{otherwise}
    \end{cases}
    \end{aligned}
    \quad\forall t \in \mathcal T, \,\forall x \in \mathcal X.
\end{equation}
{\color{blue}Note that policy~$\hat{\pi}_{n}^{\star}$ is obtained by choosing~$\mu$ as~$\hat{\mu}_{n}^\star$ in \eqref{eq: sample-based policy from Lagrangian relaxation}. 

Next, we introduce an online implementation of policy~$\hat{\pi}_{n}^{\star}$ that takes into account the fact that individuals arrive sequentially and that treatments become available over time. By construction, this policy \emph{always} satisfies the capacity constraints, even out-of-sample and in finite horizon.}

\begin{definition}[Sample-based dual-price queuing policy]
\label{def: sample dual-price queuing policy}
Establish~$m$ queues each associated with a treatment~$t \in \mathcal{T}\setminus\{0\}$. At any time, individuals in queue~$t$ are waitlisted for treatment~$t$. There are two events triggering an action: (1) arrival of an individual, (2) arrival of a treatment.
\begin{itemize}
    \item[(1)] When an individual with covariate vector~$x$ arrives, compute the treatment~$t = \min \arg \max_{t' \in \mathcal T}\; (\hat{m}_n^{t'}(x) - \hat{\mu}_{n}^{\star, t'})$, where~$t$ attains the maximum difference between the estimated conditional mean treatment outcome and the associated sample optimal dual variable. If~$t \neq 0$, assign them to queue~$t$, otherwise assign them to {no-treatment}. If they are the only individual waiting in queue~$t$, check whether there is an available treatment of type~$t$. If there is, assign treatment~$t$ to the individual. 
    \item[(2)] When a treatment of type~$t$ becomes available, assign this treatment to the first individual in the queue~$t$. If there is no one in the queue~$t$, keep the treatment until an individual is assigned to this queue. 
\end{itemize}
\end{definition}

{\color{blue}In finite horizon, the sample-based dual-price queuing policy assigns~$t=0$ to individuals still waiting in queues at the end of the time horizon. Since in practice our experiments can only occur in the finite horizon setting, we will use a `long' test horizon in the synthetic data experiments as a proxy for the asymptotic environment.}
We also note that in general, even if the historical samples are fixed,~$\hat{\mu}_{n}^\star$, the optimal solution to problem~\eqref{eq: sample-based assignment simplified}, need not be unique since the objective function of~\eqref{eq: sample-based assignment simplified} may not be strictly convex and therefore the sample-based dual-price queuing policy may not be uniquely specified. 

In the following, we denote by~$\hat{z}^{\rm{D}}_{n}$ the expected average outcome of the sample-based dual-price queuing policy, where the expectation is taken with respect to the distribution of covariates and outcomes of the individuals arriving during implementation and \emph{not} with respect to the distribution of the historical samples. Note that~$\hat{z}^{\rm{D}}_{n}$ is a random value because it relies on the historical samples, which are random. The next theorem shows that the sample-based dual-price queuing policy attains the optimal value~$z^\star$ of problem \eqref{eq:UpperBoundTrueP} as the number of data samples tends to infinity.

\begin{theorem}\label{theorem: asymptotic property of sample-based dual-price queuing policy}
Under Assumptions~\ref{ass: meanoutcomes} and~\ref{ass: estimator of m^t}, the expected average outcome~$\hat{z}^{\rm{D}}_{n}$ of the sample-based dual-price queuing policy is almost surely asymptotically at least as high as the optimal value~$z^\star$ of problem \eqref{eq:UpperBoundTrueP} as~$n \rightarrow \infty$, that is, 
\begin{equation*}
    \begin{aligned}
   \limsup_{n \rightarrow \infty} z^\star - \hat{z}^{\rm{D}}_{n} \leq 0.
    \end{aligned}
\end{equation*}
\end{theorem}
{\color{blue}To prove Theorem~\ref{theorem: asymptotic property of sample-based dual-price queuing policy}, we first construct a dual-price queuing policy assuming that the joint distribution of the covariates~$X$ and the potential outcomes~$Y^t$,~$t \in \mathcal T$, is known. We show that this dual-price queuing policy is optimal in problem~\eqref{eq:UpperBoundTrueP}; see~\ref{sec: known values} in the online appendix for further details about this construction and the corresponding result. This first part of the proof relies on establishing some structural properties of problem~\eqref{eq:UpperBoundTrueP} and applying standard techniques based on optimality conditions from network revenue management~\citep{talluri1998}, and more generally, from capacity-constrained dynamic programming~\citep{balseiro2023survey}. The idea for the rest of the proof}
is that if the allocations of a sample-based dual-price queuing policy asymptotically match those of a dual-price queuing policy almost surely, then the sample-based dual-price queuing policy will asymptotically generate expected average outcomes at least as good as those of a dual-price queuing policy almost surely. Since the dual-price queuing policy is optimal in problem~\eqref{eq:UpperBoundTrueP}, then asymptotically~$\hat{z}^{\rm{D}}_{n}$ is least as high as~$z^\star$ almost surely. We show in the online appendix that the optimal solutions to problem~\eqref{eq: sample-based assignment simplified},~$\hat{\mu}_{n}^\star$, will be arbitrarily close to an optimal solution~$\mu^\star$ of {\color{blue}the dual problem of~\eqref{eq:UpperBoundTrueP}}
as the number of samples increases,~$n \to \infty$. This will imply that the allocations of a sample-based dual-price queuing policy and dual-price queuing policy match almost surely.

\section{Fairness Extensions} \label{sec: fairness extensions}
In high-stakes settings, the allocation of scarce resources needed to satisfy basic needs raises issues of fairness with respect to some protected group(s). This is especially the case when minority groups historically discriminated against due to race, gender, or other protected characteristics disproportionately make up the population supported by the allocation system. In our motivating example, fairness and equity, or the lack thereof, are key concerns for policymakers and community members (see Section~\ref{sec: intro}).

In this section, we showcase fairness-constrained extensions of the {\color{blue}sample-based dual-price queuing policy} 
introduced in {\color{blue}Section~\ref{sec: sample-based}}
based on two fairness notions of statistical parity in allocation and statistical parity in outcomes. In the following subsections, we introduce the aforementioned fairness notions and present the resulting fair policies. Details of the derivations of the respective policies are relegated to~\ref{sec: fairness extensions derivations}. The proposed extensions and their derivations can be similarly derived for other fairness notions that can be defined as linear constraints such as conditional statistical parity, which, for example, could require that the percentage of housing resources allocated to individuals with low VI-SPDAT scores should not exceed a predefined threshold. Alternatively, one could incorporate a fairness constraint that enforces lower bounds on the percentage of resources allocated to specific protected groups, which, for example, could require that the percentage of housing allocated to minority groups should be at least as high as the percentage allocated historically. For a comprehensive overview of possible fairness notions and constraints for scarce resource allocation systems and impossibility results between these various fairness notions, see~\cite{Jo2022}.

{\color{blue}
Recall from Section~\ref{sec: problem set up} that~$X$ can include protected features~$G \in \mathcal G$ that represent sensitive attributes (e.g., race, gender, sex), i.e.,~$X = (X^{-G}, G)$, where~$X^{-G}$ collects all features of~$X$ excluding~$G$. Therefore we note that the protected feature~$G$ is also an input into our policy~$\pi$ in the following subsections.}



\subsection{Statistical Parity in Allocation} \label{sec: sp in alloc}

\emph{Statistical parity in allocation}, which requires that the treatment allocation probabilities are similar across protected groups, is an intuitive and commonly used notion by stakeholders for evaluating the fairness of high-stakes allocation systems~\citep{Jo2022}. Formally, a policy~$\pi$ satisfies statistical parity in allocation if
{\color{blue}
\begin{equation}\label{eq: sp alloc}
    \begin{aligned}
    \mathbb E \left[ \pi^t(X) \mid G=g \right] - \mathbb E \left[ \pi^t(X) \mid G=g' \right] \leq \delta \quad\forall g, g' \in \mathcal{G},\; g \neq g',\; \forall t \in \mathcal T,
    \end{aligned}
\end{equation}}
which requires that the allocation probabilities of any treatment~$t \in \mathcal T$ are approximately equal (up to a specified tolerance level~$\delta$) across different sensitive groups. Problem \eqref{eq:UpperBoundTrueP} can be extended to include the additional constraint \eqref{eq: sp alloc}. {\color{blue}Similarly, the sample approximation \eqref{eq: sample-based UpperBoundTrueP} of~\eqref{eq:UpperBoundTrueP} can be extended to include a sample approximation of the constraint~\eqref{eq: sp alloc}, i.e.,
\begin{equation}\label{eq: sample-based sp alloc}
    \begin{aligned}
    \frac{1}{n_g} \sum_{i=1}^n \mathbbm{1}[g_i = g] \pi^t(x_i) - \frac{1}{n_{g'}} \sum_{i=1}^n \mathbbm{1}[g_i = g'] \pi^t(x_i) \leq \delta \quad\forall g, g' \in \mathcal{G},\; g \neq g',\;\forall t \in \mathcal T,
    \end{aligned}
\end{equation}
where the expectations are replaced with sample averages, and $n_g$ denotes the number of individuals with protected feature $g$ in the data.}

As in 
{\color{blue}Section~\ref{sec: sample-based}}, let~$\mu \in \mathbb R_+^{m+1}$ collect the dual variables of the capacity constraints in \eqref{eq: sample-based UpperBoundTrueP}. We now denote by~$\lambda^t(g,g') \in \mathbb R_+$ 
the dual variable of the fairness constraint \eqref{eq: sample-based sp alloc} for the pair~$(g,g')$ and for treatment~$t$. Following similar steps to those in {\color{blue}Section~\ref{sec: sample-based}},
we characterize a policy 
{\color{blue}
\begin{equation*}\label{eq: policy from dual_sp}
    \begin{aligned}
    \hat{\pi}^{\star, t}_{\rm alloc, n} (x) = \begin{cases}
    1 &\text{if}\;\;\; t = \min \arg \max_{t' \in \mathcal T}\; \left( \hat{m}_n^{t'}(x) - \hat{\mu}_n^{\star,{t'}} - \frac{n}{n_{g}} \hat{\gamma}^{\star,t'}_n(g) \right)\\
    0 &\text{otherwise}
    \end{cases} 
    \end{aligned}
    \;\;\; \forall t \in \mathcal T, \;\; x \in \mathcal{X},
    \end{equation*}
where for each~$g \in \mathcal G$ and~$t \in \mathcal{T}$,~$\hat{\gamma}^{\star, t}_n(g)= \sum_{g' \in \mathcal G, g \neq g'} (\hat{\lambda}_n^{\star, t}(g,g') - \hat{\lambda}_n^{\star, t}(g', g))$ is an aggregation of dual variables across all~$g' \neq g$, and~$\hat{\mu}_n^{\star, t}$ and~$\hat{\lambda}_n^{\star, t}(g,g')$ are optimal solutions to the Lagrangian dual problem of \eqref{eq: sample-based UpperBoundTrueP} with additional fairness constraints \eqref{eq: sample-based sp alloc}.} Intuitively, if~$\hat{\mu}_n^{\star, t}$ is roughly interpreted as an opportunity cost/threshold for assigning treatment~$t$, then our new policy uses~$\hat{\gamma}_n^{\star, t}(g)$ to adjust that threshold for each group to achieve allocation parity. We can see that the adjustment~$\hat{\gamma}_n^{\star, t}(g)$ is unrestricted in sign so that it can appropriately increase (to allocate less) or decrease (to allocate more) the `opportunity cost' of each treatment for each group. Finally, our adjustment is also inversely weighted by the population proportion of each group~$g$ to appropriately account for the protected group distribution within the population. 

Using the definition of~$\hat{\pi}^{\star}_{\rm alloc,n}$, we can construct a fairness-constrained sample-based dual-price queuing policy similarly to the Section~\ref{sec: sample-based}. The fairness-constrained sample-based dual-price queuing policy assigns an individual with covariates~$x$ to the queue for treatment~$t = \min \arg \max_{t' \in \mathcal T}\; \left(\hat{m}_n^{t'}(x) - \hat{\mu}_n^{\star,{t'}} - \frac{n}{n_{g}} \hat{\gamma}^{\star,t'}_n(g) \right)$.


In certain allocation settings, policymakers may want to prioritize minority groups historically disadvantaged by racial, gender, or other forms of discrimination for allocation as a solution to mitigate inequities. While statistical parity in allocation ensures equal opportunities across protected groups, it may be insufficient for repairing existing disparities. Whether this prioritization is appropriate again depends on the fairness goals policymakers want to achieve. In this case, we can define a set~$\mathcal G^{\rm min}$ of minority groups and a set~$\mathcal G^{\rm maj}$ of majority groups such that~$\mathcal G^{\rm min}$ and~$\mathcal G^{\rm maj}$ are disjoint and~$\mathcal G^{\rm min} \cup \mathcal G^{\rm maj} = \mathcal G$. The goal of prioritizing allocation to minority groups could be formulated through the constraint
{\color{blue}
\begin{equation}\label{eq:sp_alloc__minority} 
    \begin{aligned}
    \mathbb E \left[ \pi^t(X) \mid G=g \right] \leq \mathbb E \left[ \pi^t(X) \mid G=g' \right]  \quad\forall g \in \mathcal G^{\rm maj}, \; g' \in \mathcal G^{\rm min}, \;\forall t \in \mathcal T \setminus \{0\}, 
    \end{aligned}
\end{equation}}
where the average treatment assignment probability for any minority group~$g' \in \mathcal G^{\rm min}$ is at least as high as the average treatment assignment probability for any majority group~$g \in \mathcal G^{\rm maj}$. We see that constraint \eqref{eq:sp_alloc__minority} is a variant of \eqref{eq: sp alloc}, where we have excluded the {no-treatment} option of~$t=0$, set~$\delta = 0$, and we enforce the constraint only for the pairs~$(g, g')$ such that~$g \in \mathcal G^{\rm maj}$ and~$g' \in \mathcal G^{\rm min}$. {\color{blue}We can again write a sample approximation of constraint \eqref{eq:sp_alloc__minority}, include it in \eqref{eq: sample-based UpperBoundTrueP}, and associate dual variables $\lambda^t(g, g') \in \mathbb R_+$, for $g \in \mathcal G^{\rm maj}$, $g' \in \mathcal G^{\rm min}$, $t \in \mathcal T \setminus \{0\}$, with these constraints. Following similar steps as before, we can characterize a minority prioritizing allocation policy $\hat{\pi}^{\star}_{\rm alloc,n}$ as follows. For any $x \in \mathcal X$ such that~$g \in \mathcal G^{\rm min}$, we have
\begin{equation*}\label{eq: policy from dual_sp_minority_maj_group}
    \begin{aligned}
    \hat{\pi}^{\star, t}_{\rm alloc,n} (x) = \begin{cases}
    1 &\text{if}\;\;\; t = \min \arg \max_{t' \in \mathcal T}\; \left(\hat{m}_n^{t'}(x) - \hat{\mu}_n^{\star, t'} + \frac{n}{n_g} \sum_{g' \in \mathcal G^{\rm maj}}\hat{\lambda}_n^{\star, t'} \left( g', g \right) \right)\\
    0 &\text{otherwise}
    \end{cases}
    \end{aligned}
    \;\;
    \begin{array}{r}
         \forall t \in \mathcal T,
    \end{array}   
\end{equation*}
and for any $x \in \mathcal X$ such that~$g \in \mathcal G^{\rm maj}$, we have
\begin{equation*}\label{eq: policy from dual_sp_minority_min_group}
    \begin{aligned}
    \hat{\pi}^{\star, t}_{\rm alloc,n} (x) = \begin{cases}
    1 &\text{if}\;\;\; t = \min \arg \max_{t' \in \mathcal T}\; \left( \hat{m}_n^{t'}(x) - \hat{\mu}_n^{\star, t'} - \frac{n}{n_g} \sum_{g' \in \mathcal G^{\rm min}} \hat{\lambda}_n^{\star, t'} \left( g, g' \right)\right)\\
    0 &\text{otherwise}
    \end{cases}
    \end{aligned}
    \;\;
    \begin{array}{r}
         \forall t \in \mathcal T. 
    \end{array}
\end{equation*}
For~$t = 0$, we define~$\hat{\lambda}_n^{\star, 0}(g, g') = 0$ for any pair of~$g$ and~$g'$.} 
In this case we see that the adjustments for all minority groups~$g \in \mathcal G^{\rm min}~$, given by~$\sum_{g' \in \mathcal G^{\rm maj}} \hat{\lambda}_n^{\star, t'}(g', g)$, are non-negative, and serve only to decrease the `threshold' of treatment~$t$ for group~$g$, therefore leading to more assignments. On the other-hand, since~$\sum_{g' \in \mathcal G^{\rm min}} \lambda^{\star, t'}(g, g')$ is also non-negative but has opposite sign, it serves to increase the `threshold' of treatment~$t$ for all~$g \in \mathcal G^{\rm maj}$, thereby leading to fewer assignments. {\color{blue}Similarly to before, we can define a corresponding minority allocation prioritization sample-based dual-price queuing policy.}

\subsection{Statistical Parity in Outcomes} \label{sec: sp in out}

While statistical parity in allocation ensures similar chances of receiving each treatment across different protected groups, it may not be adequate for achieving fair outcomes for the individuals. If minority groups are more vulnerable to homelessness due to 
discrimination, then statistical parity in allocation may not address existing outcome inequities. This is precisely the case and a major concern in our motivating example, where {PSH} is allocated at nearly equal rates to all racial groups, but Black individuals experience higher rates of returns to homelessness compared to other racial and ethnic groups~\citepalias{lahsa2018report}. For this reason, we will consider an alternative fairness notion called \emph{statistical parity in outcomes}, which requires that a policy~$\pi$ satisfies constraints of the form
{\color{blue}
\begin{equation}\label{eq: sp_outcome_constraint} 
    \begin{aligned}
    \mathbb E \left[ \sum_{t \in \mathcal T} \pi^t(X) m^t(X) \mid G=g \right] - \mathbb E \left[ \sum_{t \in \mathcal T} \pi^t(X) m^t(X) \mid G=g' \right] \leq \delta 
    \quad \forall g, g' \in \mathcal{G},\;\; g \neq g'.
    \end{aligned}
\end{equation}}
Constraint \eqref{eq: sp_outcome_constraint} ensures that the expected outcomes under the allocation policy~$\pi$ are approximately equal across different protected groups. {\color{blue}The sample approximation of constraint \eqref{eq: sp_outcome_constraint} takes the following form.
\begin{equation}\label{eq: sample-based sp_outcome_constraint} 
    \begin{aligned}
\frac{1}{n_g} \sum_{i=1}^n \sum_{t \in \mathcal T} \mathbbm{1}[g_i = g] \pi^t(x_i)\hat{m}_n^t(x_i) - \frac{1}{n_{g'}} \sum_{i=1}^n \sum_{t \in \mathcal T} \mathbbm{1}[g_i = g'] \pi^t(x_i)\hat{m}_n^t(x_i) \leq \delta  \quad\forall g, g' \in \mathcal{G},\;\; g \neq g'
    \end{aligned}
\end{equation}
Problem~\eqref{eq: sample-based UpperBoundTrueP} can be extended to include this constraint. As before, we associate dual variables with constraints of this extended problem and denote by~$\lambda(g,g') \in \mathbb R_+$ the dual variable of the statistical parity in outcome constraint~\eqref{eq: sample-based sp_outcome_constraint} for the group pair~$g$ and~$g'$. We denote by~$\hat{\mu}_n^\star$ and~$\hat{\lambda}_n^\star$ the optimal solutions to the Langragian dual problem of \eqref{eq: sample-based UpperBoundTrueP} with additional statistical parity in outcomes constraints~\eqref{eq: sample-based sp_outcome_constraint}. Then, we characterize a modified policy
\begin{equation}\label{eq: policy from dual_sp_out 1}
    \begin{aligned}
    \hat{\pi}^{\star, t}_{\rm out,n} (x) = \begin{cases}
    1 &\text{if}\;\;\; t = \min \arg \max_{t' \in \mathcal T}\; \left(\hat{m}_n^{t'}(x)\left(1 - \frac{n}{n_{g}} \hat{\gamma}^\star_n(g)\right) - \hat{\mu}_n^{\star,t'}\right)\\
    0 &\text{otherwise}
    \end{cases}
    \end{aligned}
    \;\;\; \forall t \in \mathcal T, \;\forall x \in \mathcal{X},
\end{equation}
where~$\hat{\gamma}_n^{\star}(g) = \sum_{g' \in \mathcal G, g \neq g'} \left(\hat{\lambda}_n^{\star}(g, g') - \hat{\lambda}_n^{\star}(g', g) \right)$.} 
Compared to statistical parity in allocation, we see now that the adjustment is directly applied to the contribution value of the estimated conditional mean treatment outcome function~$\hat{m}_n^t$. The adjustment~$\hat{\gamma}_n^{\star}(g)$, which is unrestricted in sign and inversely weighted by the population proportion of group~$g$, will scale up or down~$\hat{m}_n^t(x)$ for a given treatment~$t$ in order for assignments to achieve parity in outcomes. {\color{blue}If on average group $g$'s outcomes are not as high as other groups' outcomes, we can upweight group $g$'s conditional mean treatment outcomes so that group $g$ individuals are more likely to `pass' the `threshold' value~$\hat{\mu}_n^{\star, t}$. This leads to individuals in group $g$ being more likely to receive treatments under the statistical parity in outcomes policy, and therefore higher group-wise outcomes.} As in the case of statistical parity in allocation in Section~\ref{sec: sp in alloc}, we can similarly define a minority prioritization version of \eqref{eq: policy from dual_sp_out 1} and the respective fairness-constrained sample-based dual-price queuing policies.

\section{Empirical Results} \label{sec: emp results}

In this section, we evaluate the empirical performance of our methodology in two sets of experiments: \emph{(1)}~a synthetic example, where the data generation models for covariates, treatment assignments, and potential outcomes are known; \emph{(2)} an example based on real HMIS and VI-SPDAT data from our central motivating application of allocating scarce housing resources to people experiencing homelessness. In both cases, we split the available data into training and testing/test sets, where the training data is used to learn the dual prices and counterfactual outcomes and the testing data is used to evaluate policy performance.

\subsection{Synthetic Data Experiments} \label{sec: synthetic data experiment}


We use a simple synthetic example to evaluate the performance of our approach in a setting where counterfactuals are known and to study the effects of bias from insufficient training data, large noise terms, and model misspecification on the out-of-sample performance of our methodology. For the synthetic experiments, our primary performance metric is the ratio of the test set outcomes of the sample-based dual-price queuing policy to the outcomes of the perfect foresight test set policy (that knows both the people and resources arriving and the counterfactual outcomes). Thus, a ratio of~$1$ corresponds to the best possible performance. We also study the convergence of the performance of our policy to that of the perfect foresight test set policy as the training sample size is increased. 

\subsubsection{Data Generation.}
\label{sec:synthetic_data_generation}

For the synthetic data experiments, we adapt the data generation process of~\cite{optprestrees} by increasing the number of treatments from two to three to test our methodology in the multiple treatments case. 
In this setting, there are three possible treatments, i.e.,~$\mathcal T = \{0, 1, 2\}$, and each individual is characterized by two covariates, each following a standard normal distribution, i.e.,~$X = [X^1, X^2]$ and~$X^1, X^2 \sim \mathcal{N}(0,1)$. The potential outcomes~$Y^t$,~$t \in \mathcal T$, are expressible as 
\begin{equation*}\label{eq: sim_pot_outcome} 
    \begin{aligned}
    Y^t = \frac{1}{2} X^1 + X^2 + \left( 2 \cdot \mathbbm{1}[t=1] - 1 \right)\frac{1}{4} X^1 + \left(2 \cdot \mathbbm{1}[t=2] - 1 \right) \frac{1}{4} X^2 + \epsilon_t,
    \end{aligned}
\end{equation*}
where~$\epsilon_t \sim \mathcal{N}(0, \sigma)$,~$t \in \mathcal T$, are i.i.d.\ noise terms. In our experiments, we will vary~$\sigma$ in the set~$\{0.1, 0.5, 0.8, 1.15, 1.5 \}$ to study the impact on our method of increasing noise. Indeed, increasing the standard deviation of the noise impacting~$Y^t$,~$t\in \mathcal T$, is expected to affect the quality of the estimates of~$m^t$ and~$\mu^t$ used by our sample-based dual-price queuing policy. From the above, in our notation of conditional mean treatment outcomes, we see that
\begin{equation*}\label{eq: sim_mean_outcome} 
    \begin{aligned}
    m^0(X) = \frac{1}{4} X^1 + \frac{3}{4} X^2, \quad m^1(X) = \frac{3}{4} X^1 + \frac{3}{4} X^2, \quad \text{and} \quad m^2(X) = \frac{1}{4} X^1 + \frac{5}{4} X^2.
    \end{aligned}
\end{equation*}
The historical/training set treatment assignments~$T$ come from a discrete distribution that, in line with our assumptions, depends on the covariates~$X$ only, i.e.,~$p(X) = \left[p^0(X), p^1(X), p^2(X)\right]$, where~$p^t(x)$ is the assignment probability of treatment~$t$ given~$X = x$. In particular, we let treatment assignment probabilities depend only on the treatment that is the best in expectation. This partitions the covariate space into three disjoint regions, each with associated values for~$p(x)$: \emph{(1)} when treatment~0 is the best (or at least as good) in expectation amongst the treatments (i.e.,~$0 = \min \arg \max_{t \in \mathcal T} m^t(x)$), we set~$p(x)=[0.8, 0.1, 0.1]$; \emph{(2)} when treatment~1 is the best (or at least as good) in expectation amongst the treatments (i.e.,~$1 = \min \arg \max_{t \in \mathcal T} m^t(x)$), we set~$p(x)=[0.6, 0.3, 0.1]$; and \emph{(3)} when treatment~2 is the best (or at least as good) in expectation amongst the treatments (i.e.,~$2 = \min \arg \max_{t \in \mathcal T} m^t(x)$), we set~$p(x)=[0.6, 0.1, 0.3]$. This modeling choice mimics realistic scenarios where we anticipate that policy-/decision-makers frequently assign treatments ``correctly.'' In addition, for cases \emph{(2)} and \emph{(3)} where treatments~$1$ and~$2$ are the best in expectation and there is a higher probability of receiving those treatments relative to the probabilities in case \emph{(1)}, treatment~$0$ still has the highest probability at~$60\%$. This modeling choice mimics realistic scenarios where there is a default {no-treatment} scenario (treatment~$0$) and there are limited quantities of treatments~$1$ and~$2$. Finally, we assume that the treatment capacities are known and given by~$b^0=1$,~$b^1=0.1$, and~$b^2=0.05$.


To evaluate the performance of our method on this data, we run 25 simulations. Each of these involves~$360,000$ test set samples (individuals) as a proxy for the asymptotic setting, and a number of training samples that we vary from~$500$ to~$220,000$. 




\subsubsection{Choice of Outcome Estimators.} \label{sec: synthetic_estimators}

Since the asymptotic optimality of our  policy~$\pi$ depends on the uniform almost sure convergence of~$\hat{m}^t$ to~$m^t$ (see Assumption~\ref{ass: estimator of m^t} (ii)), we study the robustness of our method to Assumption~\ref{ass: estimator of m^t} (ii) by using a variety of estimators for learning the mean treatment outcome models. If an incorrect parametric model is employed when learning~$m^t$,~$t \in \mathcal T$, then such estimation errors may cause suboptimal policy learning. From the data generation process, we know that the true~$m^t$ functions are linear, see Section~\ref{sec:synthetic_data_generation}. Thus, we can use decision trees and lasso regression {\color{blue}(with regularization parameter $\alpha = 0.6$ so that lasso is misspecified)} as our misspecified models {\color{blue}(Assumption \ref{ass: estimator of m^t} (ii) does not hold)} and~$k$-NN as a non-parametric model with asymptotic consistency. Finally, we investigate if causal inference methods for ameliorating functional form bias in treatment effects estimation can improve policy learning in the presence of estimation errors of treatment outcomes. In particular, we use inverse propensity score weighted (IPW) and doubly robust (DR) variations of our misspecified models (decision trees and lasso regression) to see if they can mitigate bias from estimation errors. Since IPW and DR require estimating propensity scores, we use decision trees and logistic regression as the model classes for propensity score estimation. For an introduction to these methods, please see~\cite{econml}.

\subsubsection{Performance and Discussion.} \label{subsec: sim results}

\begin{figure*}[!htb]
\begin{center}
\includegraphics[scale=0.6]{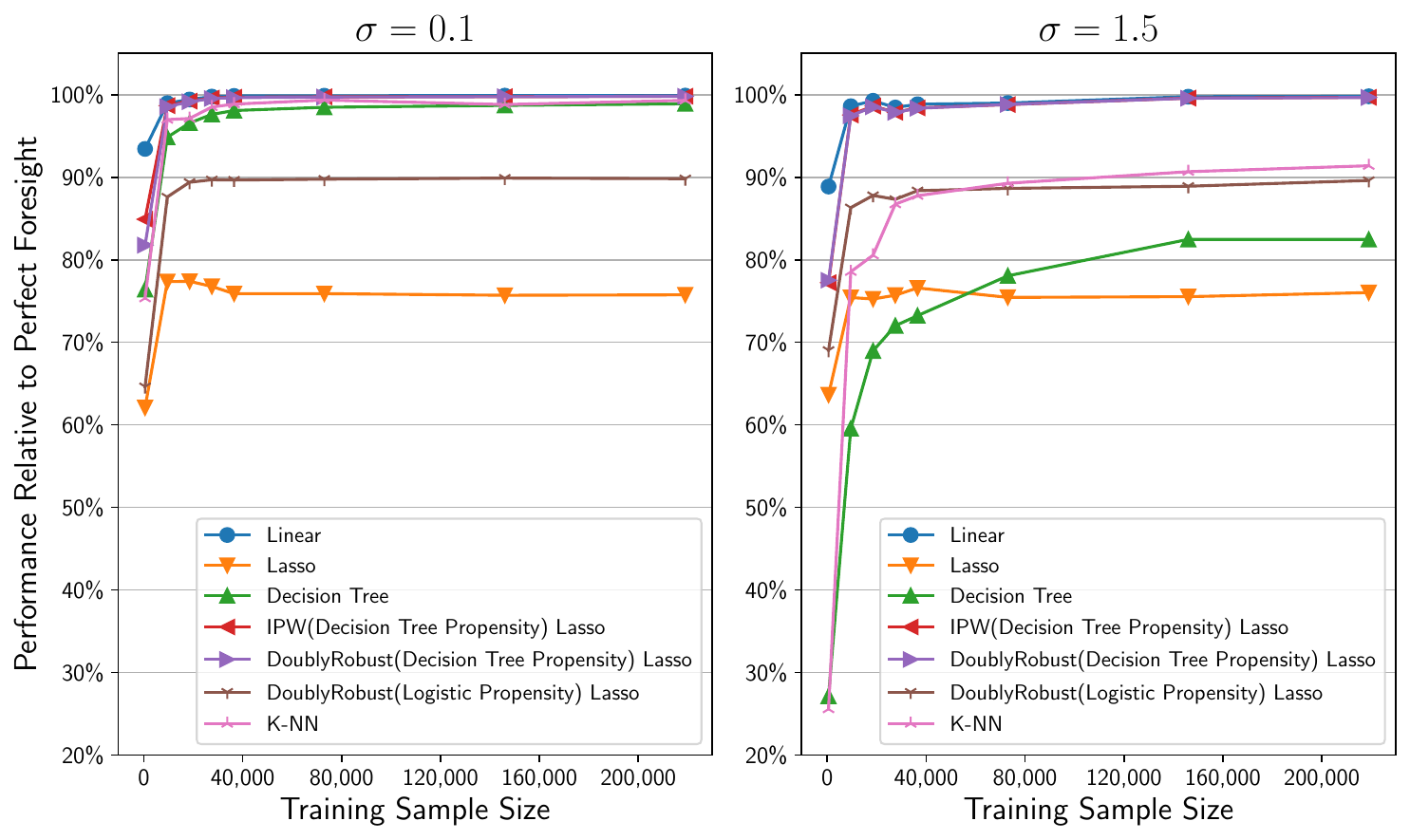}
\caption{Companion figure to the synthetic data results for~$\sigma = 0.1$ and~$\sigma = 1.5$. Both subfigures show the ratio of out-of-sample performance of the sample-based queuing policy to that of the perfect foresight policy in dependence of training set size. Each line corresponds to a policy using a different treatment outcome estimator.}
\label{fig:sim_low_var}
\end{center}
\end{figure*}

Our main results are summarized in Figure~\ref{fig:sim_low_var}, where we plot our main performance metric under policies using various estimators of treatment outcomes. The figure shows our policy performance metric for each estimator in dependence of training sample size, with the left (resp.\ right) figure focused on the low (resp.\ high) variance case where~$\sigma = 0.1$ (resp.~$\sigma = 1.5$). For brevity, in the rest of this section we refer to various dual-price queuing policies using different estimators of the conditional mean treatment outcomes by just the estimator name. For example, the linear, lasso, decision tree policies in Figure~\ref{fig:sim_low_var} are the dual-price queuing policies based on `directly' estimating the conditional mean outcomes using linear, lasso, and decision tree models, where linear is the correct model form while decision trees and lasso are misspecified models. The policies that start with `IPW' and `DoublyRobust' in their names are based on estimators of treatment outcomes that attempt to correct the bias of the lasso based policy by IPW or DR estimates of treatment outcomes. Finally, we have omitted from these results the performance of a `random' policy that randomly assigns treatments as they arrive to waitlisted individuals (which performed near~$0\%$ irrespective of training sample size). Additional results for~$\sigma$ values of~$\{0.5, 0.8, 1.15\}$ can be found in Electronic Companion~\ref{app:simdata}.

From Figure~\ref{fig:sim_low_var} (left), we see that, in the low variance setting, policy performance converges relatively fast (with under~$9,000$ training samples). This is expected: since observed outcomes affect the estimation of counterfactuals, which in turn affect the estimation of policy parameters~$\hat{\mu}^\star$, low noise in the observed outcomes leads to faster convergence of policy performance. As we would expect, the linear estimator policy attains close to perfect performance since it correctly estimates the true conditional mean treatment outcomes. On the other hand, the various modified lasso policies appear to mitigate the estimation bias of the direct lasso estimates and improve policy performance. However, these results also show that the choice of propensity score estimates used to modify the direct treatment outcome estimates can impact policy performance. We see that policies using decision tree estimated propensity scores resulted in better convergence of policy performance (close to the linear policy performance) compared to policies using logistic regression estimated propensity scores. The worse performance of policies using logistic regression estimated propensity scores, compared to that of policies using decision tree estimated propensity scores, may result from the non-linear relation between treatment assignment probabilities and covariates. We emphasize that selecting an accurate propensity score model (even when the historical treatment assignment policy is unknown) is usually feasible in practice with modern machine learning approaches. In this low noise setting, non-parametric estimator based policies such as~$k$-NN can effectively estimate the counterfactuals well enough to achieve close to perfect performance.  
Thus, as expected, in the direct estimator case, our method's performance converges to 1 if the model employed is correct. On the other hand, in policies using IPW or DR adjustments, the choice of propensity score model and estimates impacts the convergence of policy performance.
With a nonparametric estimator, convergence to 1 also occurs.

From Figure~\ref{fig:sim_low_var} (right), we see that in the high variance setting, the policies based on the non-parametric estimators of~$m^t$ all have a deterioration in performance relative to the low variance setting. In particular, the policy using~$k$-NN estimators exhibits slower performance convergence as the high variance noise terms causes the~$k$-NN mean outcome estimation error to decrease at a slower rate.
On the other hand, compared to the low variance setting, the simpler parametric lasso model maintains a similar, though suboptimal, performance. The lasso model is more `robust', in some sense, to more noise as both its mean treatment outcome estimation error and overall policy performance do not deteriorate much from the low to high variance setting. Finally, we see similar behavior in terms of improvement in policy performance from the IPW and DR correction methods to the lasso policy, though the degree of bias correction of IPW and DR methods depend on the estimated propensity scores used. These simple simulation results highlight the importance of model specification for counterfactual estimation and impact on model performance. {\color{blue} We see that when Assumption \ref{ass: estimator of m^t} (ii) does not hold, allocation policies based on misspecified estimators generally do not converge to $100\%$.} We also see that reweighting or doubly robust estimation of treatment outcomes can mitigate potential estimation bias. {\color{blue}Finally, in Electronic Companion~\ref{app:nonlinear experiments}, we investigate if our takeaways so far still hold in the nonlinear setting. We conduct an additional synthetic experiment with a similar setup but where replace the conditional mean treatment outcome functions with nonlinear quadratic functions. We find that the same qualitative takeaways still hold to a large degree.} 

\subsection{Real Data Experiments: Allocating {PSH} and {RRH} Housing in LA} \label{exp: realdata}
We next showcase the performance of our methodology on real data 
to design policies for allocating scarce resources to individuals experiencing homelessness in LA. We obtained our data from LAHSA through a Data Use Agreement and our data comes from the LA County Homeless Management Information System (HMIS) database \citepalias{lahsa_hmis} and administered VI-SPDAT assessments \citepalias{lahsa-vi-spdat}. We obtained Institutional Review Board approval for our data analysis. Additional details on the data sources and preprocessing can be found in Electronic Companion~\ref{app:realdata}.

\subsubsection{HMIS and VI-SPDAT Data Description.} \label{subsubsec: hmis data}

The HMIS data contains deanonymized data of individuals experiencing homelessness and their interactions with the system, which we refer to as \emph{enrollments}. Each enrollment represents an instance in which an individual is assigned a specific resource, where the resource types (introduced in Section~\ref{sec: intro}) are \texttt{PSH}, \texttt{RRH}, and \texttt{SO}. Note that \texttt{SO} is a broad category that includes various services that are not considered permanent housing such as street outreach, homelessness prevention, emergency shelter, supportive services like childcare, and more. We also further breakout \texttt{PSH} into two sub-types \texttt{PSH Tenant-Based} and \texttt{PSH Site-Based} that differ, broadly speaking, in terms of supportive services available on-site and process of actually receiving housing. Further details can be found in in Electronic Companion~\ref{app:realdata}.
Since an individual can have multiple enrollments, we observe for each individual and enrollment pair the relevant date, and assigned resource type. The administered VI-SPDAT surveys, which are used for assessing individual vulnerability and resource prioritization, contain individual covariates, such as disabilities
and prior housing history at the time of assessment.

We include all individuals with an assessment between the time period of 1/12/2015 (roughly when VI-SPDAT started being administered) to 12/31/2019 (to avoid idiosyncratic effects of Covid-19 during 2020) for a total of~$63,764$ samples. We use all individuals assessed between 1/12/2015 to 12/31/2017 as our training set to learn treatment outcomes and construct our policy, and then evaluate on individuals assessed between 1/1/2018 to 12/31/2019 to measure out-of-sample performance. We have a total of 4 possible resources, i.e.,~$\mathcal{T} = \{0,1,2,3\}$, which represent \texttt{no-treatment}, \texttt{RRH}, \texttt{PSH Tenant-Based}, and \texttt{PSH Site-Based}. To solve problem \eqref{eq: sample-based assignment simplified} and derive a sample-based queuing policy, we use the training set capacities for \texttt{no-treatment}, \texttt{RRH}, \texttt{PSH Tenant-Based}, and \texttt{PSH Site-Based} as the estimates of capacity per person, which are~$100\%$,~$14.2\%$,~$4.1\%$, and~$4.6\%$, respectively. 
Since racial equality and equity are important concerns in our motivating example in Section~\ref{sec: intro}, we use race as the protected characteristic for evaluating allocation and outcome fairness. The protected features set is~$\mathcal{G} = \{\texttt{BlackAfAmerican}, \texttt{Hispanic}, \texttt{Other}, \texttt{White}\}$, where \texttt{BlackAfAmerican} stands for Black or African American, and \texttt{Other} includes all other individuals who identify as American Indian or Alaska Native, Asian, Native Hawaiian or Other Pacific Islander, more than one racial group, or answered `Doesn't know'.
In this example, we consider non-White individuals to be minority groups that policymakers may want to prioritize due to historical discrimination. 

\textit{Outcome Definition.} According to the U.S. Department of Housing and Urban Development (HUD), the LA CoC should develop ``action steps to end homelessness and prevent a return to homelessness''~\citepalias{LACoC}. In line with this goal, we will use an individual's return to homelessness after an intervention as the outcome to measure the performance of a policy. This outcome will serve as an input to our optimization model for designing a dual-price policy for allocating scarce housing resources effectively. Specifically, we focus on returns to homelessness within a two-year window following intervention for our real data experiments. However `returns to homelessness' for individuals are not explicitly tracked. Instead, we construct a proxy outcome variable where an individual `returned to homelessness' following intervention if we observe a subsequent enrollment into `emergency shelter', `safe haven', or `street outreach'. We chose these enrollment categories based on discussions with domain experts like matchers and case managers working within the LA CES since these types of interactions with the system indicate a need for homeless services and suggest an individual returned to homelessness again. For individuals without a treatment enrollment into \texttt{RRH} or \texttt{PSH}, we treat their first enrollment into a \texttt{SO} resource as the start of a \texttt{no-treatment} `intervention'. Therefore, we define a positive outcome,~$Y=1$, as \emph{not} observing a subsequent return to homelessness within a two-year window of receiving intervention, and~$Y=0$ otherwise. While other observation window lengths could be chosen, there is a trade-off in terms of window length. Shorter windows are potentially biased proxies since there is not enough time to observe post-intervention outcomes while longer windows result in less individuals for whom we can observe the entire window in our dataset.

 \textit{Policy Evaluation.} We learn three different versions of our proposed policy from the training set and evaluate their performance by their overall proportion of positive outcomes as the `effectiveness' metric and their fairness properties in terms of statistical parity in allocation \eqref{eq: sp alloc} and in outcomes~\eqref{eq: sp_outcome_constraint} on the test set. {\color{blue}We emphasize that our evaluation metric is a counterfactual probability, since the true counterfactuals under assignments different from the historical ones are not available. We explain how we generate these counterfactuals, as well as how we evaluate out-of-sample performance, in Section \ref{sec: counterfactual_realdata}.} To investigate the fairness properties, we look at discrepancies in the percentage receiving housing resources and discrepancies in the proportion of positive outcomes between racial groups in the out-of-sample allocation. In our results, \texttt{Base} refers to the policy without fairness constraints, while the other two refer to fairness-constrained variations defined in Section~\ref{sec: fairness extensions}. Based on our motivating example, we chose to focus on fairness extensions that prioritized minority groups in terms of either allocation, \texttt{Alloc Min Priority}, or outcomes, \texttt{Outcome Min Priority}. To benchmark performance, we compare the three policies to the \texttt{Historical} policy, which uses the same assignments originally found in the data, in terms of the overall proportion of positive outcomes and the two statistical parity metrics. In this way, we can see if our proposed policies can improve on the `effectiveness' and fairness of the currently deployed policy. {\color{blue} Additionally, we include the performance of a \texttt{VI-SPDAT Thresholding} policy, which learns, from historical data, VI-SPDAT score thresholds for allocation so that higher vulnerability scores get more supportive resources like PSH. These thresholds are then used online to assign resources based on an individual's VI-SPDAT score. This offers a comparison point of how the VI-SPDAT score can be used for prioritization.} Finally, we also compare the performance of overall outcomes to a \texttt{Perfect Foresight} policy that implements the sample-based dual-price queuing policy while knowing both the test set distribution of arriving individuals and their counterfactuals. 


\subsubsection{Counterfactuals and Outcome Estimators.} \label{sec: counterfactual_realdata}

Unlike the synthetic experiment, we do not have access to counterfactuals for each individual under treatments other than the one received. Therefore, we take a semi-synthetic approach and use model generated counterfactuals for the test set by choosing models that are well calibrated to the entire data, therefore including data not used by the training estimators. We choose calibration as the measure of fit since we are evaluating based on expected outcomes of a binary variable, and therefore we want our predicted probabilities of a positive outcome to match the observed data. {\color{blue}In particular, to correct for selection bias in our observational data, we use weighted calibration for model validation, where we calculate average positive outcomes using the (self-normalized/Hájek) IPW estimator instead of simple averaging of observed outcomes~\citep{hernan2023causal}.} 

{\color{blue}To fit counterfactual models, we consider the three causal estimation approaches of direct, IPW, and DR that were also used in our synthetic experiment as discussed in~Section~\ref{sec: synthetic_estimators}. As some of these methods may involve estimating both the observed outcomes and propensity scores, our choice of base estimators for the outcomes and propensities include logistic regression (with and without regularization), random forests, gradient boosted trees~\citep{gradientboosting}, and XGBoost trees~\citep{xgboost}.} For model selection, we use standard cross-validation to find the best hyperparameters for each model class that minimize log-loss {\color{blue}for predicting outcomes or the assigned treatment for propensity scores} and then select the final model based on {\color{blue}weighted} calibration curves on a held-out validation set. 

{\color{blue} An additional concern of our model‑generated counterfactuals is potential bias by race because our evaluation is based on these model predicted probabilities. If the outcome models are miscalibrated for a particular racial group (e.g., they systematically over- or under-estimates a group's true positive outcome probabilities under a treatment), then the overall policy outcome numbers we report will have those errors and unfairly favor or penalize a group.}
Though this was not part of our model selection process, we investigated potential racial bias by comparing calibration curves for each racial group across treatments using our selected counterfactual-generating models. For the most part, the calibration curves across racial groups are similar except for the `Other' racial group, likely because this group contained too few samples to get accurate estimates. For additional details, please Section~\ref{app: realdata_counterfactuals_fairness}. 

While the above addresses generating counterfactuals for evaluation, we take a similar approach for actually learning~$\hat{m}^t$ functions but \emph{only} use the training data. For more details on counterfactual generation, model calibration, and racial fairness of the estimates, please see the Electronic Companions~\ref{app: realdata_counterfactuals} and~\ref{app: realdata_counterfactuals_fairness}.

\subsubsection{Experimental Results.}\label{subsec: experimentalresults}

Our efficiency results are summarized in Table~\ref{table:real_overall_outcomes}, which shows the performance of the proposed policies by their out-of-sample expected outcomes {\color{blue}and their relative percentage difference compared to the historical policy (e.g., the relative percentage change of \texttt{Base} versus \texttt{Historical} is calculated as the outcome of \texttt{Base} minus the outcome of \texttt{Historical}, divided by the outcome of \texttt{Historical}). Since the \texttt{VI-SPDAT Thresholding} performance is largely similar to \texttt{Historical}, any comparison to the \texttt{Historical} policy also applies to the \texttt{VI-SPDAT Thresholding} policy. As a point of reference for the effectiveness of different allocations, we also included the potential population expected outcome if every individual receives \texttt{no-treatment}.} 

In terms of the out-of-sample expected outcomes, we see in Table~\ref{table:real_overall_outcomes} that our base model achieves a {\color{blue} $5.16\%$ improvement above the \texttt{Historical} policy, or an absolute~$3.07$ percentage point (p.p.) increase, while having a performance
gap of~$1.55\%$, or~$0.92$ p.p., compared to the \texttt{Perfect Foresight} policy.} We found that the performance gap relative to \texttt{Perfect Foresight} comes from two sources of error: \emph{(i)} estimation error of treatment outcomes~$\hat{m}^t$, which are plug-ins into the assignment policy, and \emph{(ii)} shifts in the distribution of individuals we observe and capacities per person for each treatment from training to testing, which affects the optimal value of~$\mu^\star$ and any estimates~$\hat{\mu}^\star$ of it. These sources of error could come from insufficient training samples for estimating outcomes or observing the distribution of individuals, or from violations of assumptions needed such as stationary distribution from training to testing, or positivity assumption for causal inference. Practically speaking, we see that the \texttt{Historical} policy already performs rather well compared to the \texttt{Perfect Foresight} policy given the scarce amount of resources. {\color{blue} Our policy further improves upon \texttt{Historical}, leading to roughly $500$ more individuals exiting homelessness per year.} This shows our proposed policy can be a \emph{systematic} method to assist and supplement the experience and knowledge of domain experts.  

\begin{table}[h!]

\caption{Companion table to the real data experiment showing out-of-sample proportion of positive outcomes under each policy, and the relative percentage difference compared to the \texttt{Historical} allocation.}
\vspace*{2mm}
\centering
\begin{tabular}{ |s|c|c| }
\hline
 \rowcolor{tablecell} \textbf{Policy} & \textbf{\begin{tabular}{c}Proportion of \\ Positive Outcomes \end{tabular} } & \textbf{\begin{tabular}{c}Percentage Change\\ vs Historical \end{tabular} }  \\ 
\hline \hline
 {\texttt{No Treatment}} &~$54.25\%$ &~$-8.76\%$\\ 
 \hline
 {\texttt{Historical}} &~$59.46\%$ &~$0.00\%$\\ 
 \hline
 {\texttt{VI-SPDAT Thresholding}} &~$59.80\%$ &~$0.01\%$\\ 
 \hline
 {\texttt{Base}} &~$62.53\%$ &~$5.16\%$\\  
 \hline
 {\texttt{Alloc Min Priority}} &~$62.51\%$ &~$5.13\%$\\  
 \hline
 {\texttt{Outcome Min Priority}} &~$62.53\%$ &~$5.16\%$\\  
 \hline
 {\texttt{Perfect Foresight}} &~$63.45\%$ &~$6.71\%$\\  
 \hline 
\end{tabular}

\label{table:real_overall_outcomes}
\end{table}

\begin{figure*}[!htb]
\begin{center}
\includegraphics[scale=0.55]{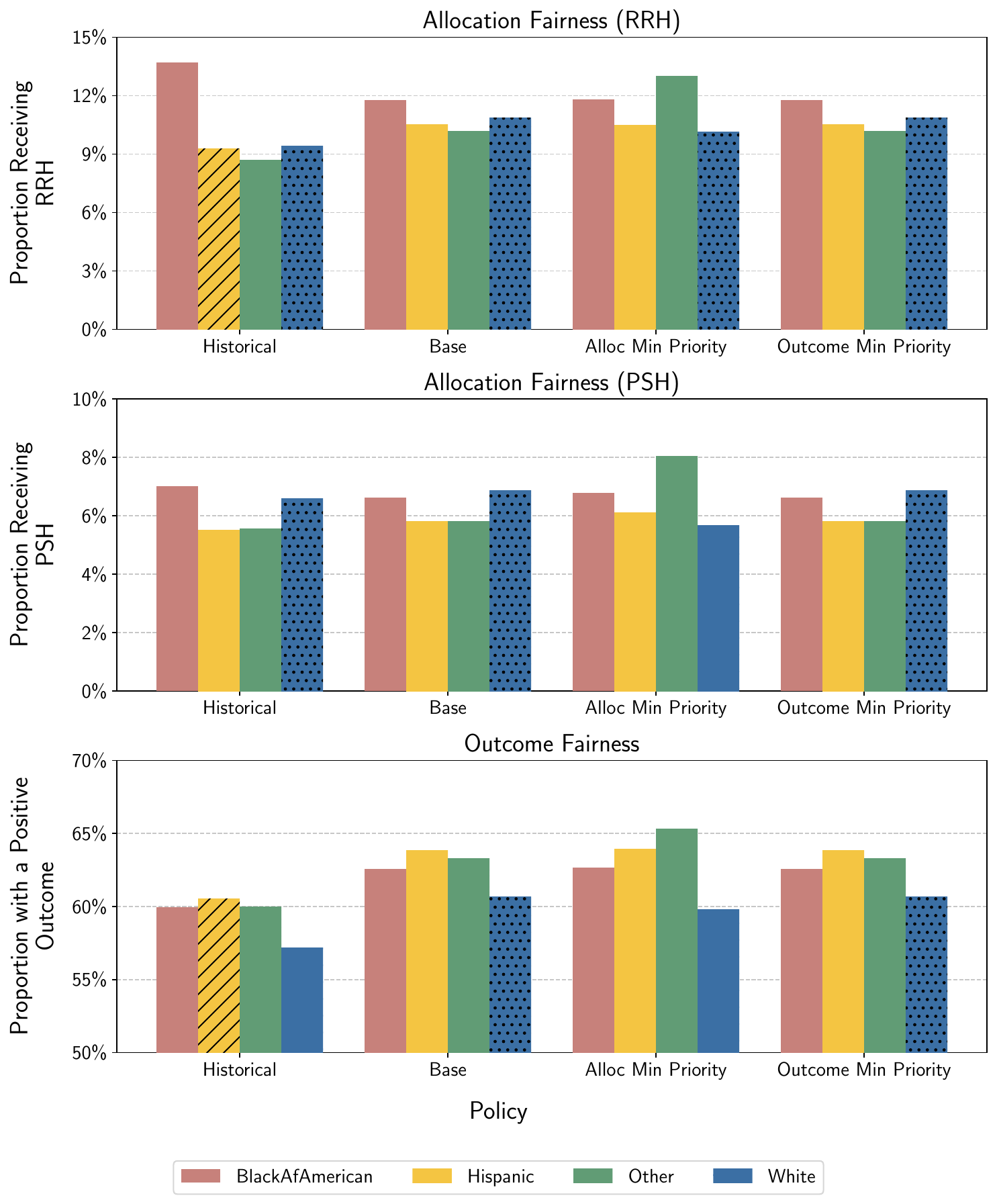}
\caption{Companion figure to the real data experiments: Out-of-sample fairness results showing the proportion of each racial group receiving \texttt{RRH}, \texttt{PSH} (Allocation Fairness) and the proportion of each racial group with a positive outcome (Outcome Fairness).}
\label{fig:real_data_fairness}
\end{center}
\end{figure*}

Figure~\ref{fig:real_data_fairness} illustrates the allocation fairness of each policy by showing the proportion of each racial group receiving  \texttt{RRH} and \texttt{PSH} and the outcome fairness of each policy by showing the expected outcome within each racial group. While \texttt{Base} improves upon the expected outcomes of \texttt{Historical}, {\color{blue}\texttt{Base} has marginally worse allocation fairness properties than \texttt{Historical} when considering the proportion of minority groups receiving \texttt{PSH} and has more uniform allocation rates across all groups for \texttt{RRH}}, as displayed in Figure~\ref{fig:real_data_fairness}. However, by adding appropriate fairness constraints, the \texttt{Alloc Min Priority} policy is able to mitigate any problems and achieve better prioritization of minority groups for allocations. The first subplot of Figure~\ref{fig:real_data_fairness} shows the proportion of each racial group receiving \texttt{RRH} under each policy. We see that the \texttt{Base} and \texttt{Outcome Min Priority} policies {\color{blue}lead to more uniform allocation rates across groups but still have \texttt{RRH}} allocation disparity gaps between \texttt{White} vs \texttt{Hispanic} and \texttt{Other} groups compared to the \texttt{Historical}. {\color{blue}\texttt{Alloc Min Priority}, however, satisfies the minority allocation prioritization fairness constraint presented in Section~\ref{sec: sp in alloc} out-of-sample.} Looking at \texttt{PSH} allocations in the second subplot of Figure~\ref{fig:real_data_fairness}, we see that the \texttt{Historical} policy does not prioritize the \texttt{Hispanic} and \texttt{Other} racial groups to be at least as likely to receive \texttt{PSH} as the majority group. Our \texttt{Base} and \texttt{Outcome Min Priority} policies {\color{blue}maintain similar allocation disparity gaps as \texttt{Historical} while also allocating slightly less now to \texttt{BlackAfAmerican} individuals relative to \texttt{White} individuals.} However, if we explicitly add allocation fairness constraints to our sample-based problem to prioritize minority groups for receiving \texttt{PSH} as shown in Section~\ref{sec: sp in alloc}, then our \texttt{Alloc Min Priority} is able to achieve our desired fairness goal.

In the last subplot of Figure~\ref{fig:real_data_fairness}, we show the proportion of positive outcomes within each racial group for each policy and focus on the prioritization of minority groups to have groupwise outcomes that are at least as good as the majority group. In this case, the \texttt{Historical} policy already achieves this notion of fairness and \texttt{Base} does not create outcome disparities between minority groups compared to the majority group. The \texttt{Outcome Min Priority} policy, which explicitly includes outcome fairness constraints in the sample problem, also achieves the desired fairness goal while also improving \texttt{BlackAfAmerican}, \texttt{Hispanic}, and \texttt{White} individual outcomes relative to \texttt{Historical} (this is also true for \texttt{Base}). We also see in these results that \texttt{Alloc Min Priority} has favorable fairness properties by achieving both fairness notions of prioritization of minority groups for allocation and outcomes, though this need not hold in other problems. In general, allocation statistical parity based fairness notions and constraints are incompatible with outcome statistical parity based fairness~\citep{Jo2022}.  Therefore, depending on the desired fairness goals, policymakers may need to choose between alternative fairness notions, or choose a combination of them to achieve a compromise. 

{\color{blue} From Figure~\ref{fig:real_data_fairness}, we can conclude that fairness-constrained policies such as \texttt{Alloc Min Priority} and \texttt{Outcome Min Priority} can improve the fairness properties of the \texttt{Base} policy and satisfy the fairness constraints specified in Section~\ref{sec: fairness extensions}. We note that although our constraints are imposed at the statistical level, they did hold for the sample path that was realized in our test data.} More interestingly, we see in Table~\ref{table:real_overall_outcomes} that our fairness-constrained policies suffer almost no performance gap in terms of expected outcomes compared to \texttt{Base}, suggesting almost no `price of fairness'. Intuitively, we would expect some loss of performance by enforcing certain constraints on allocations. This may either be due to this application or the flexibility of our methodology to find well-performing policies that satisfy various fairness goals. 

{\color{blue}
Beyond our quantitative takeaways above, we additionally analyze the qualitative behavior of our policy in the Electronic Companion~\ref{subsec: qualitative_behavior}. One potential concern for our policy is that by optimizing for overall expected outcomes, it may overlook vulnerable individuals who
are not responsive to treatments. However, we find in Section~\ref{subsec: alloc_vul} that this concern does not occur in our experimental results because the most vulnerable individuals also tend to benefit the most from treatment and therefore are prioritized by our policy. Additionally, policymakers may want to understand why our policy allocates resources differently from the existing Historical policy. In Section~\ref{subsec: policy_comparison}, we analyze cases where our policy differs from the Historical based on individual covariates. We find that allocation differences were due to covariates that did not come directly from the VI-SPDAT survey but from auxiliary HMIS data sources, such as an individual’s current employment type/status. Please see the Electronic Companion~\ref{subsec: qualitative_behavior} for additional details.
}

{\color{blue}\subsubsection{Queue Wait Times.}\label{sec: policy_queue_waittimes}
Since our dual-price queuing policy involves individuals waiting until they are at the front of the queue and the associated resource becomes available, it is important to understand the wait time properties of our policy and how it affects individuals. We clarify that the waiting time for matching \emph{prior} to taking the VI-SPDAT is an explicit input to our model and here we evaluate the wait time after the assessment. To that end, we additionally evaluate the distribution of individual waiting times within each treatment queue under our \texttt{Base} policy and the impact of increasing wait times on individual outcomes. To do so, we run $250$ simulations of individual and resource arrivals over a $10,000$ day horizon to evaluate the long-run result (our original test period in Section~\ref{subsec: experimentalresults} was $730$ days). We simulate resource arrivals using a Poisson arrival process, where the exponential inter-arrival times are parameterized by the average historical resource arrival rates per day. We similarly simulate individual arrivals using a Poisson arrival process. For each arriving individual, we sample their covariates uniformly with replacement from the historical population. }

{\color{blue}
In Figure~\ref{fig:waittime_figure_maintext}, we show wait time statistics and their progression as more individuals arrive in the system. Within each sample path, for each individual assigned to \texttt{RRH}, \texttt{PSH Site}, or \texttt{PSH Tenant} queues, we compute their wait time as the time difference between their arrival time and the time they are matched with a resource. We compute the mean, $10^{\text{th}}$, $25^{\text{th}}$, $75^{\text{th}}$, and $90^{\text{th}}$ percentile wait times of the $i^{\text{th}}$ individual to arrive in the system across all 250 simulations.

\begin{figure*}[!htb]
\begin{center}
\includegraphics[scale=0.5]{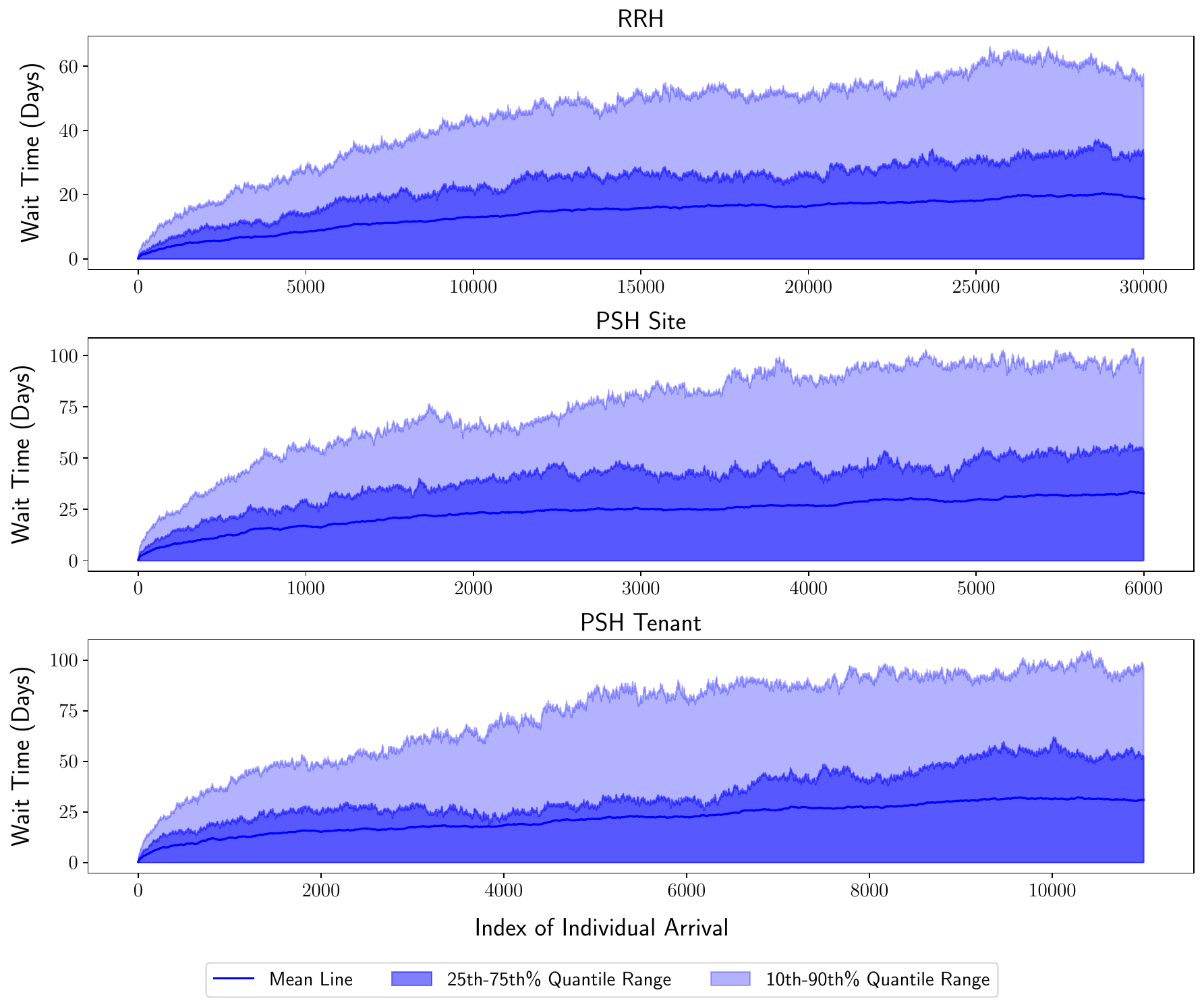}
\caption{Distribution of individual wait times for different treatments across $250$ simulated sample paths.
}
\label{fig:waittime_figure_maintext}
\end{center}
\end{figure*}

We see that the mean waiting time stabilizes at approximately 18, 30, and 30 days for \texttt{RRH}, \texttt{PSH Site-Based}, and \texttt{PSH Tenant-Based}, respectively, while $90^{\text{th}}$ percentile wait times are $55$, $100$, and $100$ days for  \texttt{RRH}, \texttt{PSH Site-Based}, and \texttt{PSH Tenant-Based}, respectively. {\color{blue} These wait times are much shorter than the historical wait times observed in our data, where the mean wait times from $2015$ to $2019$ were $213$, $373$, and $266$ days for individuals assigned to \texttt{RRH}, \texttt{PSH Site-Based}, and \texttt{PSH Tenant-Based}, respectively.} Finally, to understand how reasonable these wait times are, we investigate their effect on outcomes and find that, while the wait times negatively impact potential outcomes, the wait time levels in our numerical experiments result in relatively small decreases in expected conditional outcomes; see Electronic Companion~\ref{app_sub: waittime_causal_impact} for additional details.
For further discussions around wait time and details on the experiment setup and results, please see Electronic Companion~\ref{sec: waittime simulation app}.
}

\section*{Disclaimer}

The views and opinions of the data presented belong to the authors and do not represent the views or opinions of LAHSA.

\ACKNOWLEDGMENT{P.\ Vayanos is funded in part by the National Science Foundation under CAREER award number 2046230. B.~Tang is funded in part by the National Science Foundation Graduate Research Fellowship Program (GRFP). E.\ Rice, P.\ Vayanos, and B.\ Tang are funded in part from the Hilton C.\ Foundation, the Homeless Policy Research Institute, and the Home for Good foundation under the ``C.E.S.\ Triage Tool Research \& Refinement'' grant. The authors are grateful for the support.}




\bibliographystyle{informs2014} 
\bibliography{bibliography.bib} 



%
%
%
%

\ECSwitch


\ECHead{Online Appendix}

{\color{blue}
\section{Dual-Price Queuing Policy under Full Distributional Information} \label{sec: known values}
In this section, we assume that we know the joint distribution of the covariates~$X$ and the potential outcomes~$Y^t$,~$t \in \mathcal T$. We will characterize a dual-price queuing policy under this assumption and show that is optimal in problem \eqref{eq:UpperBoundTrueP}. This characterization and result will allow us to establish the asymptotic optimality of the sample-based dual-price queuing policy in Theorem~\ref{theorem: asymptotic property of sample-based dual-price queuing policy}.

To this end, we consider problem~\eqref{eq:UpperBoundTrueP} and its Lagrangian dual
\begin{equation}\label{eq:DualOfTrueP}
    \begin{aligned}
        \nu^\star
        &= \min_{\mu \in \mathbb R_+^{m+1}} \;\max_{\pi \in \Pi}\; \mathbb E \left[ \sum_{t\in \mathcal T} \pi^t(X) (m^t(X) - \mu^t) \right] + \sum_{t\in \mathcal T} \mu^t  b^t.
    \end{aligned}
\end{equation}
Similarly to the construction in Section \ref{sec: sample-based}, for any~$\mu \in \mathbb R_+^{m+1}$, the inner maximization problem is solved by the policy
\begin{equation}\label{eq: policy from Lagrangian relaxation}
    \begin{aligned}
    \pi^t(x) = \begin{cases}
    1 &\text{if}\quad t = \min \arg \max_{t' \in \mathcal T}\; \left( m^{t'}(x) - \mu^{t'} \right)\\
    0 &\text{otherwise}
    \end{cases}
    \quad\forall t \in \mathcal T, \,\forall x \in \mathcal X.
    \end{aligned}
\end{equation}
By substituting policy \eqref{eq: policy from Lagrangian relaxation} into the dual problem \eqref{eq:DualOfTrueP}, we obtain 
\begin{equation} \label{eq:bestupperboundproblem} \tag{$\mathcal{D}$}
    \begin{aligned}
    \nu^\star &= \min_{\mu \in \mathbb R_+^{m+1}} \; \mathbb E \left[\max_{t \in \mathcal T}\,  (m^t(X) - \mu^{t}) \right] + \sum_{t \in \mathcal T} \mu^t  b^t.
    \end{aligned}
\end{equation}
In the remainder, we denote by~$\mu^\star$ an optimal solution to \eqref{eq:bestupperboundproblem}, which exists because we can restrict the feasible set of \eqref{eq:bestupperboundproblem} to a compact set without loss of generality; see Lemma~\ref{lemma: lemma compact solution set}. By choosing~$\mu$ as~$\mu^\star$ in \eqref{eq: policy from Lagrangian relaxation}, we define a candidate policy~${\pi}^\star$ through
\begin{equation}\label{eq: policy from dual}
    \begin{aligned}
    {\pi}^{\star, t}(x) = \begin{cases}
    1 &\text{if}\quad t = \min \arg \max_{t' \in \mathcal T}\; (m^{t'}(x) - \mu^{\star, t'})\\
    0 &\text{otherwise}
    \end{cases}
    \end{aligned}
    \quad\forall t \in \mathcal T, \,\forall x \in \mathcal X.
\end{equation}
Next, we introduce the queuing implementation of policy~${\pi}^\star$.

\begin{definition}[Dual-price queuing policy]\label{def: dual-price queuing policy}
Establish~$m$ queues each associated with a treatment~$t \in \mathcal T \setminus \{0\}$. At any time, individuals in queue~$t$ are waitlisted for treatment~$t$. There are two events triggering an action: \emph{(i)} arrival of an individual, and \emph{(ii)} arrival of a treatment.
\begin{itemize}
    \item[(i)] When an individual with covariate vector~$x$ arrives, compute the treatment~$t = \min \arg \max_{t' \in \mathcal T}\;{m}^{t'}(x) - \mu^{\star, t'}$, 
    where~$t$ attains the maximum difference between the conditional mean treatment outcome and the associated dual variable. If~$t \neq 0$, assign them to queue~$t$, otherwise assign them {no-treatment}, meaning they will not wait for a particular treatment and will exit the system. If they are the only individual waiting in queue~$t$, check whether there is an available treatment of type~$t$. If there is, assign treatment~$t$ to the individual.
    \item[(ii)] When a treatment of type~$t$ becomes available, assign this treatment to the first individual in queue~$t$. If there is no one in queue~$t$, keep the treatment until an individual is assigned to this queue. 
\end{itemize}
\end{definition}
Similarly to the sample-based dual-price queuing policy, the dual-price queuing policy assigns~$t=0$ to individuals still waiting in queues at the end of the time horizon when implemented in finite horizon. 
We note that in general,~$\mu^\star$, the optimal solution to \eqref{eq:bestupperboundproblem}, may not be unique since the objective function of \eqref{eq:bestupperboundproblem} may not be \emph{strictly} convex; this depends on the form of the joint distribution of~$(X, \{Y^{t}\}_{t \in \mathcal T})$. This means that the candidate policy~$\pi^\star$ defined in \eqref{eq: policy from dual} and the dual-price queuing policy may not be uniquely specified. 

The next theorem shows that the dual-price queuing policy satisfies the capacity constraints in problem~\eqref{eq:UpperBoundTrueP}, and its asymptotic expected average outcome matches the optimal value~$z^\star$ of \eqref{eq:UpperBoundTrueP}. It thus solves problem~\eqref{eq:UpperBoundTrueP}.
\begin{theorem} \label{theorem: asymptotic opt with known values}
Under Assumption~\ref{ass: meanoutcomes}, the dual-price queuing policy is optimal in problem~\eqref{eq:UpperBoundTrueP}.
\end{theorem}
The proof of Theorem~\ref{theorem: asymptotic opt with known values} relies on establishing some structural properties of problem~\eqref{eq:UpperBoundTrueP} and applying standard techniques based on optimality conditions from network revenue management~\citep{talluri1998}, and more generally, from capacity-constrained dynamic programming~\citep{balseiro2023survey}. The idea of the proof is that it is sufficient to show that~$\pi^\star$ is feasible and optimal in~\eqref{eq:UpperBoundTrueP}. Since strong duality holds for problem~\eqref{eq:bestupperboundproblem}, the optimal solutions of~\eqref{eq:bestupperboundproblem} and its dual must satisfy the Karush-Kuhn-Tucker (KKT) conditions. The KKT conditions imply that~$\pi^\star$, which is constructed from an optimal solution~$\mu^\star$ of~\eqref{eq:bestupperboundproblem}, is a feasible solution of~\eqref{eq:UpperBoundTrueP} and yields an objective value equal to~$v^\star = z^\star$. Full proof can be found in the next appendix section.

In contrast to the sample-based dual-price queuing policy, the  dual-price queuing policy in Definition~\ref{def: dual-price queuing policy} is unidentifiable because its definition relies on the joint distribution of~$(X, \{Y^{t}\}_{t \in \mathcal T})$, which is unknown in reality.}

\section{Proofs}
In this section, we provide the proofs of the results in Sections~\ref{sec: known values} and~\ref{sec: sample-based}.

\subsection{Proofs for Section~\ref{sec: known values}}

The proof of Theorem~\ref{theorem: asymptotic opt with known values} requires the following two technical lemmas. 

\begin{lemma}\label{lemma: z r.v. cont dist}
For any~$\mu \in \mathbb R_+^{m+1}$, denote by~$\mu^{-t}$ its subvector~$(\mu^0, \dots, \mu^{t-1}, \mu^{t+1}, \dots, \mu^{m}) \in \mathbb R_+^{m}$ excluding~$\mu^t$, and define the random variable~$Z^{\mu^{-t}} = m^t(X) - \max_{t' \in \mathcal T \setminus \{t\}} \; (m^{t'}(X) -  \mu^{t'})$ for all~$t \in \mathcal{T}$. Under Assumption~\ref{ass: meanoutcomes},~$Z^{\mu^{-t}}$ is continuously distributed for any~$t \in \mathcal{T}$ and~$\mu^{-t} \in \mathbb R_+^{m}$.
\end{lemma}
\proof{Proof.}
We will show that~$\mathbb{P}(Z^{\mu^{-t}} = z) = 0$ for all~$t \in \mathcal{T}$,~$\mu^{-t} \in \mathbb R_+^{m}$, and~$z \in \mathbb{R}$. To this end, fix arbitrary~$t \in \mathcal{T}$,~$\mu^{-t} \in \mathbb R_+^{m}$, and~$z \in \mathbb{R}$. We have~$\mathbb{P}(Z^{\mu^{-t}} = z) \geq 0$ so we only need to show that~$\mathbb{P}(Z^{\mu^{-t}} = z) \leq 0$. We have
\begin{equation*}
\begin{aligned} 
     \mathbb{P}(Z^{\mu^{-t}} = z) &= \mathbb P (m^t(X) - \textstyle \max_{t' \in \mathcal T \setminus \{t\}} \; (m^{t'}(X) -  \mu^{t'}) = z) \\
     &= \mathbb P (\textstyle \min_{t' \in \mathcal T \setminus \{t\}} (m^t(X) - m^{t'}(X) +  \mu^{t'}) = z) \\
     &\leq \mathbb P (\cup_{t' \in \mathcal T \setminus \{t\}} \{m^t(X) - m^{t'}(X) +  \mu^{t'} = z\}) \\
     &\leq \sum_{t' \in \mathcal T \setminus \{t\}} \mathbb P (m^t(X) - m^{t'}(X) +  \mu^{t'} = z) = 0,
    \end{aligned}
\end{equation*}
where 
the first inequality holds because the probability of ~$\min_{t' \in \mathcal T \setminus \{t\}} (m^t(X) - m^{t'}(X) +  \mu^{t'}) = z$ is less than or equal to the probability that~$m^t(X) - m^{t'}(X) +  \mu^{t'} = z$ for at least one~$t' \in \mathcal{T}$, the second inequality follows from the union bound, and the last equality holds because~$m^t(X) - m^{t'}(X)$ is continuously distributed for all~$t$,~$t' \in \mathcal T$ by Assumption~\ref{ass: meanoutcomes}~(ii).
\Halmos
\endproof

In the remainder, we denote by~$\mathbb P_{\operatorname{x}}$ the marginal distribution of~$X$ and by~$\mathbb P_{\operatorname{z}^{\mu^{-t}}}$ the marginal distribution of~$Z^{\mu^{-t}}$ defined in Lemma~\ref{lemma: z r.v. cont dist} when necessary for clarity. The next lemma shows that the objective function of the dual problem \eqref{eq:bestupperboundproblem} is differentiable under Assumption~\ref{ass: meanoutcomes}, and this result will later allow us to establish the optimality conditions for problem~\eqref{eq:bestupperboundproblem}.

\begin{lemma}\label{lemma: differentiability of obj func in upper bound prob}
Define the auxiliary function~$f:\mathbb R^{m+1}\rightarrow \mathbb R$ through
\begin{equation}\label{eq: f}
\begin{aligned}
    f(\mu) = \mathbb E \left[\max_{t' \in \mathcal T}\;\left( m^{t'}(X) - \mu^{t'} \right) \right] = \int_{\mathcal X} \max_{t' \in \mathcal T}\; \left( m^{t'}(x) - \mu^{t'} \right) \,\operatorname{d} \mathbb P_{\operatorname{x}}(x),
    \end{aligned}
\end{equation} 
which appears in the objective function of the dual problem \eqref{eq:bestupperboundproblem}. Under Assumption~\ref{ass: meanoutcomes}, the partial derivative of function~$f$ with respect to~$\mu^t$,~$t \in \mathcal T$, is given by 
\begin{equation*}
    \begin{aligned}
    &\frac{\partial}{\partial \mu^t} f(\mu) = - \mathbb P_{\operatorname{x}} \left(m^t(X) - \mu^t \geq \max_{t' \in \mathcal T \setminus \{t\}} \;  \left( m^{t'}(X) - \mu^{t'} \right) \right).
    \end{aligned}
\end{equation*}
Furthermore, each partial derivative is continuous in~$\mu$, and~$f$ is therefore differentiable. 
\end{lemma}
\proof{Proof.}
We first derive the partial derivative of~$f$ with respect to each~$\mu^t$ and then show that each partial derivative is continuous in~$\mu \in \mathbb{R}_+^{m+1}$. Given a real number~$a$, we use~$a^+$ to denote the positive part function~$\max(a, 0)$. For any~$t \in \mathcal T$, the partial derivative of~$f$ with respect to~$\mu^t$ is given by 
\begin{equation*} 
\begin{aligned} 
     &\frac{\partial}{\partial \mu^t} f(\mu) &&= \frac{\partial}{\partial \mu^t} \int_{\mathcal X} \max_{t' \in \mathcal T}\;\left( m^{t'}(x) - \mu^{t'} \right) \,\text{d} \mathbb P_{\operatorname{x}}(x) \\
     &&&= \frac{\partial}{\partial \mu^t} \int_{\mathcal X} \max \left\{m^{t}(x) - \mu^{t}, \max_{t' \in \mathcal T \setminus \{t\}} \; \left( m^{t'}(x) - \mu^{t'} \right) \right\} \,\text{d} \mathbb P_{\operatorname{x}}(x) \\
     &&&= \frac{\partial}{\partial \mu^t} \int_{\mathcal X} \left(m^{t}(x) - \max_{t' \in \mathcal T \setminus \{t\}} \; \left( m^{t'}(x) - \mu^{t'} \right) - \mu^{t} \right)^+ + \max_{t' \in \mathcal T \setminus \{t\}} \; \left( m^{t'}(x) - \mu^{t'} \right) \,\text{d} \mathbb P_{\operatorname{x}}(x) \\
    &&&= \frac{\partial}{\partial \mu^t} \int_{-\infty}^{\infty} \left(z - \mu^{t} \right)^+  \,\text{d} \mathbb P_{\operatorname{z}^{\mu^{-t}}}(z) \\
    &&&= \frac{\partial}{\partial \mu^t} \int_{\mu^t}^{\infty} \left(z - \mu^{t} \right)  \,\text{d} \mathbb P_{\operatorname{z}^{\mu^{-t}}}(z) \\
    &&&=  \int_{\mu^t}^{\infty} \frac{\partial}{\partial \mu^t} \left(z - \mu^{t} \right)  \,\text{d} \mathbb P_{\operatorname{z}^{\mu^{-t}}}(z),\\
    \end{aligned}
\end{equation*}
where the fourth equality follows from a variable substitution by using the definition of~$Z^{\mu^{-t}}$ from Lemma~\ref{lemma: z r.v. cont dist} and that~$\frac{\partial}{\partial \mu^t} \int_{\mathcal X} \max_{t' \in \mathcal T \setminus \{t\}} \; \left( m^{t'}(x) - \mu^{t'} \right) \,\text{d} \mathbb P_{\text{x}}(x) = 0$, and the last equality follows from the Leibniz integral rule, which applies because for all~$z \in \mathbb{R}$,~$\frac{\partial}{\partial \mu^t} (z-\mu^t) = -1$ exists and is continuous for all~$\mu^t \in \mathbb R_+$, and it is also uniformly bounded by the constant~$1$. We thus have 
\begin{equation}\label{eq:partial der} 
    \begin{aligned}
        \frac{\partial}{\partial \mu^t} f(\mu) &= \int_{\mu^t}^{\infty} \frac{\partial}{\partial \mu^t} \left(z - \mu^{t} \right)  \,\text{d} \mathbb P_{\text{z}^{\mu^{-t}}}(z) \\
        &= - \mathbb P_{\text{z}^{\mu^{-t}}}(Z^{\mu^{-t}} \geq \mu^t) \\
        &= - \mathbb P_{\text{x}} \left(m^t(X) - \mu^t \geq \max_{t' \in \mathcal T \setminus \{t\}} \;  \left( m^{t'}(X) - \mu^{t'} \right) \right),\\
    \end{aligned}
\end{equation}
where the last equality follows again from the definition of~$Z^{\mu^{-t}}$.

Next, we show that~$\frac{\partial}{\partial \mu^t} f(\mu)$ is continuous in~$\mu$, i.e.,~$\lim_{{\mu} \to \mu_0} \frac{\partial}{\partial \mu^t} f({\mu}) = \frac{\partial}{\partial \mu^t} f(\mu_0)$ for every~$\mu_0 \in \mathbb R^{m+1}_+$. To this end, consider an arbitrary~$\mu_0 \in \mathbb R^{m+1}_+$, and note that we have
\begin{equation*}
    \begin{aligned}
    &\lim_{{\mu} \to \mu_0} \frac{\partial}{\partial \mu^t} f({\mu}) &&= \lim_{{\mu} \to \mu_0} - \int_{\mathcal{X}} \mathbbm{1}[h^t(x, \mu) > 0] \,\text{d} \mathbb P_{\text{x}}(x),\\
    \end{aligned}
\end{equation*}
where ~$h^t(x, \mu) = m^{t}(x) - \max_{t' \in \mathcal T \setminus \{t\}} \; \left( m^{t'}(X) - \mu^{t'} \right) - \mu^{t} = Z^{\mu^{-t}} - \mu^{t}$, and the equality follows from \eqref{eq:partial der}. The indicator function~$\mathbbm{1} [h^t(x, \mu) > 0]$ is uniformly bounded by the constant~$1$.
By Lemma~\ref{lemma: z r.v. cont dist} ,~$h^t(x, \mu)$ is continuously distributed for every~$\mu \in \mathbb R^{m+1}_+$ since~$Z^{\mu^{-t}}$ is continuously distributed for any~$\mu^{-t} \in \mathbb R_+^{m}$, and therefore~$\mathbb P_{\text x}(h^t(x, \mu) = 0) = 0$. As the function~$h^t(x, \mu)$ is continuous in~$\mu$ for every~$x \in \mathcal X$, the indicator function~$\mathbbm{1}[h^t(x, \mu) > 0]$ is discontinuous at~$\mu_0$ only if~$h(x, \mu_0) = 0$. As~$\mathbb P_{\text x}(h(X, \mu_0) = 0) = 0$, we have~$\lim_{{\mu} \to \mu_0} \mathbbm{1}[h^t(X, \mu) > 0] = \mathbbm{1}[h(X, \mu_0) > 0]$~$\mathbb{P}_{\text x}$-almost surely. We thus obtain 
\begin{equation*}
    \begin{aligned}
    \lim_{{\mu} \to \mu_0} \frac{\partial}{\partial \mu^t} f({\mu}) &= \lim_{{\mu} \to \mu_0} - \int_{\mathcal{X}} \mathbbm{1}[h^t(x, \mu) > 0] \,\text{d} \mathbb P_{\text{x}}(x) \\
    &=  - \int_{\mathcal{X}} \lim_{{\mu} \to \mu_0} \mathbbm{1}[h^t(x, \mu) > 0] \,\text{d} \mathbb P_{\text{x}}(x) \\
    &= - \int_{\mathcal{X}} \mathbbm{1}[h(x, \mu_0) > 0] \,\text{d} \mathbb P_{\text{x}}(x) \\
    &= \frac{\partial}{\partial \mu^t} f({\mu_0}),
    \end{aligned}
\end{equation*}
where the second equality follows from the Dominated Convergence Theorem. This proves that the partial derivative~$\frac{\partial}{\partial \mu^t} f({\mu})$ is continuous in~$\mu$. Function~$f$ is thus differentiable. \Halmos
\endproof

Lemma~\ref{lemma: differentiability of obj func in upper bound prob} implies that the objective function of problem~\eqref{eq:bestupperboundproblem} is differentiable. As problem~\eqref{eq:bestupperboundproblem} is a convex optimization problem and satisfies Slater’s condition, strong duality holds, i.e., the optimal value of \eqref{eq:bestupperboundproblem} coincides with the one of its dual problem. Denote by~$\lambda$ the dual variable of the constraint~$\mu \geq 0$ in problem~\eqref{eq:bestupperboundproblem}. If~$\mu^\star$ and~$\lambda^\star$ are solutions to problem \eqref{eq:bestupperboundproblem} and its dual, respectively, then~$\mu^\star, \lambda^\star$ satisfy the Karush–Kuhn–Tucker (KKT) conditions 
\begin{equation} \label{eq: KKT conditions}
    \begin{aligned}
     \mathbb P \left(m^t(X) - \mu^{\star, t} \geq \max_{t' \in \mathcal T \setminus \{t\}} \;  (m^{t'}(X) - \mu^{\star, t'}) \right) + \lambda^{\star, t} &= b^t  &\text{(Stationarity Condition)},\\
     \lambda^{\star, t} \mu^{\star, t} &= 0 \quad &\text{(Complementary Slackness)},\\
     \mu^{\star, t}\geq 0, \quad \lambda^{\star, t} &\geq 0 \quad &\text{(Primal~$\&$ Dual Feasibility)},
    \end{aligned}
\end{equation}
for all~$t \in \mathcal T$.
Intuitively, the KKT conditions imply that policy~${\pi}^\star$ defined in \eqref{eq: policy from dual} satisfies the capacity constraints in expectation. 
To see this, note that~$\mathbb P \left(m^t(X) - \mu^{\star, t} \geq \max_{t' \in \mathcal T \setminus \{t\}} \;  (m^{t'}(X) - \mu^{\star, t'}) \right)$
is the expected number of individuals assigned treatment~$t$ under policy~${\pi}^\star$. 
As~$\lambda^{\star, t} \geq 0$ by dual feasibility, this quantity cannot exceed the available capacity~$b_t$ by the stationarity condition. 

We are now ready to prove Theorem~\ref{theorem: asymptotic opt with known values}.
\proof{Proof of Theorem~\ref{theorem: asymptotic opt with known values}.} 
For any~$\mu^\star$ optimal in \eqref{eq:bestupperboundproblem}, the treatment assignment probabilities of~$\pi^\star$ defined in \eqref{eq: policy from dual} and that of the dual-price queuing policy defined in Definition~\ref{def: dual-price queuing policy} coincide. Indeed, the dual-price queuing policy is nothing more than an online implementation of~$\pi^\star$ since both make assignments using the same treatment assignment criteria. The difference is that the dual-price queuing policy waitlists individuals for treatments by assigning them to queues since there may be a mismatch in timing between the arrival of individuals and resources.
It is thus sufficient to show that~$\pi^\star$ is optimal in \eqref{eq:UpperBoundTrueP}. We first prove that~$\pi^\star$ is feasible in \eqref{eq:UpperBoundTrueP} and then show that its expected average outcome matches the optimal value~$z^\star$ of \eqref{eq:UpperBoundTrueP}. 

Policy~$\pi^\star$ is feasible in \eqref{eq:UpperBoundTrueP} because
\begin{equation*}
    \begin{aligned}
     \mathbb E [{\pi}^{\star, t} (X)] &= \mathbb P\left( t = \min \operatorname{argmax}_{t' \in \mathcal{T}} m^{t'}(X) - \mu^{\star, t'} \right) \\
     &\leq \mathbb P \left(m^t(X) - \mu^{\star, t} \geq \textstyle \max_{t' \in \mathcal T \setminus \{t\}} \;  (m^{t'}(X) - \mu^{\star, t'}) \right) \leq b^t \quad \forall t \in \mathcal T,
    \end{aligned}
\end{equation*}
where the first inequality follows from the fact that if the event~$t = \min \operatorname{argmax}_{t' \in \mathcal{T}} m^{t'}(X) - \mu^{\star, t'}$ occurs, then the inequality~$m^t(X) - \mu^{\star, t} \geq \textstyle \max_{t' \in \mathcal T \setminus \{t\}} \;  (m^{t'}(X) - \mu^{\star, t'})$ must hold, and the second inequality follows from the stationarity and dual feasibility conditions in \eqref{eq: KKT conditions}.

As~$\pi^\star$ is feasible in \eqref{eq:UpperBoundTrueP}, its expected average outcome cannot exceed~$z^\star$. We next show that the expected average outcome of~$\pi^\star$ is at least as high as~$z^\star$, which implies that the two values coincide. The expected average outcome of~$\pi^\star$ amounts to 
\begin{equation}\label{eq: lowerboundbyalternatesystem}
    \begin{aligned}
    &\mathbb E\left[ \sum_{t \in \mathcal T} \mathbbm{1}\left[ t = \min \operatorname{argmax}_{t' \in \mathcal{T}} m^{t'}(X) - \mu^{\star, t'} \right] m^t (X) \right] \\ 
    &= \mathbb E\left[ \max_{t \in \mathcal{T}}\; (m^t(X) - \mu^{\star, t}) + \sum_{t \in \mathcal T} \mathbbm{1} \left[ t = \min \operatorname{argmax}_{t' \in \mathcal T} {m}^{t'}(X) - \mu^{\star, t'} \right] \mu^{\star, t} \right] \\
     &= \mathbb E\left[\max_{t \in \mathcal{T}}\; (m^t(X) - \mu^{\star, t}) \right] + \sum_{t \in \mathcal T}  \mathbb P\left( t = \min \operatorname{argmax}_{t' \in \mathcal{T}} m^{t'}(X) - \mu^{\star, t'} \right) \mu^{\star, t} \\
     &= \mathbb E\left[\max_{t \in \mathcal{T}}\; (m^t(X) - \mu^{\star, t}) \right] + \sum_{t \in \mathcal T}  b^t \mu^{\star, t} \\
     &= \nu^\star  \geq z^\star.
    \end{aligned}
\end{equation}
The third equality above follows from the KKT conditions~\eqref{eq: KKT conditions}, which imply that 
%
\begin{equation*}
    \begin{aligned}
    \mu^{\star, t} b^t &= \mu^{\star, t}  \left [ \mathbb P \left(m^t(X) - \mu^{\star, t} \geq \max_{t' \in \mathcal T \setminus \{t\}} \;  (m^{t'}(X) - \mu^{\star, t'}) \right) + \lambda^{\star, t} \right]\\
    &= \mu^{\star, t}\, \mathbb P \left(m^t(X) - \mu^{\star, t}\geq \max_{t' \in \mathcal T \setminus \{t\}} \;  (m^{t'}(X) - \mu^{\star, t'}) \right)\\
    &= \mu^{\star, t}\,  \mathbb P \left( t = \min \operatorname{argmax}_{t' \in \mathcal T} m^{t'}(X) - \mu^{\star, t'} \right),
    \end{aligned}
\end{equation*}
where the first and second equalities follows from the stationarity and complementary slackness conditions, respectively, and the third equality holds because~$\operatorname{argmax}_{t' \in \mathcal{T}} m^{t'}(X) - \mu^{\star, t}$ is a singleton almost surely.
The last equality in \eqref{eq: lowerboundbyalternatesystem} holds because~$\mu^\star$ is optimal in \eqref{eq:bestupperboundproblem}, and the inequality follows from weak duality. The claim thus follows.
\Halmos
\endproof

\subsection{Proofs for Section~\ref{sec: sample-based}} \label{app: sample-based}

 The proof of Theorem~\ref{theorem: asymptotic property of sample-based dual-price queuing policy} requires the following technical lemmas: Lemmas~\ref{lemma: lemma compact solution set} --~\ref{lemma: opt mu close}. 

We first define the following functions~$\nu:\mathbb R^{m+1}\rightarrow \mathbb R$ and~$\hat{\nu}_{n}:\mathbb R^{m+1}\rightarrow \mathbb R$ for denoting the objective functions of the dual problem~\eqref{eq:bestupperboundproblem} and its sample approximation \eqref{eq: sample-based assignment simplified} as functions of~$\mu$. These definitions will be used throughout Section~\ref{app: sample-based}.
\begin{equation} \label{eq: J quantities}
    \begin{aligned}
    & \nu(\mu) = \mathbb E \left[\max_{t \in \mathcal T}  (m^{t}(X) - \mu^{t}) \right] + \sum_{t \in \mathcal T} \mu^t  b^t \\ 
    & \hat{\nu}_{n}(\mu) =  \frac{1}{n} \sum_{i=1}^n \max_{t \in \mathcal{T}}\;\left( \hat{m}_n^t(x_i) -  \mu^t \right) + \sum_{t\in \mathcal T}  \mu^t  b^t.
    \end{aligned}
\end{equation} 
By definition, we have~$\nu^\star = \min_{\mu \in \mathbb{R}_+^{m+1}} \nu(\mu)$ and~$\hat{\nu}_{n}^\star = \min_{\mu \in \mathbb{R}_+^{m+1}} \hat{\nu}_{n}(\mu)$. Denote by~$\mathcal{S}^\star$ and~$\hat{\mathcal{S}}_{n}^\star$ the sets of optimal solutions to problems \eqref{eq:bestupperboundproblem} and \eqref{eq: sample-based assignment simplified}, respectively. The next lemma shows that~$\mathcal{S}^\star$ is contained in a compact set as long as the function~$m^t$ is uniformly bounded for all~$t \in \mathcal T$, which is guaranteed under Assumption~\ref{ass: meanoutcomes}. 
\begin{lemma} \label{lemma: lemma compact solution set} 
Under Assumption~\ref{ass: meanoutcomes}, the set~$\mathcal{S}^\star$ of optimal solutions to problem \eqref{eq:bestupperboundproblem} is non-empty and contained in the compact set~$[0, C]^{m+1} \subset \mathbb{R}_+^{m+1}$, where~$C$ is defined as in Assumption~\ref{ass: meanoutcomes}(i), i.e.,~$\vert m^t(x) \vert \leq C$ for all~$x \in \mathcal X$ and~$t \in \mathcal T$. 
\end{lemma}
\proof{Proof.}
We will first show that the optimal solution set~$\mathcal{S}^\star$ is contained in the compact set~$[0, C]^{m+1}$ and then show that the set~$\mathcal{S}^\star$ is non-empty.

Suppose for contradiction that~$\mathcal{S}^\star \nsubseteq [0, C]^{m+1}$, i.e., there exists a~$\mu \in \mathcal{S}^\star$ such that~$\mu \not\in [0, C]^{m+1}$. As~$\mu \not\in [0, C]^{m+1}$, there exists at least one~$t \in \mathcal T$ such that~$\mu^t > C$. Denote by~$\mathcal T' = \{t \in \mathcal T \,:\, \mu^t > C\}$ the nonempty set of treatments~$t$ for which~$\mu^t > C$. We construct another feasible solution~$\mu' \in [0, C]^{m+1}$ via~$(\mu')^t = \mu^t$ for all~$t \in \mathcal T \setminus \mathcal T'$ (all~$t$'s where~$\mu^t$ is in the interval~$[0, C]$), and~$(\mu')^t = C$ for all~$t \in \mathcal T'$ (all~$t$'s where~$\mu^t$ is not in the interval~$[0, C]$). We will show that the objective value of~$\mu'$ in problem \eqref{eq:bestupperboundproblem} is strictly lower than that of~$\mu$. This implies that~$\mu$ cannot be optimal and results in a contradiction.

We will first show that~$\max_{t \in \mathcal T}  (m^{t}(x) - \mu^{t}) = \max_{t \in \mathcal T}  (m^{t}(x) - (\mu')^{t})$ for all~$x \in \mathcal X$. To this end, recall that treatment~$0$ represents the no-treatment option, and~$b^0=1$, which makes the capacity constraint corresponding to this treatment redundant. This implies that we can assume without loss of generality~$\mu^0 = 0$, and therefore~$0 \not \in \mathcal T'$. For any~$t \in \mathcal T'$ and~$x \in \mathcal X$, we have
\begin{equation}  \label{eq: mu upperbound}
    \begin{aligned}
    m^t(x) - \mu^{t} < m^t(x) - C = m^t(x) - (\mu')^{t} \leq m^0(x) = m^0(x) - \mu^{0} = m^0(x) - (\mu')^{0},
    \end{aligned}
\end{equation}
where the first inequality follows from the definition of~$\mathcal T'$, the first equality follows from the definition of~$\mu'$, and the second inequality holds because~$m^t(x) \leq C$ and~$m^0(x)$ is non-negative. The last two equalities hold because~$\mu^0 = (\mu')^0  = 0$. By equation \eqref{eq: mu upperbound}, we have~$m^t(x) - \mu^{t} < m^0(x) - \mu^{0}$ and~$m^t(x) - (\mu')^{t} \leq m^0(x) - (\mu')^{0}$ for all~$t \in \mathcal T'$. This implies that for any~$x \in \mathcal{X}$, we have
\begin{equation*}
    \begin{aligned}
        \max_{t \in \mathcal T}  (m^{t}(x) - \mu^{t}) = \max_{t \in \mathcal T \setminus \mathcal T'}  (m^{t}(x) - \mu^{t}) = \max_{t \in \mathcal T \setminus \mathcal T'}  (m^{t}(x) - (\mu')^{t}) = \max_{t \in \mathcal T}  (m^{t}(x) - (\mu')^{t}),
    \end{aligned}
\end{equation*}
where the first and third equality hold because no~$t \in \mathcal T'$ can attain the maximum in~$\max_{t \in \mathcal T}  (m^{t}(x) - \mu^{t})$ and~$\max_{t \in \mathcal T}  (m^{t}(x) - (\mu')^{t})$ by \eqref{eq: mu upperbound}, and
the second equality holds because~$(\mu')^t = \mu^t$ for all~$t \in \mathcal T \setminus \mathcal T'$.

The objective value~$\nu(\mu)$ of~$\mu$ thus exceeds the objective value~$\nu(\mu')$ of~$\mu'$ as
\begin{equation*}
    \begin{aligned}
        \nu(\mu) = \mathbb E \left[\max_{t \in \mathcal T}  (m^{t}(X) - \mu^{t}) \right] + \sum_{t \in \mathcal T} \mu^t  b^t &= \mathbb E \left[\max_{t \in \mathcal T}  (m^{t}(X) - (\mu')^{t}) \right] + \sum_{t \in \mathcal T} \mu^t  b^t \\
        &> \mathbb E \left[\max_{t \in \mathcal T}  (m^{t}(X) - (\mu')^{t}) \right] + \sum_{t \in \mathcal T} (\mu')^t  b^t \\
        &= \nu(\mu'),
    \end{aligned}
\end{equation*}
where the second equality holds because we showed that~$\max_{t \in \mathcal T}  (m^{t}(x) - \mu^{t}) = \max_{t \in \mathcal T}  (m^{t}(x) - (\mu')^{t})$ for all~$x \in \mathcal X$, and the inequality holds because~$(\mu')^t \leq \mu^t$ for all~$t \in \mathcal{T}$ and the inequality is strict for all~$t \in \mathcal{T}'$ by definition of~$\mu'$. This implies that~$\mu$ cannot be optimal and results in a contradiction. Thus, we must have~$\mathcal{S}^\star \subseteq [0, C]^{m+1}$.

We now show that~$\mathcal{S}^\star$ is always non-empty. As~$\mathcal{S}^\star \subseteq [0, C]^{m+1}$, the dual problem~\eqref{eq:bestupperboundproblem} is equivalent to~$\min_{\mu \in [0, C]^{m+1}} \nu(\mu)$, where we only minimize over the compact set~$[0, C]^{m+1}$. Since~$\nu(\mu)$ is a convex function, and hence continuous, in~$\mathbb{R}^n$, it attains a minimum over any compact set (Proposition A.2.7 of~\cite{bertsekas2009convex}), which implies~$\mathcal{S}^\star$ is non-empty.
\Halmos
\endproof

Similar to Lemma~\ref{lemma: lemma compact solution set}, the next lemma shows that the optimal solution set~$\hat{\mathcal{S}}_{n}^\star$ must also lie within a compact set.
\begin{lemma} \label{lemma: sample compact solution set alt}
Under Assumption~\ref{ass: estimator of m^t}, for any~$n \in \mathbb{N}$, the set~$\hat{\mathcal{S}}_{n}^\star$ is non-empty and contained in the compact set~$[0, \hat{C}]^{m+1} \subset \mathbb{R}_+^{m+1}$, where~$\hat{C}$ is defined as in Assumption~\ref{ass: estimator of m^t}(i), i.e.,~$|\hat{m}^t_n(x)| \leq \hat{C}$ for all~$x \in \mathcal{X}$,~$t \in \mathcal{T}$, and~$n \in \mathbb{N}$. 
\end{lemma}

\proof{Proof.}
Fix an arbitrary~$\omega \in \Omega$,~$i.e$, fix historical samples. By Assumption~\ref{ass: estimator of m^t}(i),~$\hat{C}$ is a uniform upper bound on~$|\hat{m}^t_n(x)|$ for all ~$x \in \mathcal{X}$,~$t \in \mathcal{T}$, and~$n \in \mathbb{N}$. 
Following the same arguments used in the proof of Lemma~\ref{lemma: lemma compact solution set}, we can show that any solution~$\mu \not\in [0, \hat{C}]^{m+1}$ will be sub-optimal and there exists a solution in~$[0, \hat{C}]^{m+1}$ with strictly lower objective value. Therefore, the feasible set of problem~\eqref{eq: sample-based assignment simplified} can be restricted to the compact set~$[0, \hat{C}]^{m+1}$. For all~$n \in \mathbb{N}$,~$\hat{\nu}_{n}(\mu)$ is a convex function and hence attains a minimum over any compact set. Since the choice of~$\omega$ was arbitrary, for any~$\omega \in \Omega$,~$\hat{\mathcal{S}}_{n}^\star$ is non-empty for all~$n \in \mathbb{N}$ and contained in~$[0, \hat{C}]^{m+1}$. \Halmos
\endproof

The next lemma shows that the objective function~$\hat{\nu}_{n}(\mu)$ of the sample-approximate dual problem~\eqref{eq: sample-based assignment simplified} converges almost surely to the objective function~${\nu}(\mu)$ of the true dual problem \eqref{eq:bestupperboundproblem} uniformly in~$\mu$.
\begin{lemma} \label{lemma: lemma unif converg} 
Under Assumptions~\ref{ass: meanoutcomes} and~\ref{ass: estimator of m^t}, the objective function~$\hat{\nu}_{n}(\mu)$ of problem \eqref{eq: sample-based assignment simplified} converges almost surely to the objective function~${\nu}(\mu)$ of problem \eqref{eq:bestupperboundproblem} uniformly on any nonempty compact set~$\mathcal{S} \subseteq \mathbb R_+^{m+1}$ as~$n$ tends to infinity, that is, 
%
\begin{equation*}
    \begin{aligned}
    \lim_{n \rightarrow \infty} \sup_{\mu \in \mathcal{S}} \left| \nu(\mu) - \hat{\nu}_{n}(\mu) \right| = 0.
    \end{aligned}
\end{equation*}
\end{lemma}
\proof{Proof.}
Consider an arbitrary nonempty compact set~$\mathcal S \subseteq \mathbb R_+^{m+1}$. By triangle inequality, we have 
\begin{equation*}
    \begin{aligned}
    \sup_{\mu \in \mathcal{S}} \left| \nu(\mu) - \hat{\nu}_{n}(\mu) \right| &\leq  \sup_{\mu \in \mathcal{S}} \biggl( \left| \nu(\mu) - \nu'_{n}(\mu) \right| + \left| \nu'_{n}(\mu) - \hat{\nu}_{n}(\mu) \right| \biggr)\\
    &\leq  \sup_{\mu \in \mathcal{S}} \left| \nu(\mu) - \nu'_{n}(\mu) \right| + \sup_{\mu \in \mathcal{S}} \left| \nu'_{n}(\mu) - \hat{\nu}_{n}(\mu) \right| ,
    \end{aligned}
\end{equation*}
where~$\nu'_{n}$ is an sample average approximation of~$\nu$ using the true~$m^t$ functions, i.e.,
\begin{equation*}
    \begin{aligned}
    & \nu'_{n}(\mu) &&= \frac{1}{n} \sum_{i=1}^n \max_{t \in \mathcal{T}}\;\left( m^t(x_i) -  \mu^t \right) + \sum_{t \in \mathcal T} \mu^t  b^t.
    \end{aligned}
\end{equation*} 
In the following, we investigate the terms~$\sup_{\mu \in \mathcal{S}} \left| \nu(\mu) - \nu'_{n}(\mu) \right|$ and~$\sup_{\mu \in \mathcal{S}} \left| \nu'_{n}(\mu) - \hat{\nu}_{n}(\mu) \right|$ one by one.

We first consider the term~$\sup_{\mu \in \mathcal{S}} \left| \nu(\mu) - \nu'_{n}(\mu) \right|$ and show that 
\begin{equation} \label{eq: unif converg term 1}
    \begin{aligned}
    \lim_{n \rightarrow \infty} \sup_{\mu \in \mathcal{S}}\left| \nu(\mu) - \nu'_{n}(\mu) \right| = 0.
    \end{aligned}
\end{equation}
To this end,
let~$F(x, \mu) = \max_{t \in \mathcal{T}}\;\left( m^t(x) -  \mu^t \right)$, and note that, for any~$x \in \mathcal X$,~$F(x, \mu)$ is continuous in~$\mu$ as it is a piecewise-linear function of~$\mu$. For any~$\mu \in \mathcal{S}$, the function~$F(x, \mu)$ is dominated by the integrable function~$g(x) = C$
as~$C$ is a uniform bound on~$m^t$ for each~$t \in \mathcal{T}$ by Assumption~\ref{ass: meanoutcomes}(i). 
We can thus invoke Theorem~7.53 of~\cite{shapiro} and conclude that~$\frac{1}{n} \sum_{i=1}^n F(X_i, \mu)$ converges uniformly on~$\mathcal S$ and almost surely to~$\mathbb E \left[ F(X, \mu) \right]$ as~$n$ tends to infinity. This implies that ~$\nu'_{n}(\mu)$ converges uniformly on~$\mathcal S$ and almost surely to~$\nu(\mu)$. We thus proved \eqref{eq: unif converg term 1}.

Next, we consider the term~$\sup_{\mu \in \mathcal{S}} \left| \nu'_{n}(\mu) - \hat{\nu}_{n}(\mu) \right|$ and similarly show that 
\begin{equation}\label{eq: unif converg term 2}
    \begin{aligned}
    \lim_{n \rightarrow \infty} \sup_{\mu \in \mathcal{S}} \left| \nu'_{n}(\mu) - \hat{\nu}_{n}(\mu) \right| = 0.
    \end{aligned}
\end{equation}
For any~$\mu \in \mathcal{S}$, we have
\begin{equation}\label{eq: unif converg 2}
    \begin{aligned}
    \left| \nu'_{n}(\mu) - \hat{\nu}_{n}(\mu) \right| &= \left| \frac{1}{n} \sum_{i=1}^n \max_{t \in \mathcal{T}}\;\left( m^t(x_i) -  \mu^t \right) - \max_{t \in \mathcal{T}}\;\left( \hat{m}_n^t(x_i) -  \mu^t \right) \right|\\ 
    & \leq \frac{1}{n} \sum_{i=1}^n \left|  \max_{t \in \mathcal{T}}\;\left( m^t(x_i) -  \mu^t \right) - \max_{t \in \mathcal{T}}\;\left( \hat{m}_n^t(x_i) -  \mu^t \right) \right| \\
    & \leq \frac{1}{n} \sum_{i=1}^n \left|  \max_{t \in \mathcal{T}}\; \left( m^t(x_i) -  \hat{m}_n^t(x_i) \right) \right| \\
    & \leq \frac{1}{n} \sum_{i=1}^n \max_{t \in \mathcal{T}}\; \left|   m^t(x_i) - \hat{m}_n^t(x_i)   \right| \\
    & \leq \sup_{x \in \mathcal{X}}\; \max_{t \in \mathcal{T}}\; \left|   m^t(x) - \hat{m}_n^t(x)   \right|.
    \end{aligned}
\end{equation} 
The first inequality follows by triangle inequality. To see that the second inequality holds, fix an arbitrary~$x \in \mathcal X$ and~$\mu \in \mathcal{S}$ and let~$t^\star = \min \arg \max_{t \in \mathcal T}\; (m^{t}(x) - \mu^t)$. We then have
\begin{equation*}
    \begin{aligned}
    \max_{t \in \mathcal{T}}\;\left( m^t(x) -  \mu^t \right) - \max_{t \in \mathcal{T}}\;\left( \hat{m}_n^t(x) -  \mu^t \right) &\leq \left( m^{t^\star}(x) -  \mu^{t^\star} \right) - \left( \hat{m}_n^{t^\star}(x) -  \mu^{t^\star} \right) \\
    &= m^{t^\star}(x) - \hat{m}_n^{t^\star}(x) \\
    & \leq \max_{t \in \mathcal{T}}\; \left[ m^t(x) -  \hat{m}_n^t(x) \right],
    \end{aligned}
\end{equation*} 
where the first inequality holds by definition of~$t^\star$ and~$\hat{m}_n^{t^\star}(x) -  \mu^{t^\star} \leq \max_{t \in \mathcal{T}}\;\left( \hat{m}_n^t(x) -  \mu^t \right)$. Since the choice of~$x$ and~$\mu$ were arbitrary, the above holds for all~$x \in \mathcal X$ and~$\mu \in \mathcal{S}$, thereby proving the second inequality in \eqref{eq: unif converg 2}. Finally, by Assumption~\ref{ass: estimator of m^t}(ii), we have~$\lim_{n \rightarrow \infty} \sup_{x \in \mathcal{X}} \max_{t \in \mathcal{T}}\; \left|   m^t(x) - \hat{m}_n^t(x)   \right| = 0$. As this claim holds irrespective of the value of~$\mu \in \mathcal{S}$ and in view of \eqref{eq: unif converg 2}, we conclude that \eqref{eq: unif converg term 2} holds. 

The findings above imply that~$\lim_{n \rightarrow \infty} \sup_{\mu \in \mathcal{S}} \left| \nu(\mu) - \hat{\nu}_{n}(\mu) \right| = 0$. The claim thus follows. \Halmos
\endproof

Next, we use Lemmas~\ref{lemma: lemma compact solution set}--\ref{lemma: lemma unif converg} to show that the optimal value~$\hat{\nu}_{n}^\star$ of the sample-approximate dual problem \eqref{eq: sample-based assignment simplified} converges almost surely to the optimal value~$\nu^\star$ of the true dual problem \eqref{eq:bestupperboundproblem}. We will also prove a convergence result in terms of their respective optimal solution sets. The following definitions will be relevant for this result. We use~$\Vert \cdot \Vert$ to denote the Euclidean norm and define 
\begin{equation*}
   \begin{aligned}
   & \text{dist}(x, \mathcal{A}) = \inf_{x^{'} \in \mathcal{A}} \Vert x - x' \Vert,
   && \mathbb{D}(\mathcal{A}, \mathcal{B}) = \sup_{x \in \mathcal{A}} \text{dist}(x, \mathcal{B}),
   \end{aligned}
\end{equation*}
where~$\text{dist}(x, \mathcal{A})$ denotes the distance between a point~$x \in \mathbb{R}^d$ and a set~$\mathcal{A} \subseteq \mathbb{R}^d$, and~$\mathbb{D}(\mathcal{A},\mathcal{B})$ denotes the deviation between two sets~$\mathcal{A} \subseteq \mathbb{R}^d$ and~$\mathcal{B} \subseteq \mathbb{R}^d$. 

\begin{lemma} \label{lemma: convergence of sample estim}
Under Assumptions~\ref{ass: meanoutcomes} and~\ref{ass: estimator of m^t}, the optimal value~$\hat{\nu}_{n}^\star$ of problem \eqref{eq: sample-based assignment simplified} converges almost surely to the optimal value~$\nu^\star$ of problem \eqref{eq:bestupperboundproblem}, and the deviation~$\mathbb{D}(\hat{\mathcal{S}}_{n}^\star, \mathcal{S}^\star)$ between their respective optimal solution sets~$\hat{\mathcal{S}}_{n}^\star$ and~$\mathcal{S}^\star$ converges almost surely to zero as~$n$ tends to infinity, that is,~$\lim_{n \rightarrow \infty} \vert \nu^\star - \hat{\nu}_{n}^\star \vert = 0$ and~$\lim_{n \rightarrow \infty} \mathbb{D}(\hat{\mathcal{S}}_{n}^\star, \mathcal{S}^\star) = 0$.
\end{lemma}

\proof{Proof.}

Consider the compact set~$[0,  \max(C, \hat{C})]^{m+1} \subset \mathbb{R}_+^{m+1}$, and note that: (i) the optimal solution set~$\mathcal{S}^\star$ is nonempty and contained in~$[0, \max(C, \hat{C})]^{m+1}$ by Lemma~\ref{lemma: lemma compact solution set}; (ii) the function~${\nu}(\mu)$ is finite valued by Assumption~\ref{ass: meanoutcomes}(i) and continuous on~$[0, \max(C, \hat{C})]^{m+1}$; (iii)~$\hat{\nu}_{n}(\mu)$ converges uniformly on~$[0, \max(C, \hat{C})]^{m+1}$ and almost surely to~${\nu}(\mu)$ by Lemma~\ref{lemma: lemma unif converg}; and, (iv)~$\hat{\mathcal{S}}_{n}^\star$ is nonempty and contained in~$[0, \max(C, \hat{C})]^{m+1}$ by Lemma~\ref{lemma: sample compact solution set alt}. We can thus invoke Theorem 5.3 of~\cite{shapiro} and conclude that~$\hat{\nu}_{n}^\star$ converges almost surely to~$\nu^\star$, and~$\mathbb{D}(\hat{\mathcal{S}}_{n}^\star, \mathcal{S}^\star)$ converges almost surely to zero.
\Halmos
\endproof

The proof of Theorem~\ref{theorem: asymptotic property of sample-based dual-price queuing policy} will require comparing the treatment assignments made by a dual-price queuing policy and its sample approximation defined in Definitions~\ref{def: dual-price queuing policy} and~\ref{def: sample dual-price queuing policy}, respectively. To that end, we extend the sample solution set convergence results of Lemma~\ref{lemma: convergence of sample estim} by showing that for~$n$ large enough, any optimal sample-based solution~$\hat{\mu}_{n}^\star$ will be arbitrarily close to some element within the optimal solution set~$\mathcal{S}^\star$. This will help ensure that almost surely for~$n$ large enough, a sample approximation policy will make the same treatment assignments as a dual-price queuing policy. 

\begin{lemma} \label{lemma: opt mu close}
Under Assumptions~\ref{ass: meanoutcomes} and~\ref{ass: estimator of m^t}, there exists~$\Omega' \subseteq \Omega$, where~$\mathbb P(\Omega') = 1$, with the following property. For any~$\omega \in \Omega'$ and~$\epsilon > 0$, there exists~$N(\omega) \in \mathbb{N}$ such that
\begin{equation*}
    \begin{aligned}
    \forall \hat{\mu}^\star_{n}(\omega) \in \hat{\mathcal{S}}_{n}^\star(\omega),\, \exists \mu^\star \in \mathcal{S}^\star \,:\, \Vert  \hat{\mu}^\star_{n}(\omega) - \mu^\star \Vert < \epsilon \quad  \forall n \geq N(\omega).
    \end{aligned}
\end{equation*}
\end{lemma}
\proof{Proof.}
 By Lemma~\ref{lemma: convergence of sample estim},~$\mathbb D (\hat{\mathcal{S}}_{n}^\star, \mathcal{S}^\star)$ converges almost surely to zero as~$n \rightarrow \infty$. In other words, there exists~$\Omega' \subseteq \Omega$, where~$\mathbb P(\Omega') = 1$, with the following property. For any~$\omega \in \Omega'$ and~$\epsilon > 0$, there exists~$N(\omega) \in \mathbb{N}$ such that 
\begin{equation*}
    \begin{aligned}
     \sup_{\hat{\mu}_{n}(\omega) \in \hat{\mathcal{S}}_{n}^\star(\omega)}\, \inf_{\mu \in \mathcal{S}^\star} \, \Vert \hat{\mu}_{n}(\omega) - \mu \Vert < \epsilon \quad\forall n \geq N(\omega),
    \end{aligned}
\end{equation*}
which follows by the definition of~$\mathbb D(A, B)$. The above implies that
\begin{equation}\label{eq:proof of Lemma 7 eq 1}
    \begin{aligned}
    \forall \hat{\mu}_{n}(\omega) \in \hat{\mathcal{S}}_{n}^\star(\omega) \,:\, \inf_{\mu \in \mathcal{S}^\star} \, \Vert \hat{\mu}_{n}(\omega) - \mu \Vert < \epsilon \quad \forall n \geq N(\omega).
    \end{aligned}
\end{equation}
Next, we show that for every~$\omega \in \Omega$ (i.e., for fixed historical sample),~$\inf_{\mu \in \mathcal{S}^\star} \Vert \hat{\mu}_{n}(\omega) - \mu \Vert$ has an optimal solution for all~$\hat{\mu}_{n}(\omega) \in \hat{\mathcal{S}}_{n}^\star(\omega)$. To this end, we first show that~$\mathcal{S}^\star$ is a compact set. Recall that~$\mathcal{S}^\star$ denotes the optimal solution set to problem \eqref{eq:bestupperboundproblem}. The set~$\mathcal{S}^\star$ is bounded as it is contained in the compact set~$[0, C]^{m+1}$ by Lemma~\ref{lemma: lemma compact solution set}, where~$C$ is defined as in Assumption~\ref{ass: meanoutcomes}(i). The objective function~$\nu(\mu)$ of problem \eqref{eq:bestupperboundproblem} is convex and thus continuous on~$\mathbb R^{m+1}$. This implies that it has closed sublevel sets of the form~$\{\mu \in \mathbb{R}^{m+1}_+ \mid \nu(\mu) \leq c\}$, where~$c$ is any scalar. We now have~$\mathcal{S}^\star = [0, C]^{m+1} \cap \{\mu \in \mathbb{R}^{m+1}_+ \mid \nu(\mu) \leq \nu^\star\}$, which is closed as it is an intersection of closed sets. The set~$\mathcal{S}^\star$ is thus compact.
 
 Since all norms are continuous on~$\mathbb{R}^{m+1}_+$ and~$\mathcal{S}^\star$ is compact, by Weierstrass' Extreme Value Theorem (Proposition A.2.7 of~\cite{bertsekas2009convex}),~$\inf_{\mu \in \mathcal{S}^\star} \Vert \hat{\mu}_{n}(\omega) - \mu \Vert$ has an optimal solution for any input~$\hat{\mu}_{n}(\omega) \in \hat{\mathcal{S}}_{n}^\star(\omega)$. This together with \eqref{eq:proof of Lemma 7 eq 1} implies that for any~$\omega \in \Omega'$ and~$\epsilon > 0$, there exists~$N(\omega) \in \mathbb{N}$ such that 
\begin{equation*}
    \begin{aligned}
    \forall \hat{\mu}^\star_{n}(\omega) \in  \hat{\mathcal{S}}_{n}^\star(\omega),\, \exists \mu^\star \in \mathcal{S}^\star \,:\, \Vert  \hat{\mu}^\star_{n} - \mu^\star \Vert < \epsilon \quad \forall n \geq N(\omega).
    \end{aligned}
\end{equation*}
The claim thus follows as~$\mathbb P (\Omega') = 1$. \Halmos
\endproof
We are now ready to prove Theorem~\ref{theorem: asymptotic property of sample-based dual-price queuing policy}.

\proof{Proof of Theorem~\ref{theorem: asymptotic property of sample-based dual-price queuing policy}.}
We use the following notation throughout the proof. Define the function~$h^{t}$ through~$h^{t}(x, \mu) = m^t(x) - \mu^t$ for all~$x \in \mathcal X$,~$\mu \in \mathbb R_+^{m+1}$, and~$t \in \mathcal{T}$. For any~$\omega \in \Omega$ (i.e., fixed historical samples), define the function~$\hat{h}_n^{t}(\cdot, \omega)$ through~$\hat{h}_n^{t}(x, \mu, \omega) = \hat{m}_n^t(x, \omega) - \mu^t$ for all~$x \in \mathcal X$,~$\mu \in \mathbb R_+^{m+1}$,~$t \in \mathcal{T}$, and~$n \in \mathbb N$. 
Denote by~$z^{\rm{D}}$ the expected average outcome of a dual-price queuing policy. Recall that~$\hat{z}^{\rm{D}}_n$ and~$z^\star$ denote the expected average outcome of a sample-based dual-price queuing policy,  where the expectation is taken with respect to the distribution of covariates and outcomes of the individuals arriving during implementation and \emph{not} with respect to the distribution of the historical samples, and the optimal value of problem \eqref{eq:UpperBoundTrueP}, respectively. Note that, by definition,~$\hat{z}^{\rm{D}}_n$ is a random variable as it depends on the historical samples observed. 

The proof is divided into two steps. In Step 1, we show that the allocations of a sample-based dual-price queuing policy match almost surely those of a dual-price queuing policy as~$n \rightarrow \infty$. In Step 2, we then use this fact to show that a sample-based dual-price queuing policy asymptotically generates at least as high expected average outcome as that of a dual-price queuing policy as~$n \rightarrow \infty$ almost surely. As the expected average outcome of any dual-price queuing policy is equal to~$z^\star$ by Theorem~\ref{theorem: asymptotic opt with known values}, we will thus conclude that the claim holds.

\textbf{Step 1.} 
For any~$\omega \in \Omega$,~$n \in \mathbb N$,~$\hat{\mu}_{n}^\star(\omega) \in \hat{\mathcal{S}}_{n}^\star(\omega)$, and~$\mu^\star \in \mathcal{S}^\star$, by triangle inequality we have
\begin{equation}\label{eq: proof of main theorem eq 1}
    \begin{aligned}
    \sup_{t \in \mathcal{T}} \,\sup_{x \in \mathcal{X}} \, \vert \hat{h}^t_{n}(x, \hat{\mu}_{n}^{\star}(\omega), \omega) - h^t(x, \mu^\star) \vert &= \sup_{t \in \mathcal{T}} \,\sup_{x \in \mathcal{X}} |\hat{m}_n^t(x, \omega) - \hat{\mu}_{n}^{\star, t}(\omega) - m^t(x) + \mu^{\star,t}|\\
    &\leq \sup_{t \in \mathcal{T}} \,\sup_{x \in \mathcal{X}} |\hat{m}_n^t(x, \omega) - m^t(x)| + \sup_{t \in \mathcal{T}} |\hat{\mu}_{n}^{\star, t}(\omega) - \mu^{\star,t}|.
    \end{aligned}
\end{equation}
For the first term in the last line of \eqref{eq: proof of main theorem eq 1}, by Assumption~\ref{ass: estimator of m^t}(ii), there exists~$\Omega' \subseteq \Omega$, where~$\mathbb P(\Omega') = 1$, with the following property. For every~$\omega \in \Omega'$ and~$\epsilon > 0$, there exists~$N(\omega) \in \mathbb{N}$ such that~$$\sup_{t \in \mathcal T}\,\sup_{x \in \mathcal{X}}\, |\hat{m}_n^t(x, \omega) - m^t(x)| < \epsilon \quad\forall n \geq N(\omega).$$ 
For the second term in the last line of \eqref{eq: proof of main theorem eq 1}, Lemma~\ref{lemma: opt mu close} and the property~$\sup_{t \in \mathcal{T}} |\mu^{t}| \leq \Vert \mu \Vert$ imply that there exists~$\Omega'' \subseteq \Omega$, where~$\mathbb P(\Omega'') = 1$, with the following property. For every~$\omega \in \Omega''$ and~$\epsilon > 0$, there exists~$N(\omega) \in \mathbb{N}$ such that
\begin{equation*}
    \begin{aligned}
    \forall \hat{\mu}_{n}^\star(\omega) \in \hat{\mathcal{S}}_{n}^\star(\omega), \exists \mu^\star \in \mathcal{S}^\star\,:\, \sup_{t \in \mathcal{T}} |\hat{\mu}_{n}^{\star, t}(\omega) - \mu^{\star,t}| < \epsilon  \quad \forall n \geq N(\omega),
    \end{aligned}
\end{equation*}
Let~$\hat \Omega = \Omega' \cap \Omega''$ and note that~$\mathbb P(\hat \Omega) = 1$ because
\begin{equation*}
    \begin{aligned}
    \mathbb P(\hat \Omega) \geq 1 - \mathbb P((\Omega \setminus \Omega') \cup (\Omega \setminus \Omega'')) \geq 1 - \mathbb P((\Omega \setminus \Omega')) - \mathbb P((\Omega \setminus \Omega'')) = 1,
    \end{aligned}
\end{equation*}
where the second inequality follows from the union bound, and the equality holds because~$\mathbb P(\Omega') = 1$ and~$\mathbb P(\Omega'') = 1$. Assumption~\ref{ass: estimator of m^t}(ii) and  Lemma~\ref{lemma: opt mu close} thus imply that for every~$\omega \in \hat \Omega$ and~$\epsilon > 0$, there exists~$N(\omega) \in \mathbb{N}$ such that
\begin{equation*}
    \begin{aligned}
    \forall \hat{\mu}_{n}^\star(\omega) \in \hat{\mathcal{S}}_{n}^\star(\omega), \exists \mu^\star \in \mathcal{S}^\star\,:\, 
    \sup_{t \in \mathcal{T}} \,\sup_{x \in \mathcal{X}} |\hat{m}_n^t(x, \omega) - m^t(x)| + \sup_{t \in \mathcal{T}} |\hat{\mu}_{n}^{\star, t}(\omega) - \mu^{\star,t}| < \epsilon  \quad \forall n \geq N(\omega).
    \end{aligned}
\end{equation*}
The above expression, in conjunction with equation \eqref{eq: proof of main theorem eq 1}, imply that for every~$\omega \in \hat \Omega$ and~$\epsilon > 0$, there exists~$N(\omega) \in \mathbb{N}$ such that
\begin{equation}\label{eq: sup over t sup norm func}
    \begin{aligned}
    \forall \hat{\mu}_{n}^\star(\omega) \in \hat{\mathcal{S}}_n^\star(\omega), \exists \mu^\star \in \mathcal{S}^\star\,:\, \sup_{t \in \mathcal{T}} \,\sup_{x \in \mathcal{X}} \, \vert \hat{h}^t_{n}(x, \hat{\mu}_{n}^{\star}(\omega), \omega) - h^t(x, \mu^\star) \vert < \epsilon \quad \forall n \geq N(\omega).
    \end{aligned}
\end{equation}
Decomposing the expression~$\sup_{t \in \mathcal{T}} \,\sup_{x \in \mathcal{X}} \, \vert \hat{h}^t_{n}(x, \hat{\mu}_{n}^{\star}(\omega), \omega) - h^t(x, \mu^\star) \vert < \epsilon$ in Equation \eqref{eq: sup over t sup norm func}, we have that for every~$\omega \in \hat \Omega$ and~$\epsilon > 0$, there exists~$N(\omega) \in \mathbb{N}$ such that
\begin{equation}\label{eq: split_abs_val}
\begin{aligned}
    \forall \hat{\mu}_{n}^\star(\omega) \in \hat{\mathcal{S}}_{n}^\star(\omega), \exists \mu^\star \in \mathcal{S}^\star \,:\,
    \hat{h}^t_{n}(x, \hat{\mu}_{n}^{\star}(\omega), \omega) - h^t(x, \mu^\star) < \epsilon, \;\;
    h^t(x, \mu^\star) - \hat{h}^t_{n}(x, \hat{\mu}_{n}^{\star}(\omega), \omega) < \epsilon \\
    \forall t \in \mathcal T, \forall x \in \mathcal X, \forall n \geq N(\omega).
       \end{aligned}
\end{equation}
For any~$t, t' \in \mathcal T$, by summing the inequalities from above, i.e.,
\begin{equation*}
     h^{t}(x, \mu^\star) - \hat{h}^{t}_{n}(x, \hat{\mu}_{n}^{\star}(\omega), \omega) < \epsilon \quad \text{and} \quad \hat{h}^{t'}_{n}(x, \hat{\mu}_{n}^{\star}(\omega), \omega) - h^{t'}(x, \mu^\star) < \epsilon,
\end{equation*}
we obtain
\begin{equation}\label{eq: h_func}
    h^t(x, \mu^\star) - \hat{h}^t_{n}(x, \hat{\mu}_{n}^{\star}(\omega), \omega) + \hat{h}^{t'}_{n}(x, \hat{\mu}_{n}^{\star}(\omega), \omega) - h^{t'}(x, \mu^\star) < 2 \epsilon.
\end{equation}
Using \eqref{eq: split_abs_val} and \eqref{eq: h_func}, we have that for every~$\omega \in \hat \Omega$ and~$\epsilon > 0$, there exists~$N(\omega) \in \mathbb{N}$ such that
\begin{equation}\label{eq: main proof eq 2}
    \begin{aligned}
    \forall \hat{\mu}_{n}^\star(\omega) \in \hat{\mathcal{S}}_{n}^\star(\omega), \exists \mu^\star \in \mathcal{S}^\star
    \,:\, h^t(x, \mu^\star) - \hat{h}^t_{n}(x, \hat{\mu}_{n}^{\star}(\omega), \omega) +  
    \hat{h}^{t'}_{n}(x, \hat{\mu}_{n}^{\star}(\omega), \omega) - h^{t'}(x, \mu^\star) < 2 \epsilon \\
    \forall t, t' \in \mathcal T, \forall x \in \mathcal X, \forall n \geq N(\omega).
    \end{aligned}
\end{equation}
Let~$t(x, \mu) = \min \arg \max_{t \in \mathcal T} h^{t}(x, {\mu})$, and define~$\delta(x, \mu) = h^{t(x, \mu)}(x, {\mu}) - \max_{t \neq t(x, \mu)} h^{t}(x, {\mu})$ for all~$x \in \mathcal X$ and~$\mu \in \mathbb R_+^{m+1}$. Note that~$\delta(x, \mu)$ represents the difference between the highest and second highest value of the set~$\{ h^{t}(x, {\mu}) \}_{t \in \mathcal{T}}$ and that we have~$\delta(x, \mu) > 0$ if and only if~$\arg \max_{t \in \mathcal T} h^{t}(x, {\mu}) = \{t(x, \mu)\}$, i.e.,~$t(x, \mu)$ is the unique treatment that attains the maximum. 
Next, we write \eqref{eq: main proof eq 2} not for all treatment pairs~$t, t' \in \mathcal T$ but only for the treatment pairs of the form~$(t(x, \mu^\star), t)$, where~$t \in \mathcal T \setminus \{t(x, \mu^\star)\}$. Using the new notation we introduced and \eqref{eq: main proof eq 2}, we have that for every~$\omega \in \hat \Omega$ and~$\epsilon > 0$, there exists~$N(\omega) \in \mathbb{N}$
\begin{equation}\label{eq: main proof eq 3}
    \begin{aligned}
    \forall \hat{\mu}_{n}^\star(\omega) \in \hat{\mathcal{S}}_{n}^\star(\omega), \exists \mu^\star \in \mathcal{S}^\star
    \,:\, \hat{h}^{t}_{n}(x, \hat{\mu}_{n}^{\star}(\omega), \omega) - \hat{h}^{t(x, \mu^\star)}_{n}(x, \hat{\mu}_{n}^{\star}(\omega), \omega)   < 2 \epsilon  - \delta(x, \mu^\star) \\
    \forall t \in \mathcal T \setminus \{t(x, \mu^\star)\}, \forall x \in \mathcal X, \forall n \geq N(\omega),
    \end{aligned}
\end{equation}
where we used the fact that 
\begin{equation*}
     \delta(x, \mu^\star) = h^{t(x, \mu)}(x, {\mu^\star}) - \max_{t \neq t(x, \mu^\star)} h^{t}(x, {\mu^\star}) \leq h^{t(x, \mu^\star)}(x, {\mu^\star}) - h^{t}(x, {\mu^\star}) \quad \forall t \in \mathcal T \setminus \{t(x, \mu^\star)\}, \forall x \in \mathcal X.
\end{equation*}
Equation \eqref{eq: main proof eq 3} implies that for a fixed~$\omega \in \hat{\Omega}$ and~$x \in \mathcal{X}$, if~$\delta(x, \mu^\star) > 0$, then for~$\epsilon > 0$ small enough (an therefore~$n$ large enough),~$2 \epsilon  - \delta(x, \mu^\star) < 0$. This in turn implies~$\hat{h}^{t}_{n}(x, \hat{\mu}_{n}^{\star}(\omega), \omega) - \hat{h}^{t(x, \mu^\star)}_{n}(x, \hat{\mu}_{n}^{\star}(\omega), \omega) < 0$ for all~$\forall t \in \mathcal T \setminus \{t(x, \mu^\star)\}$. Since~$t(x, \mu^\star)$ is also the unique maximizer of~$\arg \max_{t \in \mathcal T} h^{t}(x, {\mu})$ when~$\delta(x, \mu^\star) > 0$, this implies that the allocation of a sample-based dual-price queuing policy for an individual with covariates~$x$ match that of a dual-price queuing policy for~$n$ large enough for any~$\omega \in \hat \Omega$, where~$\mathbb P(\hat \Omega) = 1$. By Assumption~\ref{ass: meanoutcomes} and Lemma~\ref{lemma: z r.v. cont dist}, we have~$\delta(X, \mu^\star) > 0$ almost surely. We will use this observation in Step 2 to prove that the expected average outcome of a sample-based dual-price policy almost surely will be as high as that of a true dual-price policy as~$n \rightarrow \infty$.

\textbf{Step 2.} 
By Step 1, for any~$\omega \in \hat \Omega$ and sequence~$(\hat \mu^\star_{n}(\omega))_{n \in \mathbb N}$ of sample-based solutions, where~$\hat \mu^\star_{n}(\omega) \in \hat{\mathcal{S}}_{n}^\star(\omega)$ for all~$n \in \mathbb N$, there exists a sequence~$(\mu^\star_{n}(\omega))_{n \in \mathbb N}$ of solutions, where ~$\mu^\star_{n}(\omega) \in \mathcal{S}^\star$ for all~$n \in \mathbb N$, with the following property. For every~$\epsilon > 0$, there exists~$N(\omega) \in \mathbb{N}$ such that 
\begin{equation}\label{eq: main proof eq 4}
    \begin{aligned}
    \hat{h}^{t}_{n}(x, \hat{\mu}_{n}^{\star}(\omega), \omega) - \hat{h}^{t(x, \mu_{n}^\star(\omega))}_{n}(x, \hat{\mu}_{n}^{\star}(\omega), \omega)   < 2 \epsilon  - \delta(x, \mu_{n}^\star(\omega)) \\
    \forall t \in \mathcal T \setminus \{t(x, {\mu}_{n}^\star(\omega))\}, \forall x \in \mathcal X, \forall n \geq N(\omega).
    \end{aligned}
\end{equation}
Fix a~$\omega \in \hat \Omega$ and sequence couple~$(\hat \mu^\star_{n}(\omega))_{n \in \mathbb N}$ and~$(\mu^\star_{n}(\omega))_{n \in \mathbb N}$ such that \eqref{eq: main proof eq 4} holds. From now on, we slightly abuse the notation and drop~$\omega$ notation. In the following, intuitively, we fix a path of historical samples, which implies that our estimates~$\hat \mu^\star_{n}$,~$\hat m^t_n$,~$\hat h^t_{n}$ and~$\hat{z}^{\rm{D}}_{n}$ are no longer random. We still perceive the covariates and outcomes of the individuals arriving during the implementation as random and will denote by~$\mathbb E_{\text x}$ the expectation taken with respect to the distribution of these individuals' covariates.
Recalling that~$z^{\rm{D}} = z^\star$ by Theorem~\ref{theorem: asymptotic opt with known values}, we have 
\begin{equation}\label{eq: main proof eq 5}
    \begin{aligned}
    &\limsup_{n \rightarrow \infty} \;z^\star - \hat{z}^{\rm{D}}_{n} = 
    \limsup_{n \rightarrow \infty} \;z^{\rm{D}} - \hat{z}^{\rm{D}}_{n}\\
    &\leq \limsup_{n \rightarrow \infty} \;\mathbb E_{\text x}\left[C \mathbbm{1} \left\{\min \arg \max_{t \in \mathcal T} h^{t}(X, {\mu}_{n}^\star) \neq \min \arg \max_{t \in \mathcal T} \hat h_{n}^{t}(X, \hat{\mu}_{n}^\star) \right\}  \right]\\
    &\leq \;\mathbb E_{\text x}\left[ \limsup_{n \rightarrow \infty} C \mathbbm{1} \left\{\min \arg \max_{t \in \mathcal T} h^{t}(X, {\mu}_{n}^\star) \neq \min \arg \max_{t \in \mathcal T} \hat h_{n}^{t}(X, \hat{\mu}_{n}^\star) \right\}  \right]\\
    &\leq \;\mathbb E_{\text x}\left[ \limsup_{n \rightarrow \infty} C \mathbbm{1} \left\{\hat h_{n}^{t}(X, \hat{\mu}_{n}^\star)  - \hat h_{n}^{t(X, {\mu}_{n}^\star)}(X, \hat{\mu}_{n}^\star) \geq 0 \;\;\text{for some}\; t \in \mathcal T \setminus \{t(X, {\mu}_{n}^\star)\} \right\}  \right]\\
    &\leq \;\mathbb E_{\text x}\left[ \limsup_{n \rightarrow \infty} C \mathbbm{1} \left\{\hat h_{n}^{t}(X, \hat{\mu}_{n}^\star)  - \hat h_{n}^{t(X, {\mu}_{n}^\star)}(X, \hat{\mu}_{n}^\star) + \delta(X, \mu_{n}^\star) > 0 \;\;\text{for some}\; t \in \mathcal T \setminus \{t(X, {\mu}_{n}^\star)\} \right\}  \right],
    \end{aligned}
\end{equation}
where the first inequality holds because a gap between~$z^{\rm{D}}$ and~$z^{\rm{D}}_{n}$ occurs only when the dual-price queuing policy and its sample-based variant assign different treatments to an individual with covariates~$X$ and because~$m^t$ is uniformly bounded by~$C$ by Assumption~\ref{ass: meanoutcomes} (i), which implies that the gap is also bounded by the same value. The second inequality follows from the reverse Fatou's lemma, which applies because the indicator function is bounded above by~$1$. The third inequality holds because for every~$x \in \mathcal X$,~$\min \arg \max_{t \in \mathcal T} h^{t}(x, {\mu}_{n}^\star) \neq \min \arg \max_{t \in \mathcal T} \hat h_{n}^{t}(x, \hat{\mu}_{n}^\star)$ implies that there exists a~$t \in \mathcal T \setminus \{t(x, {\mu}_{n}^\star)\}$ such that~$\hat h_{n}^{t}(x, \hat{\mu}_{n}^\star)  - \hat h_{n}^{t(x, {\mu}_{n}^\star)}(x, \hat{\mu}_{n}^\star) \geq 0$, where~$t(x, {\mu}_{n}^\star)$ is defined as before, i.e.,~$t(x, {\mu}_{n}^\star) = \min \arg \max_{t \in \mathcal T} h^{t}(x, {\mu}_{n}^\star)$. To see this, for a fixed~$x \in \mathcal{X}$, if~$t(x, {\mu}_{n}^\star) \not \in \arg \max_{t \in \mathcal T} \hat h_{n}^{t}(x, \hat{\mu}_{n}^\star)$, then there exists~$t \in \mathcal T \setminus \{t(x, {\mu}_{n}^\star)\}$ such that~$\hat h_{n}^{t}(x, \hat{\mu}_{n}^\star)  - \hat h_{n}^{t(x, {\mu}_{n}^\star)}(x, \hat{\mu}_{n}^\star) > 0$. On the other hand, if~$t(x, {\mu}_{n}^\star) \in \arg \max_{t \in \mathcal T} \hat h_{n}^{t}(x, \hat{\mu}_{n}^\star)$ but is not the minimum index treatment, then there exists~$t \in \mathcal T \setminus \{t(x, {\mu}_{n}^\star)\}$ such that~$\hat h_{n}^{t}(x, \hat{\mu}_{n}^\star)  - \hat h_{n}^{t(x, {\mu}_{n}^\star)}(x, \hat{\mu}_{n}^\star) = 0$. Finally, the last inequality holds because~$\hat h_{n}^{t}(X, \hat{\mu}_{n}^\star)  - \hat h_{n}^{t(X, {\mu}_{n}^\star)}(X, \hat{\mu}_{n}^\star) \geq 0$ implies that~$\hat h_{n}^{t}(X, \hat{\mu}_{n}^\star)  - \hat h_{n}^{t(X, {\mu}_{n}^\star)}(X, \hat{\mu}_{n}^\star) + \delta(X, \mu_{n}^\star) > 0$ almost surely as~$\delta(X, \mu_{n}^\star) > 0$ almost surely by Assumption~\ref{ass: meanoutcomes} and Lemma~\ref{lemma: z r.v. cont dist}.

We will next show that the last expression in \eqref{eq: main proof eq 5} equals zero. By \eqref{eq: main proof eq 4}, we have
\begin{equation} \label{eq: limsup_h}
    \begin{aligned}
    \limsup_{n \rightarrow \infty} \hat{h}^{t}_{n}(x, \hat{\mu}_{n}^{\star}) - \hat{h}^t_{n}(x, \hat{\mu}_{n}^{\star}) + \delta(x, {\mu}_{n}^\star) \leq 0  \quad \forall t \in \mathcal T \setminus \{t(x, {\mu}_{n}^\star)\}, x \in \mathcal X\newck{.}
    \end{aligned}
\end{equation}
%
Equation \eqref{eq: limsup_h} implies that for any~$x \in \mathcal{X}$, we also have~$$\liminf_{n \rightarrow \infty} \mathbbm{1} \left\{\hat{h}^{t}_{n}(x, \hat{\mu}_{n}^{\star}) - \hat{h}^t_{n}(x, \hat{\mu}_{n}^{\star}) + \delta(x, {\mu}_{n}^\star) \leq 0  \quad \forall t \in \mathcal T \setminus \{t(x, {\mu}_{n}^\star)\} \right\} = 1.$$ 
We then have
\begin{equation*}
    \begin{aligned}
    &\limsup_{n \rightarrow \infty} \mathbbm{1} \left\{\hat h_{n}^{t}(x, \hat{\mu}_{n}^\star)  - \hat h_{n}^{t(x, {\mu}_{n}^\star)}(x, \hat{\mu}_{n}^\star) + \delta(x, \mu_{n}^\star) > 0 \;\;\text{for some}\; t \in \mathcal T \setminus \{t(x, {\mu}_{n}^\star)\newck{\}} \right\} \\
    =\; &1 - \liminf_{n \rightarrow \infty} \mathbbm{1} \left\{\hat{h}^{t}_{n}(x, \hat{\mu}_{n}^{\star}) - \hat{h}^t_{n}(x, \hat{\mu}_{n}^{\star}) + \delta(x, {\mu}_{n}^\star) \leq 0  \quad \forall t \in \mathcal T \setminus \{t(x, {\mu}_{n}^\star)\} \right\} \\
    =\; &0.
    \end{aligned}
\end{equation*}
As the above equality holds true for all~$x \in \mathcal X$, the last expression in \eqref{eq: main proof eq 5} therefore amounts to zero. By \eqref{eq: main proof eq 5}, we thus have~$\limsup_{n \rightarrow \infty} \;z^\star - \hat{z}^{\rm{D}}_n \leq 0$ for the fixed path of historical samples. The claim now follows from the observation that our choice of sequence couple~$(\hat \mu^\star_n)_{n \in \mathbb N}$ and~$(\mu^\star_n)_{n \in \mathbb N}$ was arbitrary and that~$\limsup_{n \rightarrow \infty} \;z^\star - \hat{z}^{\rm{D}}_n \leq 0$ holds for every sequence couple~$(\hat \mu^\star_n)_{n \in \mathbb N} = (\hat \mu^\star_n)_{n \in \mathbb N}(\omega)$ and~$(\mu^\star_n)_{n \in \mathbb N} = (\mu^\star_n)_{n \in \mathbb N}(\omega)$, where~$\omega \in \hat \Omega$ and~$\mathbb P(\hat \Omega) = 1$. In other words, we have~$\limsup_{n \rightarrow \infty} \;z^\star - \hat{z}^{\rm{D}}_n \leq 0$ almost surely.
\Halmos
\endproof

\section{Derivations of Fairness-Constrained Policies for Section~\ref{sec: fairness extensions}}\label{sec: fairness extensions derivations}

\subsection{Allocation-Parity-Constrained Policies}
{\color{blue} In this section, we show the derivation of policy~$\hat{\pi}^{\star}_{\rm alloc,n}$ of Section~\ref{sec: sp in alloc}. The derivation for the minority prioritization case is similar and omitted. Recall that we now explicitly consider the protected feature~$G$ as an input into the policy~$\pi$, i.e.,~$X = (X^{-G}, G)$. The sample approximation of the allocation-parity-constrained policy design problem is given by
\begin{equation}\label{eq:UpperBoundTrueP_sp} 
    \begin{aligned}
    \hat{z}_{\rm alloc, n}^\star = &\max_{\pi \in \Pi} &&\frac{1}{n} \sum_{i=1}^n \sum_{t \in \mathcal T} \pi^t(x_i) \hat{m}_n^t(x_i) \\
    &\;\;\,\text{s.t.} &&\frac{1}{n} \sum_{i=1}^n  \pi^t(x_i) \leq b^t \quad\forall t \in \mathcal T \\
    &&& \frac{1}{n_g} \sum_{i=1}^n \mathbbm{1}[g_i = g] \pi^t(x_i) - \frac{1}{n_{g'}} \sum_{i=1}^n \mathbbm{1}[g_i = g'] \pi^t(x_i) \leq \delta \;\;\forall g, g' \in \mathcal{G}, g \neq g', t \in \mathcal T,
    \end{aligned}
\end{equation}
where $n_g$ denotes the number of individuals with protected feature $g$ in the data. Let~$\mu \in \mathbb R_+^{m+1}$ collect the dual variables of the capacity constraints and denote by~$\lambda^t(g,g') \in \mathbb R_+$ the dual variable of the fairness constraint for pair~$(g,g')$ and for treatment~$t$. 
The Lagrangian dual of problem~\eqref{eq:UpperBoundTrueP_sp} is given by 
\begin{equation}\label{eq:relaxation_sp}
    \begin{aligned}
    \hat{v}^\star_{\rm alloc, n} &= \min_{\substack{\mu \in \mathbb R_+^{m+1} \\ \lambda \in {\mathcal L}_\infty(\mathcal G \times \mathcal G, \mathbb R_+^{m+1}) } } \max_{\pi \in \Pi} \;\; L_{\rm alloc,n}(\pi, \mu, \lambda),
    \end{aligned}
\end{equation}
where
\begin{align*}
    L_{\rm alloc,n}(\pi, \mu, \lambda)
    &= \frac{1}{n} \sum_{i=1}^n \sum_{t \in \mathcal T} \pi^t(x_i) \hat{m}_n^t(x_i) + \sum_{t \in \mathcal T} \mu^t \left( b^t - \frac{1}{n} \sum_{i=1}^n  \pi^t(x_i) \right) + \sum_{t \in \mathcal T} \sum_{g \in \mathcal G} \sum_{g' \in \mathcal G, g \neq g'} \lambda^t(g,g') \delta \\
    &- \sum_{t \in \mathcal T} \sum_{g \in \mathcal G} \sum_{g' \in \mathcal G, g \neq g'} \lambda^t(g,g') \left( \frac{1}{n_g} \sum_{i=1}^n \mathbbm{1}[g_i = g] \pi^t(x_i) - \frac{1}{n_{g'}} \sum_{i=1}^n \mathbbm{1}[g_i = g'] \pi^t(x_i) \right).
\end{align*}
We can rewrite the last term as
\begin{equation*}
    \begin{aligned}
        &\sum_{t \in \mathcal T} \sum_{g \in \mathcal G} \sum_{g' \in \mathcal G, g \neq g'} \lambda^t(g,g') \left( \frac{1}{n_g} \sum_{i=1}^n \mathbbm{1}[g_i = g] \pi^t(x_i) - \frac{1}{n_{g'}} \sum_{i=1}^n \mathbbm{1}[g_i = g'] \pi^t(x_i) \right)\\
    &= \sum_{t \in \mathcal T} \sum_{g \in \mathcal G} \left(\frac{1}{n_g} \sum_{i=1}^n \mathbbm{1}[g_i = g] \pi^t(x_i) \right) \left( \sum_{g' \in \mathcal G, g \neq g'} (\lambda^t(g,g') - \lambda^t(g',g)) \right).
    \end{aligned}
\end{equation*}
For ease of notation, we let~$\gamma^t(g) = \sum_{g' \in \mathcal G, g \neq g'} (\lambda^t(g,g') - \lambda^t(g',g))$. We now have
\begin{equation*}
    \begin{aligned}
    L_{\rm alloc,n}(\pi, \mu, \lambda)
    &= \frac{1}{n} \sum_{i=1}^n \sum_{t \in \mathcal T} \pi^t(x_i) \hat{m}_n^t(x_i) + \sum_{t \in \mathcal T} \mu^t \left( b^t - \frac{1}{n} \sum_{i=1}^n  \pi^t(x_i) \right) \\
    &+ \sum_{t \in \mathcal T} \sum_{g \in \mathcal G} \sum_{g' \in \mathcal G, g \neq g'} \lambda^t(g,g') \delta - \sum_{t \in \mathcal T} \sum_{g \in \mathcal G} \left(\frac{1}{n_g} \sum_{i=1}^n \mathbbm{1}[g_i = g] \pi^t(x_i) \right)  \gamma^t(g)\\
    &= \frac{1}{n} \sum_{i=1}^n \sum_{t \in \mathcal T} \pi^t(x_i) \left(\hat{m}_n^t(x_i) - \mu^t - \frac{n}{n_{g_i}} \gamma^t(g_i) \right) + \sum_{t\in \mathcal T} \mu^t  b^t + \sum_{t \in \mathcal T} \sum_{g \in \mathcal G} \sum_{g' \in \mathcal G, g \neq g'} \lambda^t(g,g') \delta.
    \end{aligned}
\end{equation*}
For any fixed~$\mu$ and~$\lambda$, the inner maximization problem of \eqref{eq:relaxation_sp} is solved by 
\begin{equation}\label{eq: policy from Lagrangian relaxation sp}
    \begin{aligned}
    \pi^{t}(x) = \begin{cases}
    1 &\text{if}\;\;\; t = \min \arg \max_{t' \in \mathcal T}\; \left( \hat{m}_n^{t'}(x) - \mu^{t'} - \frac{n}{n_{g}} \gamma^{t'}(g) \right)\\
    0 &\text{otherwise}
    \end{cases}
    \;\forall t \in \mathcal T,\, x \in \mathcal{X}.
    \end{aligned}
\end{equation}
We suppress the dependence of the above policy~$\pi$ on~$\mu, \lambda$ notationally in order to avoid clutter and tie-break using a lexicographic tie-breaker. By substituting the above policy into the dual problem~\eqref{eq:relaxation_sp}, we obtain the convex program
\begin{equation*}\label{eq:bestupperboundproblem_sp}
    \begin{aligned}
    \hat{v}^\star_{\rm alloc,n} &= \min_{\substack{\mu \in \mathbb R_+^{m+1}; \\ \lambda \in {\mathcal L}_\infty(\mathcal G \times \mathcal G, \mathbb R_+^{m+1}) } } \overline{L}_{\rm alloc,n}(\mu, \lambda),
    \end{aligned}
\end{equation*}
where
\begin{equation*} \label{eq:bestupperboundobj_alloc}
    \begin{aligned}
    \overline{L}_{\rm alloc,n}(\mu, \lambda) = &\frac{1}{n} \sum_{i=1}^n \max_{t \in \mathcal T} \left(\hat{m}_n^t(x_i) - \mu^t - \frac{n}{n_{g_i}} \gamma^t(g_i) \right) + \sum_{t\in \mathcal T} \mu^t  b^t + \sum_{t \in \mathcal T} \sum_{g \in \mathcal G} \sum_{g' \in \mathcal G, g \neq g'} \lambda^t(g,g') \delta.
    \end{aligned}
\end{equation*}
We denote by~$\hat{\mu}_n^\star, \hat{\lambda}_n^\star$ an optimal solution to this finite convex program. We can now define our allocation-parity-constrained policy~$\hat{\pi}^\star_{\rm alloc,n}$ by substituting~$\hat{\mu}_n^\star, \hat{\lambda}_n^\star$ into~\eqref{eq: policy from Lagrangian relaxation sp}.}

\subsection{Outcome-Parity-Constrained Policies}
{\color{blue}We now derive the policy~$\hat{\pi}^{\star}_{\rm out, n}$ of Section~\ref{sec: sp in out}. The sample approximation of the outcome-parity-constrained policy design problem is given by
\begin{equation}\label{eq:UpperBoundTrueP_sp_out} 
    \begin{aligned}
    \hat{z}_{\rm out, n}^\star = &\max_{\pi \in \Pi} &&\frac{1}{n} \sum_{i=1}^n \sum_{t \in \mathcal T} \pi^t(x_i) \hat{m}_n^t(x_i) \\
    &\;\;\,\text{s.t.} &&\frac{1}{n} \sum_{i=1}^n  \pi^t(x_i) \leq b^t \quad\forall t \in \mathcal T \\
    &&& \frac{1}{n_g} \sum_{i=1}^n \sum_{t \in \mathcal T} \mathbbm{1}[g_i = g] \pi^t(x_i)\hat{m}_n^t(x_i) - \frac{1}{n_{g'}} \sum_{i=1}^n \sum_{t \in \mathcal T} \mathbbm{1}[g_i = g'] \pi^t(x_i)\hat{m}_n^t(x_i) \leq \delta \\
&&&\qquad\qquad\qquad\qquad\qquad\qquad\qquad\qquad\qquad\qquad\qquad\qquad\;\,\forall g, g' \in \mathcal{G}, g \neq g'.
    \end{aligned}
\end{equation}
Again let~$\mu \in \mathbb R_+^{m+1}$ be the dual variable for the capacity constraints, and denote by~$\lambda(g,g') \in \mathbb R_+$ the dual variable of the statistical parity in outcome constraint for the group pair~$g$ and~$g'$. The Lagrangian relaxation of problem \eqref{eq:UpperBoundTrueP_sp_out} is given by
\begin{align}\label{eq:relaxation_sp_out}
    \hat{v}_{\rm out,n}^\star &= \min_{\substack{\mu \in \mathbb R_+^{m+1} \\ \lambda \in {\mathcal L}_\infty(\mathcal G \times \mathcal G, \mathbb R_+) } } \max_{\pi \in \Pi} \; L_{\rm out,n} (\pi, \mu, \lambda),
\end{align}
where
\begin{align*}
    &L_{\rm out,n} (\pi, \mu, \lambda) = \frac{1}{n} \sum_{i=1}^n \sum_{t \in \mathcal T} \pi^t(x_i) \hat{m}_n^t(x_i) + \sum_{t \in \mathcal T} \mu^t \left( b^t - \frac{1}{n} \sum_{i=1}^n  \pi^t(x_i) \right) + \sum_{g \in \mathcal G} \sum_{g' \in \mathcal G, g \neq g'} \lambda(g,g') \delta \\
    &- \sum_{g \in \mathcal G} \sum_{g' \in \mathcal G, g \neq g'} \lambda(g,g') \left( \frac{1}{n_g} \sum_{i=1}^n \sum_{t \in \mathcal T} \mathbbm{1}[g_i = g] \pi^t(x_i)\hat{m}_n^t(x_i) - \frac{1}{n_{g'}} \sum_{i=1}^n \sum_{t \in \mathcal T} \mathbbm{1}[g_i = g'] \pi^t(x_i)\hat{m}_n^t(x_i) \right).
\end{align*}
We can rewrite the last term as
\begin{align*}
    & \sum_{g \in \mathcal G} \sum_{g' \in \mathcal G, g \neq g'} \lambda(g,g') \left( \frac{1}{n_g} \sum_{i=1}^n \sum_{t \in \mathcal T} \mathbbm{1}[g_i = g] \pi^t(x_i)\hat{m}_n^t(x_i) - \frac{1}{n_{g'}} \sum_{i=1}^n \sum_{t \in \mathcal T} \mathbbm{1}[g_i = g'] \pi^t(x_i)\hat{m}_n^t(x_i) \right)\\
    &= \sum_{g \in \mathcal G} \left(\frac{1}{n_g} \sum_{i=1}^n \sum_{t \in \mathcal T} \mathbbm{1}[g_i = g] \pi^t(x_i)\hat{m}_n^t(x_i) \right) \left( \sum_{g' \in \mathcal G, g \neq g'} (\lambda(g,g') - \lambda(g', g)) \right).
\end{align*}
For ease of notation, we let~$\gamma(g) = \sum_{g' \in \mathcal G, g' \neq g} (\lambda(g,g') - \lambda(g', g))$. We now have
\begin{equation*}
    \begin{aligned}
    L_{\rm out,n}(\pi, \mu, \lambda)  &=&& \frac{1}{n} \sum_{i=1}^n \sum_{t \in \mathcal T} \pi^t(x_i)\left(\hat{m}_n^t(x_i)\left(1 - \frac{n}{n_{g_i}} \gamma(g_i)\right) - \mu^t  \right)+\sum_{t\in \mathcal T} \mu^t  b^t + \sum_{g \in \mathcal G} \sum_{g' \in \mathcal G, g \neq g'} \lambda(g,g') \delta.
    \end{aligned}
\end{equation*}
For any fixed~$\mu$ and~$\lambda$, the inner maximization problem of \eqref{eq:relaxation_sp_out} is thus solved by 
\begin{equation}\label{eq: policy from Lagrangian relaxation_sp out}
    \begin{aligned}
    \pi^{t}(x) = \begin{cases}
    1 &\text{if}\;\;\; t = \min \arg \max_{t' \in \mathcal T}\; \left(\hat{m}_n^{t'}(x)\left(1 - \frac{n}{n_{g}} \gamma(g)\right) - \mu^{t'}\right)\\
    0 &\text{otherwise}
    \end{cases}
    \;\;\; \forall t \in \mathcal T, \;\; x \in \mathcal{X}.
    \end{aligned}
\end{equation}
We suppress the dependence of policy~$\pi$ on~$\mu, \lambda$ notationally and tie-break using a lexicographic tie-breaker. By substituting the above policy into the dual problem \eqref{eq:relaxation_sp_out}, we obtain the convex program
\begin{equation*}\label{eq:bestupperboundproblem_sp_out}
    \begin{aligned}
    \hat{v}_{\rm out,n}^\star &= \min_{\substack{\mu \in \mathbb R_+^{m+1} \\ \lambda \in  {\mathcal L}_\infty(\mathcal G \times \mathcal G, \mathbb R_+) } } \overline{L}_{\rm out,n}(\mu, \lambda),
    \end{aligned}
\end{equation*}
where
\begin{equation*} \label{eq:bestupperboundobj_out}
    \begin{aligned}
    \overline{L}_{\rm out,n}(\mu, \lambda) = &\frac{1}{n} \sum_{i=1}^n \max_{t \in \mathcal T} \left(\hat{m}_n^t(x_i)\left(1 - \frac{n}{n_{g_i}} \gamma(g_i)\right) - \mu^t  \right)+\sum_{t\in \mathcal T} \mu^t  b^t + \sum_{g \in \mathcal G} \sum_{g' \in \mathcal G, g \neq g'} \lambda(g,g') \delta.
    \end{aligned}
\end{equation*}
We denote by~$\hat{\mu}_n^\star, \hat{\lambda}_n^\star$ an optimal solution to this convex program. We define the outcome-parity-constrained policy~$\hat{\pi}^\star_{\rm out,n}$ by substituting~$\hat{\mu}_n^\star, \hat{\lambda}_n^\star$ into \eqref{eq: policy from Lagrangian relaxation_sp out}.}

\section{Empirical Results Details}

\subsection{Linear Outcomes Simulated Data} \label{app:simdata}

\begin{figure*}[!htb]
\begin{center}
\includegraphics[scale=0.43]{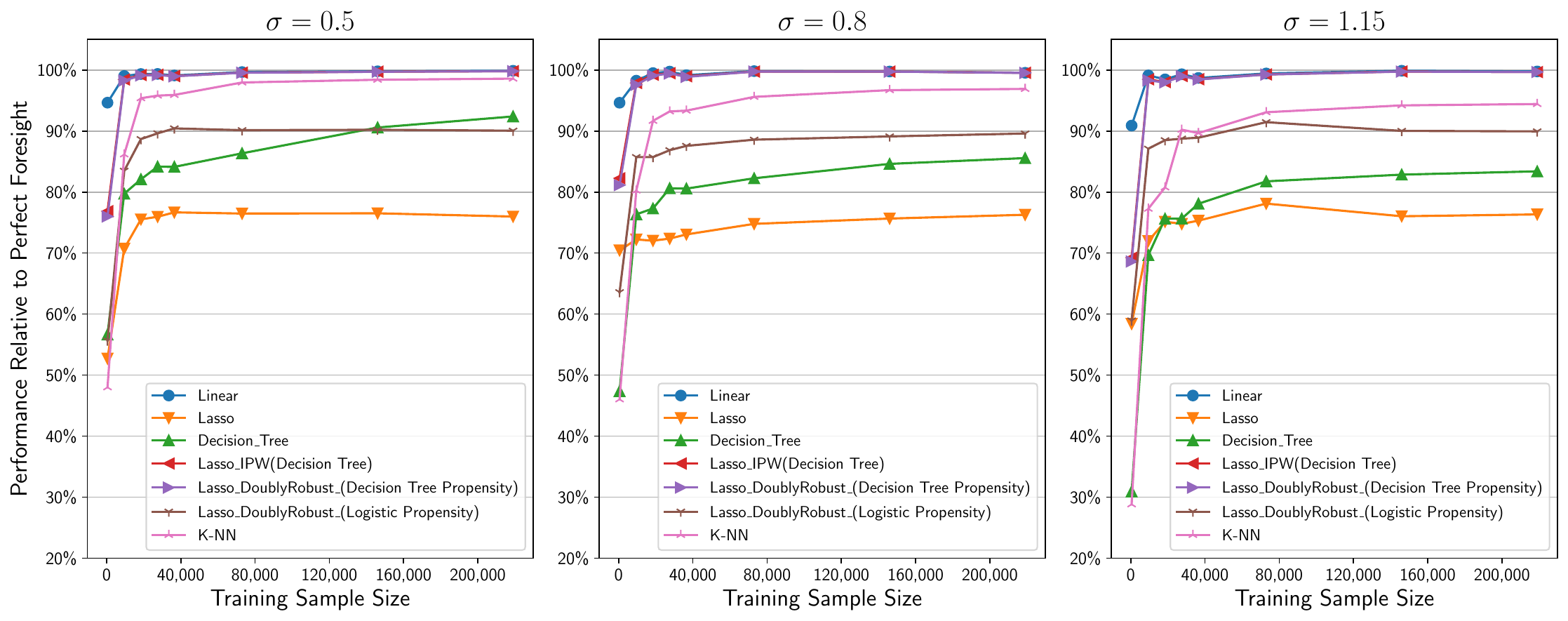}
\caption{Synthetic data results for~$\sigma$ values~$\{0.5, 0.8, 1.15\}$. All subfigures show the ratio of out-of-sample performance of the sample-based queuing policy to that of the perfect foresight policy in dependence of training set size. Each line corresponds to a policy using a different outcome estimator.}
\label{fig:appendix_simulation}
\end{center}
\end{figure*}

In Figure~\ref{fig:appendix_simulation}, we show additional simulation results for values of~$\sigma$ from the set~$\{0.5, 0.8, 1.15\}$ for the noise term. The qualitative takeaways are the same as those in Section~\ref{subsec: sim results} and serve primarily to show the degradation of non-parametric methods of decision tree and KNN from added noise in the data.

\begin{figure*}[!htb]
\begin{center}
\includegraphics[scale=0.43]{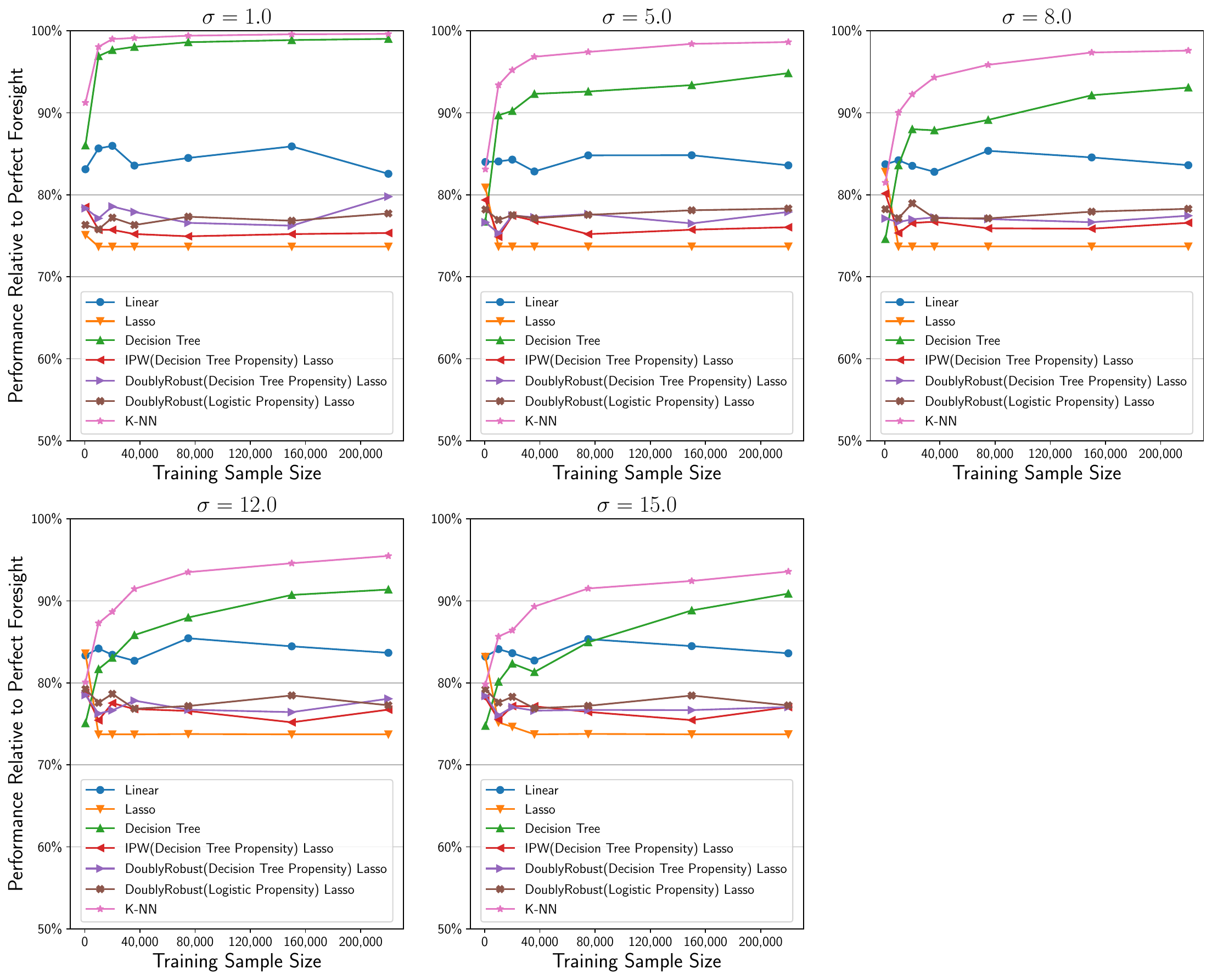}
\caption{Synthetic data results for nonlinear outcome function for~$\sigma$ values~$\{1, 5, 8, 12, 15\}$. All subfigures show the ratio of out-of-sample performance of the sample-based queuing policy to that of the perfect foresight policy in dependence of training set size. Each line corresponds to a policy using a different outcome estimator.}
\label{fig:appendix_nonlinear_sim_res}
\end{center}
\end{figure*}

{\color{blue}
\subsection{Nonlinear Outcomes Simulated Data} \label{app:nonlinear experiments}

In this section, we investigate a variant of the setup from Section~\ref{sec:synthetic_data_generation}, where the outcomes are now nonlinear in the covariates. In particular, the potential outcomes $Y^t$ and conditional mean treatment outcomes $m^t$, for each $t \in \{0, 1, 2\}$, can be written as:
\begin{equation*}\label{eq: sim_mean_outcome_nonlinear} 
    \begin{aligned}
    &Y^0 = m^0(X) + \epsilon_0, &m^0(X) = 0.25 X_1^2 + 0.75 X_2^2\\
    &Y^1 = m^1(X) + \epsilon_1, &m^1(X) = 0.75 X_1^2 + 0.75 X_2^2\\
    &Y^2 = m^2(X) + \epsilon_2, &m^2(X) = 0.25 X_1^2 + 1.25 X_2^2\\
    \end{aligned}
\end{equation*}
where $\epsilon_t \sim \mathcal{N}(0, \sigma)$, $t \in \mathcal{T}$, are i.i.d. noise terms. We now vary $\sigma$ over a larger range of values $\{1, 5, 8, 12, 15\}$ since the quadratic outcome functions $m^0, m^1, m^2$ now take on a larger range of values (using the range of $\sigma$ values from the linear case of Section~\ref{sec:synthetic_data_generation} does not yield noise terms with meaningful differences compared to the scale of outcome values).

Figure~\ref{fig:appendix_nonlinear_sim_res} shows the simulation results for this setting. Compared to the case of linear conditional mean treatment outcomes, as one would expect, no policy based on a linear outcome estimator reaches close to $100\%$ of the perfect performance, regardless of sample size, since they are all misspecified now. In the low noise regime ($\sigma = 1$), the policy based on the non-parametric $k$-NN estimator  comes close to $100\%$ (as does using decision trees), but increasing levels of noise cause a deterioration in the performance of $k$-NN. We also see that simpler parametric models like linear and lasso, though misspecified, are relatively robust to increasing levels of noise. The causal adjusted versions of lasso also improve their performance as well here, though only to a small degree.} 

\subsection{HMIS and VI-SPDAT Data} \label{app:realdata}

\subsubsection{Data Details and Preparation}
In this subsection, we provide further details on the current system, raw data, and how we defined outcomes for our experimental results. 

\begin{figure*}[!htb]
\begin{center}
\includegraphics[scale=0.6]{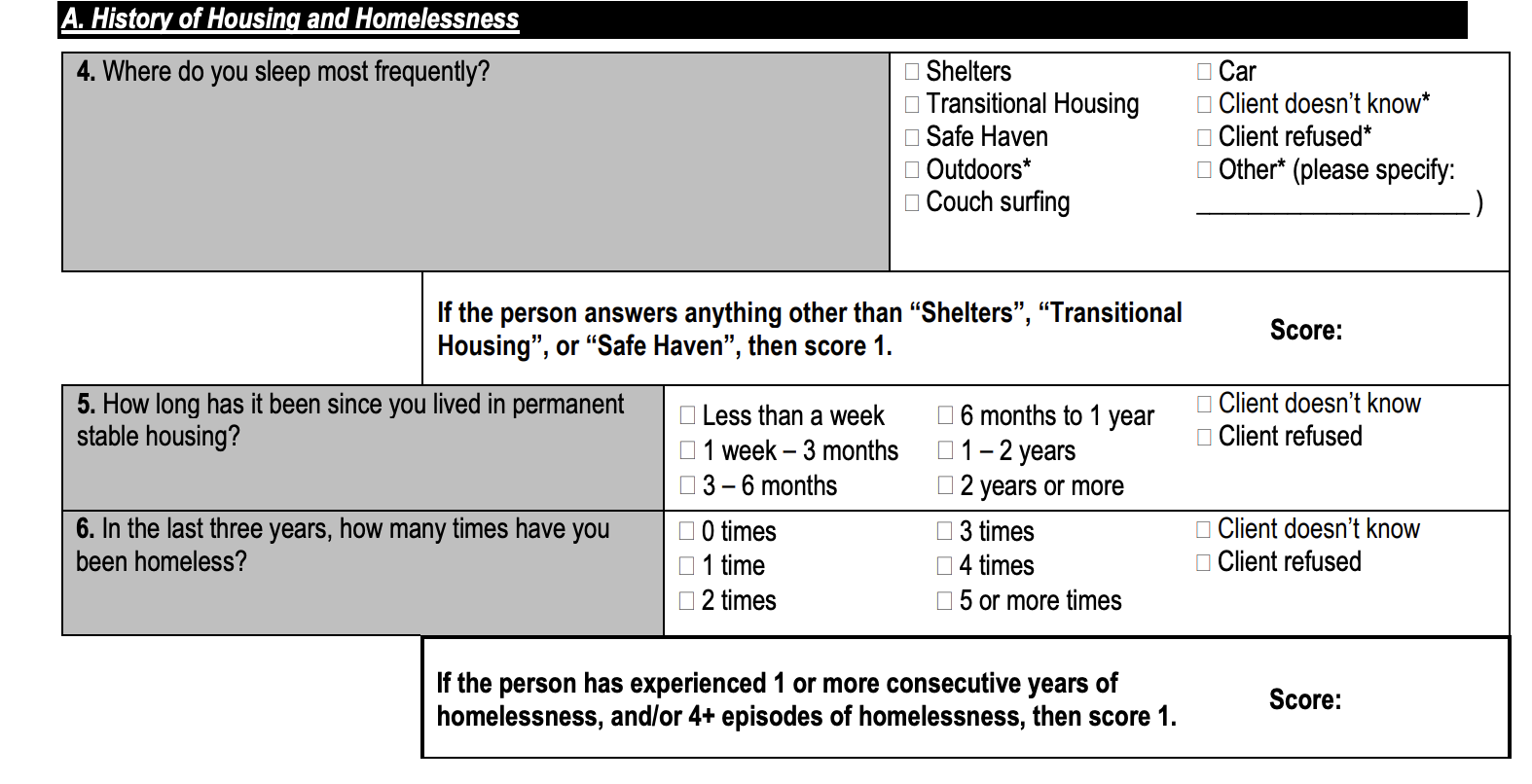}
\caption{VI-SPDAT sample}
\label{fig:vi-spdat}
\end{center}
\end{figure*}

When an individual arrives to the system, they are assessed for vulnerability with the VI-SPDAT tool, which contains a series of questions such as prior length of and number of experiences of homelessness. Adverse answers, such as sleeping situations indicative of homelessness, are given a weight of 1 so that an individual's vulnerability score ranges from 0 to 17, where a higher score indicates more vulnerability. The question and answers from the tool are used as features in our estimation methods and a sample of the tool is shown in Figure~\ref{fig:vi-spdat}~\citepalias{lahsa-vi-spdat}. To determine treatment assignments and outcomes after individuals receive their VI-SPDAT assessments, we can match them with their enrollment history from the raw HMIS data, which details their interactions with the system. Full details on the HMIS data and variables can be found in the HMIS Data Standards Data Dictionary~\citepalias{hmis_data_dict}. For a given individual, we observe a sequence of enrollments, each characterized by a type of service and the corresponding date information. 

Finally note that there is not a one-to-one correspondence between individuals found in VI-SPDAT assessment data and the HMIS data. After merging the two datasets, and removing missing data issues, we obtain our final dataset that consists of~$63,764$ individuals assessed between 1/12/2015 and 12/31/2019, their assessment information, treatment received, and observed outcome. Recall that we use all individuals assessed between 1/12/2015 to 12/31/2017 as our training set to learn treatment outcomes and construct our policy, and then evaluate on individuals assessed between 1/1/2018 to 12/31/2019 to measure out-of-sample performance.

\textit{Historical Treatment Assignment Definition.} Since the original HMIS data does not explicitly contain ``Treatment'', we use the enrollment history to create these variables. In Table~\ref{table:HMIS Enrollment} below, we breakout the~$11$ types of enrollments found within HMIS data, where the enrollments `PH-PSH' and `PH-RRH' are considered permanent interventions or treatments, and the remaining enrollments are considered no-treatment, or `Services Only', either due to their temporary nature, shelter type, or nature of service provided. 
We also further breakout `PH-PSH' into `tenant' and `site' based on an additional column, where the two sub-types differ, broadly speaking, in terms of supportive services available on-site and process of actually receiving housing. This categorization of the enrollment types was chosen based on discussions with matchers working within the LA CES system. 

\begin{table}[h!]
\caption{\color{blue} HMIS enrollment types grouped by enrollments considered as permanent interventions and no treatments. Percentage breakdowns come from our data used for estimation and are based on individuals who received a VI-SPDAT assessment.
}
\vspace*{2mm}
\centering
\begin{tabular}{ |s|s|c| }
\hline
 \rowcolor{tablecell} \textbf{Treatment Type} & \textbf{Enrollment Type} & \textbf{Proportion of All Enrollments} \\ 
\hline \hline
 {Permanent Intervention} & {{PH - RRH}} &~$14.16\%$ \\  
 \hline
 {Permanent Intervention} & {{PH - PSH Site}} &~$4.63\%$ \\  
 \hline
 {Permanent Intervention} & {{PH - PSH Tenant}} &~$4.14\%$ \\  
 \hline
 {No Treatment} &{{Emergency Shelter}} &~$26.74\%$ \\  
 \hline
 {No Treatment} & {{Street Outreach}}&~$26.23\%$ \\  
 \hline
 {No Treatment} & {Services} &~$18.73\%$ \\  
 \hline
 {No Treatment} & {Other} &~$3.10\%$ \\  
 \hline  
 {No Treatment} & {Transitional Housing} &~$1.31\%$ \\  
 \hline 
 {No Treatment} & {Day Shelter} &~$0.55\%$ \\  
 \hline 
 {No Treatment} & {Homelessness Prevention} &~$0.30\%$ \\  
 \hline 
 {No Treatment} & {PH - Other} &~$0.09\%$ \\  
 \hline 
 {No Treatment} & {{Safe Haven}} &~$0.05\%$ \\  
 \hline
\end{tabular}

\label{table:HMIS Enrollment}
\end{table}

\textit{Outcome Definition.} Recall from Section~\ref{subsubsec: hmis data} that we focus on measuring returns to homelessness after an intervention as the outcome variable, in line with the goals listed by HUD. Since such an outcome is not explicitly tracked within the data, we construct a proxy variable as follows. To determine an individual's outcome, we first define an two-year observation window depending on their treatment received. For those receiving `PH-RRH', we also observe a `Move-In Date' column signifying when an individual was recorded to have officially received the resource. We define the `PH-RRH' window to start from the `Move-In Date'. Similarly for `PH-PSH', we use the `Move-In Date' to determine the starting point, but add a 100 day lag to the `Move-In Date' because we concluded the recorded `Move-In Date' within HMIS is likely not accurate. Through discussions with matchers within the LA CES system, we learned that the `Move-In Date' within HMIS is an initial date inputted by case managers, but often does not get updated after an individual eventually moves in. Furthermore, we analyzed an auxiliary dataset, the Resource Management System (RMS) data, that supports this conclusion. For a subset of PSH units found within the HMIS data, the RMS data contains specific information such as their eligibility requirements and the updated move-in dates. For units found within both the HMIS and RMS data, we found discrepancies between the reported move-in dates, where the RMS `Move-In Date' occurs often months after the corresponding HMIS `Move-In Date', indicating that the HMIS `Move-In Date' column for `PH-PSH' have inaccuracies. This presents a problem where we often found individuals receiving PSH within the HMIS data to subsequently enroll in enrollments indicating an individual was experiencing homelessness (`Street Outreach', `Emergency Shelter', and `Safe Haven') shortly after their reported HMIS PSH 'Move-In Date'. If the HMIS `Move-In Date' was incorrect and those individuals had actually not yet moved-in yet to their PSH units, then we would incorrectly conclude they experienced an adverse outcome after receiving PSH. 

\begin{figure*}[!htb]
\begin{center}
\includegraphics[scale=0.45]{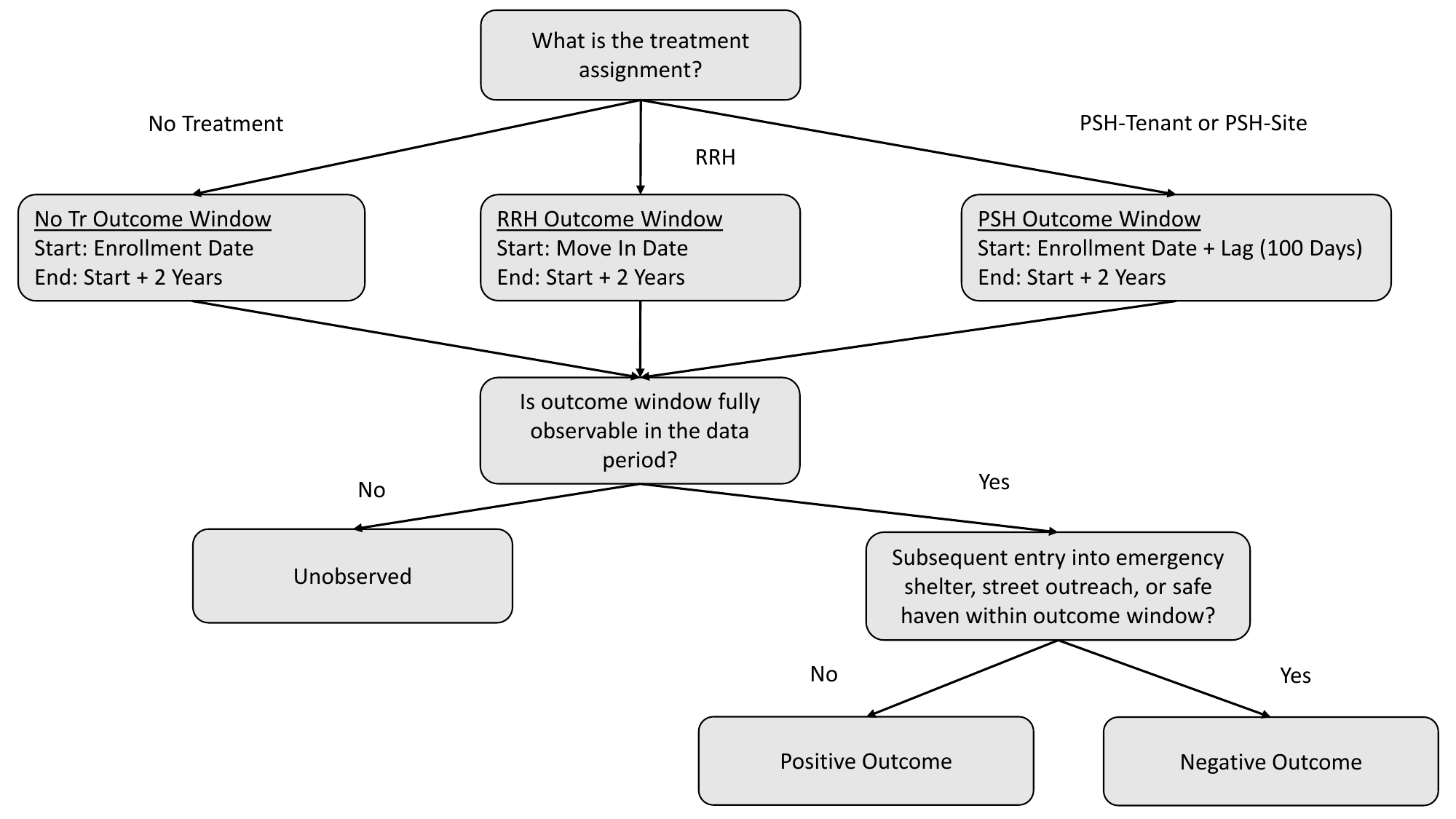}
\caption{Construction of treatment outcomes.}
\label{fig:outcome def chart}
\end{center}
\end{figure*}

To address the problems with the HMIS PSH `Move-In Date' column, we add a 100 day lag to the `Move-In Date' to determine the starting point for the `PH-PSH' observation window. We determined this 100 day lag based on the distribution of positive differences between the RMS and HMIS `Move-In Date'. Finally for those receiving \texttt{no-treatment}, we consider their first enrollment date of any type to be the start of their observation window. Our choice of a two-year observation window for determining if an individual returned to homelessness after an intervention involved a trade-off where shorter windows resulted in insufficient time to observe post-intervention outcomes while longer windows resulted in fewer individuals for whom we can observe their full outcome window. If we observe any occurrence of the three proxies of homelessness enrollments (`Street Outreach', `Emergency Shelter', or `Safe Haven') within that observation window, then we consider an individual to have experienced a negative outcome. Otherwise we consider them to have a positive outcome. We made the above choices of treatment types, window length, start dates, and more considering the trade-offs in terms of number of samples, `bias' of the proxies, and data quality. For individuals for whom we do not fully observe their outcome window, i.e., the outcome window ends after the last observation date in our data, we remove them from our analysis. Keeping those indivisuals would bias the data towards positive outcomes since individuals may experience negative outcomes after the observable time in our data. For a summary of our outcome column construction, see Figure~\ref{fig:outcome def chart}.

\subsubsection{Generating Counterfactuals and Fitting Training Set Estimators} \label{app: realdata_counterfactuals}
To obtain counterfactuals for evaluating performance within the held out test set, we took a semi-synthetic approach and used model-generated outcomes based on the entire dataset. {\color{blue}We consider three causal estimation approaches for generating counterfactual outcomes, the direct (or T-Learner method), IPW, and DR methods. While the direct method simply fits an observed outcome estimator for each treatment group, the IPW and DR methods involve fitting an additional propensity model that estimates the probability of receiving each treatment conditional on observed covariates. The inverse propensity scores are then used as sample weights in a second stage counterfactual estimator in order to correct for selection bias. Since many causal estimation methods involve a two step procedure where the first step estimates nuisance components like outcome and propensity models and the second step combines these components into a final counterfactual model, recent works often refer to the final model as a ``final'' or ``second'' stage estimator~\citep{Kunzel2019}. Therefore, we also adopt this terminology going forward and refer to the final counterfactual estimator as a second stage estimator. We largely follow similar procedures as those described in \cite{econml} and \cite{jacob2021cate}. For counterfactual outcome model selection for each treatment,} we want our model output probabilities to be as well-calibrated to the observed outcome distribution as possible. Since our outcome is binary, the probability of positive outcome for individual with features~$X$ under treatment~$t$ is exactly the function~$m^t(X)$. Therefore, requiring calibrated probabilities under each treatment is the same as requiring~$\hat{m}^t(x) \sim m^t(x)$ for all~$x \in \mathcal{X}$ for our chosen treatment outcome generating model~$\hat{m}^t$. Therefore, we generate counterfactuals from models that empirically match the outcome distribution of each treatment as well as possible. To ensure generalization as well, we take a standard cross-validation approach to select the best model for each treatment group. Finally, we noticed individuals without any indication of a disability in the data were very unlikely to receive \texttt{PSH} and, generally speaking, \texttt{PSH} units required some form of disability to be eligible. Therefore modeling the \texttt{PSH} outcomes of non-disabled individuals would possibly violate the positivity assumption in Assumption~\ref{ass: causal} and practically we have insufficient samples. Therefore, we model disabled and non-disabled individuals separately and did not model non-disabled individual outcomes under \texttt{PSH} treatment. 

Note that the procedure outlined in this section was applied to the entire dataset and does not relate to the time-based splitting of training (1/12/2015 to 12/31/2017 assessments) and testing sets (1/1/2018 to 12/31/2019 assessments) used to evaluate performance in our experiment. Here, we are simply trying to find models that empirically fit our data and outcome distributions well enough to serve as counterfactuals.

\textit{Cross Validation.}
Since we want to have well-calibrated probabilistic predictions, we use log-loss during the cross validation phase as the evaluation metric for tuning hyper-parameters of {\color{blue}each model class for estimating observed outcomes, propensities, and any second stage estimator.} We choose to minimize log-loss since log-likelihood is a strictly proper scoring rule, i.e., optimizing this metric will yield well-calibrated probability predictions~\citep{strictlyproperscoringrule}. For {\color{blue}estimating each treatment group's outcomes, treatment assignment propensities, and second stage components, we use 5-fold cross validation and select the hyper-parameters with the best average log-loss across all folds for each model class we considered (logistic regression with and without regularization, random forests, gradient boosted trees, and XGBoost trees).}

\begin{figure*}
\begin{center}
\includegraphics[width=\textwidth]{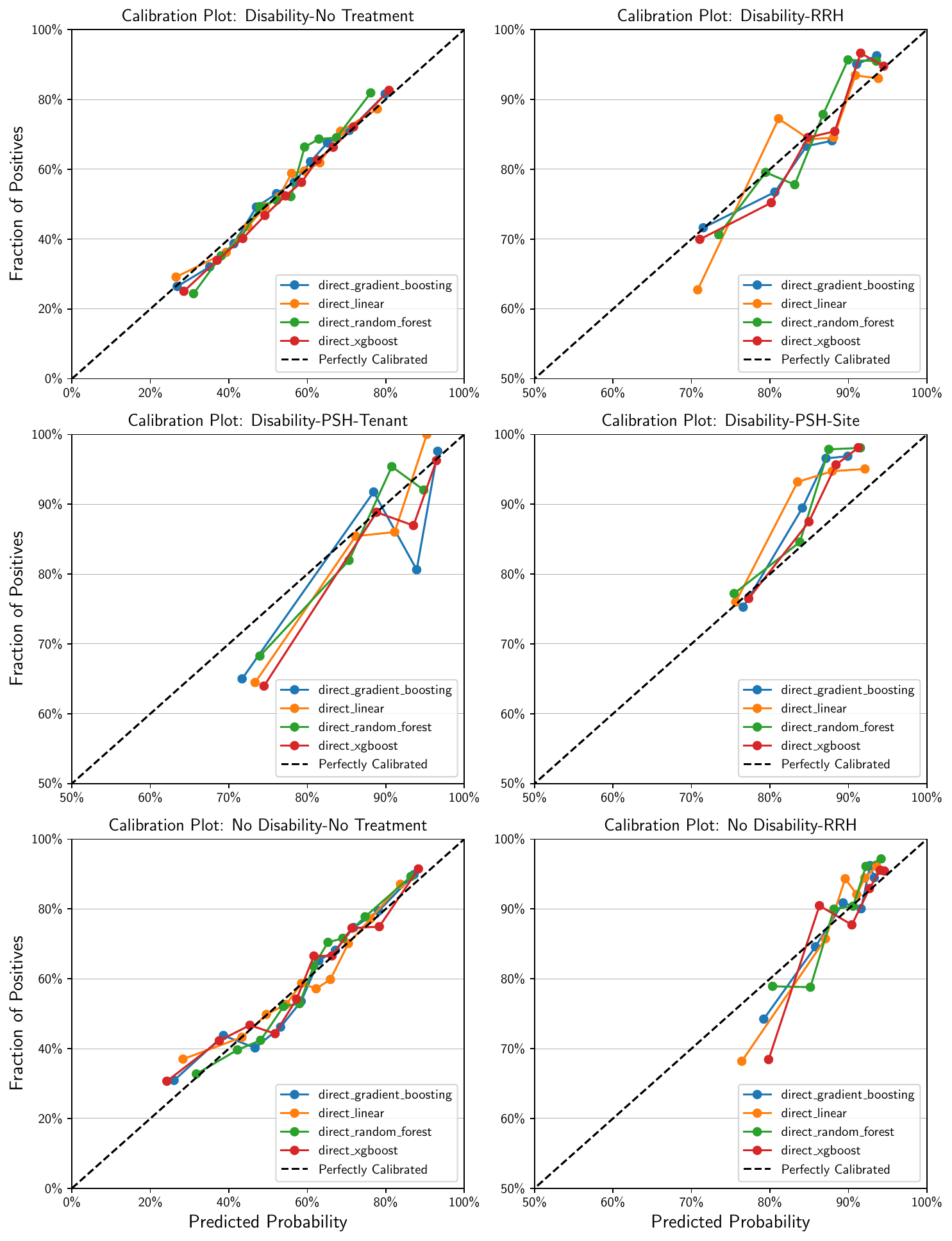}
\caption{Calibration plots on validation set for each direct method outcome model class and different combinations of disability and treatment group.}
\label{fig: validation set calibration direct}
\end{center}
\end{figure*}

\begin{figure*}
\begin{center}
\includegraphics[width=\textwidth]{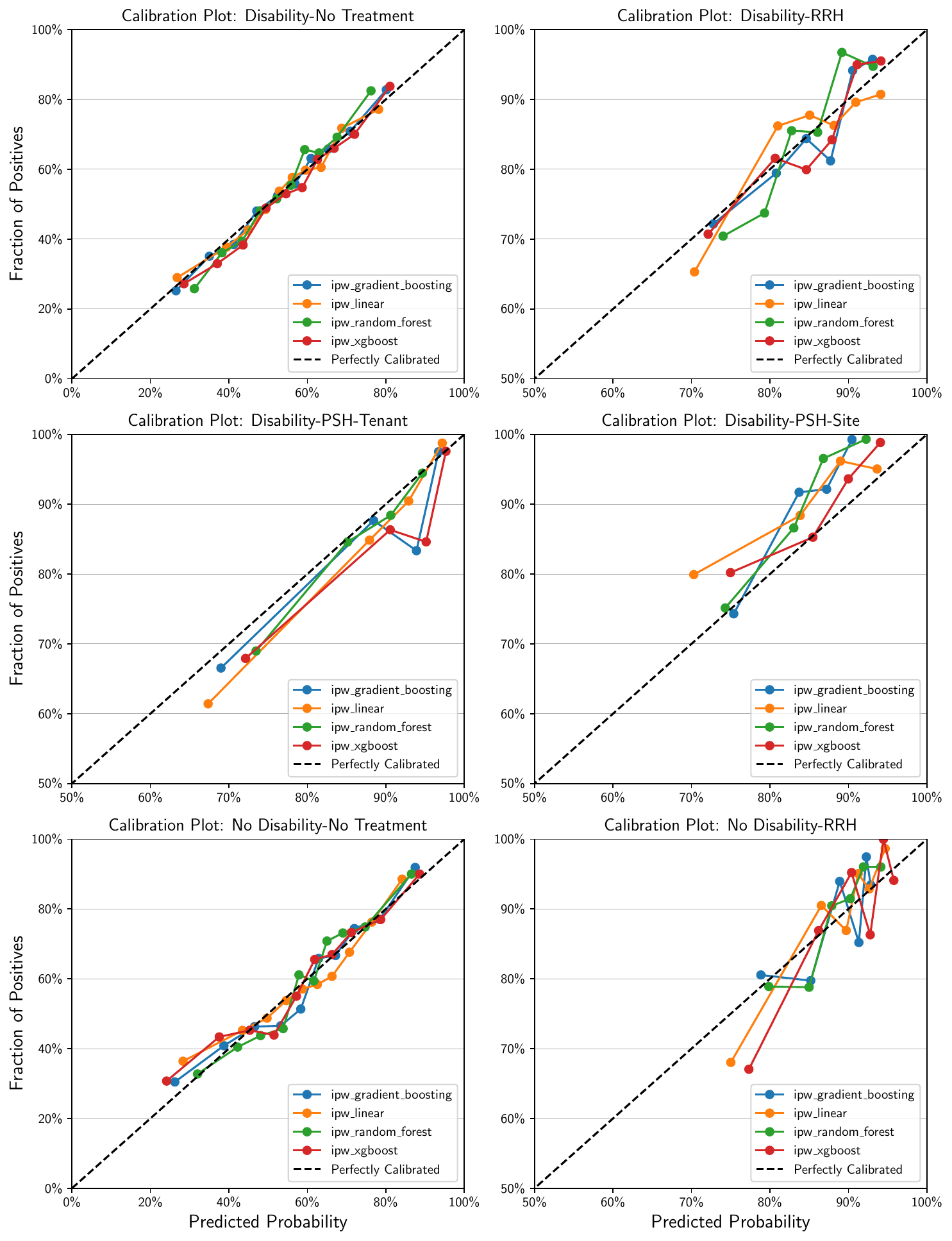}
\caption{Calibration plots on validation set for each IPW method outcome model class and different combinations of disability and treatment group.}
\label{fig: validation set calibration ipw}
\end{center}
\end{figure*}

\begin{figure*}
\begin{center}
\includegraphics[width=\textwidth]{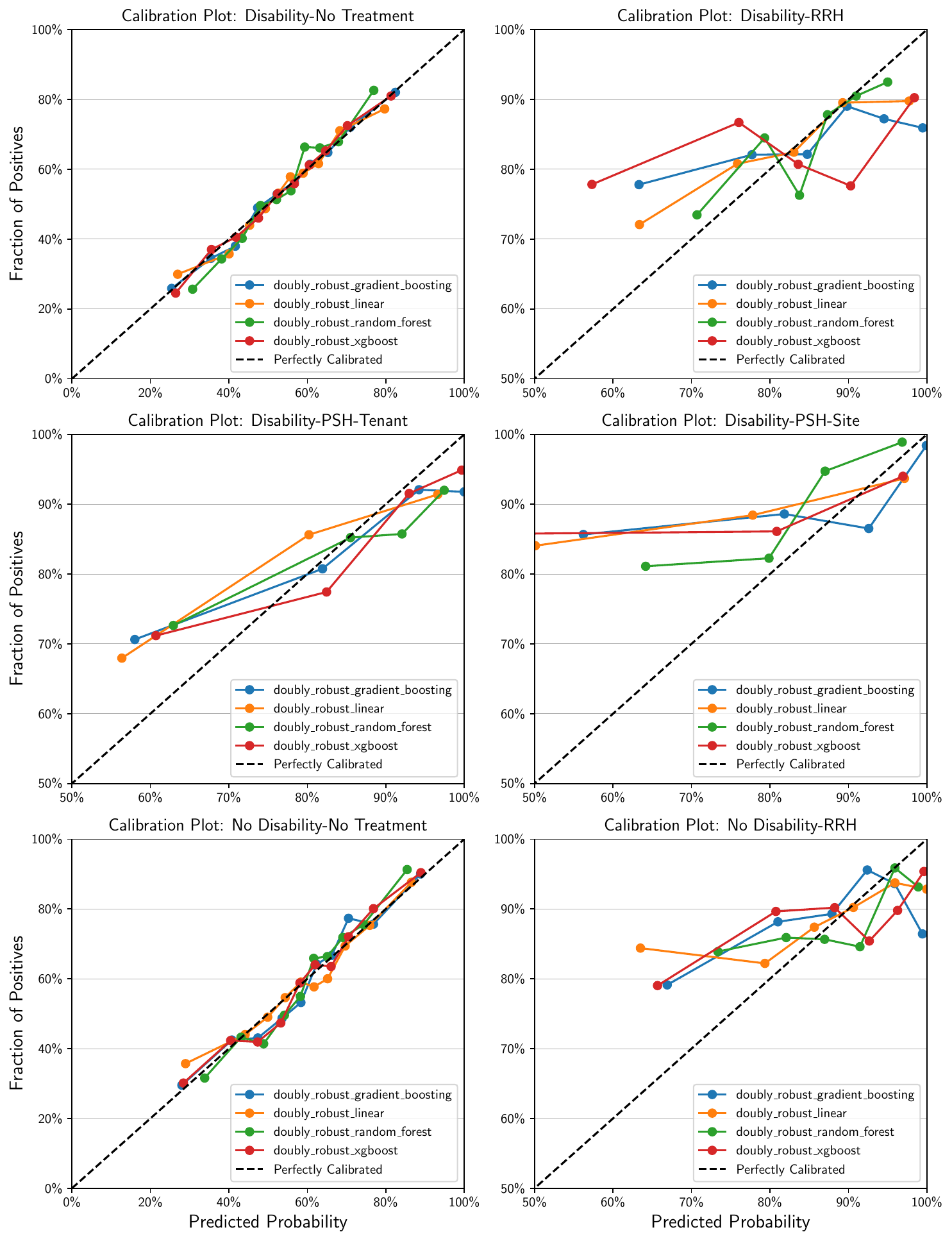}
\caption{Calibration plots on validation set for each Doubly Robust method outcome model class and different combinations of disability and treatment group.}
\label{fig: validation set calibration doubly robust}
\end{center}
\end{figure*}

\textit{Validation Set of Calibration.}
For each treatment group and {\color{blue}hyper-parameter tuned counterfactual outcome model,} we generate a calibration plot to visually pick the model that is best calibrated. {\color{blue}In Figures~\ref{fig: validation set calibration direct}, \ref{fig: validation set calibration ipw}, \ref{fig: validation set calibration doubly robust}, we plot the calibration curves from a held-out validation set for the direct, IPW, and DR methods, respectively. For example, in each of the subplots of Figure~\ref{fig: validation set calibration direct}}, the dotted diagonal line represents probability outputs of a perfectly calibrated model while the colored solid lines with markers are the calibration curves of each model class. {\color{blue}In particular, we calculate weighted calibration curves. To generate each calibration curve, we split the probability predictions for a given model into equal sized buckets and find the average predicted probability (x-axis) and plot it against the weighted average proportion of positive outcomes (y-axis), where the weights are the inverse propensity weights. In essence, because we can only plot and observe outcomes for those who received a particular treatment, we need to correct for observational data selection bias. One way to correct for selection bias is to calculate each bucket's average outcome using the IPW estimator \citep{hernan2023causal}. 
} 
Also, because the treatment groups \texttt{RRH}, \texttt{PSH}-Tenant, and \texttt{PSH}-Site had much smaller numbers of samples relative to the \texttt{no-treatment} group, we used a smaller number of buckets to evaluate calibration to avoid each bucket being too small in size and resulting in noise. Generally speaking, {\color{blue}there was not a big difference between the calibration curves of the direct and IPW methods, while the DR curves on the smaller treatment groups look more miscalibrated. The DR miscalibration may be due to the weighted residuals used in the doubly robust second stage estimation, as the residuals, which are functions of a binary outcome subtracting a predicted probability, may be better suited for debiasing continuous outcomes. Given that overall the best IPW calibration curve within each disability and treatment group appeared the most calibrated, we defaulted to using the IPW models for generating counterfactual outcomes.}

\begin{figure*}
\begin{center}
\includegraphics[width=\textwidth]{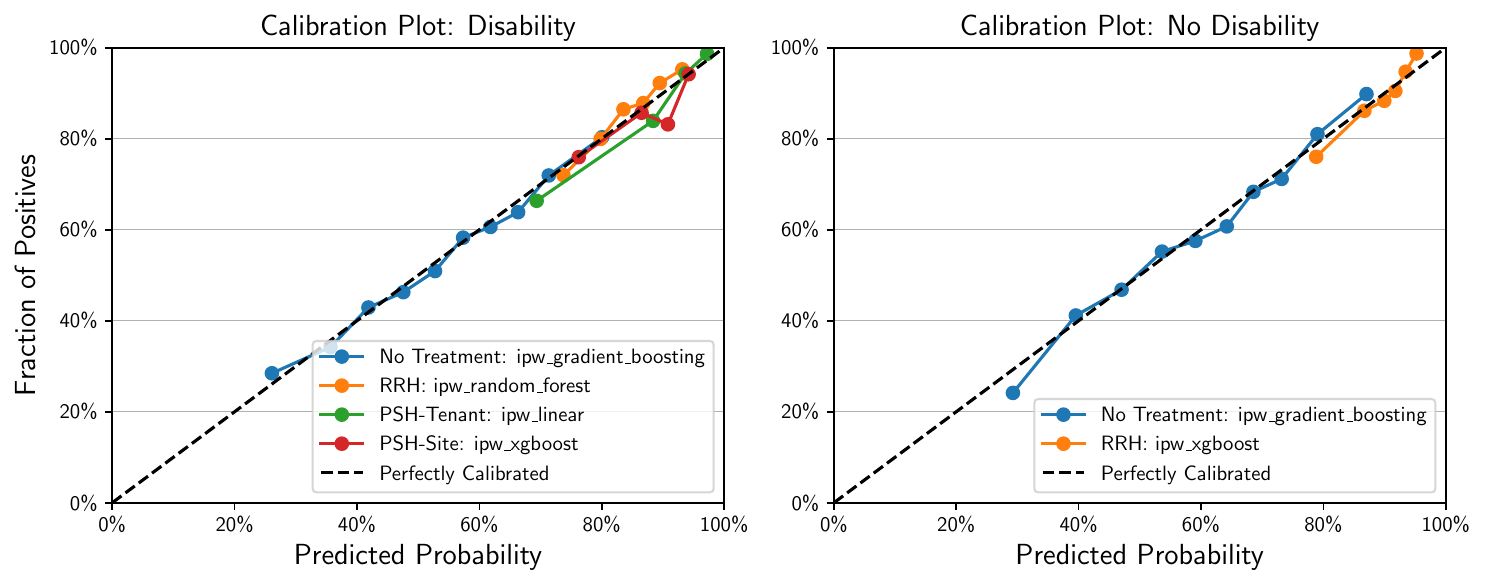}
\caption{Out-of-sample calibration plots of final counterfactual generating models chosen for each treatment group and for individuals with and without disabilities.}
\label{fig: test set calibration}
\end{center}
\end{figure*}

After selecting the appropriate estimator for each treatment group of disability and non-disability groups, we refitted the chosen models on the training and validation sets and checked for out-of-sample generalization on the held-out test set.  In Figure~\ref{fig: test set calibration}. we plot the final test-set calibration curves for each treatment group for individuals with and without disabilities. In general, the final chosen models for each treatment group appear well-calibrated visually with some slight deviations from the diagonal line. Finally, while the above procedure describes how we chose the counterfactual generating models, we still need to generate predictions of treatment outcomes for the training set only (all individuals assessed 1/12/2015 to 12/31/2017) to learn a sample-based dual-price queuing policy. We repeat the above procedure but only on the training data to generate~$\hat{m}^t$. 
\subsubsection{Fairness of Estimates} \label{app: realdata_counterfactuals_fairness}
\begin{figure*}
\begin{center}
\includegraphics[width=\textwidth]{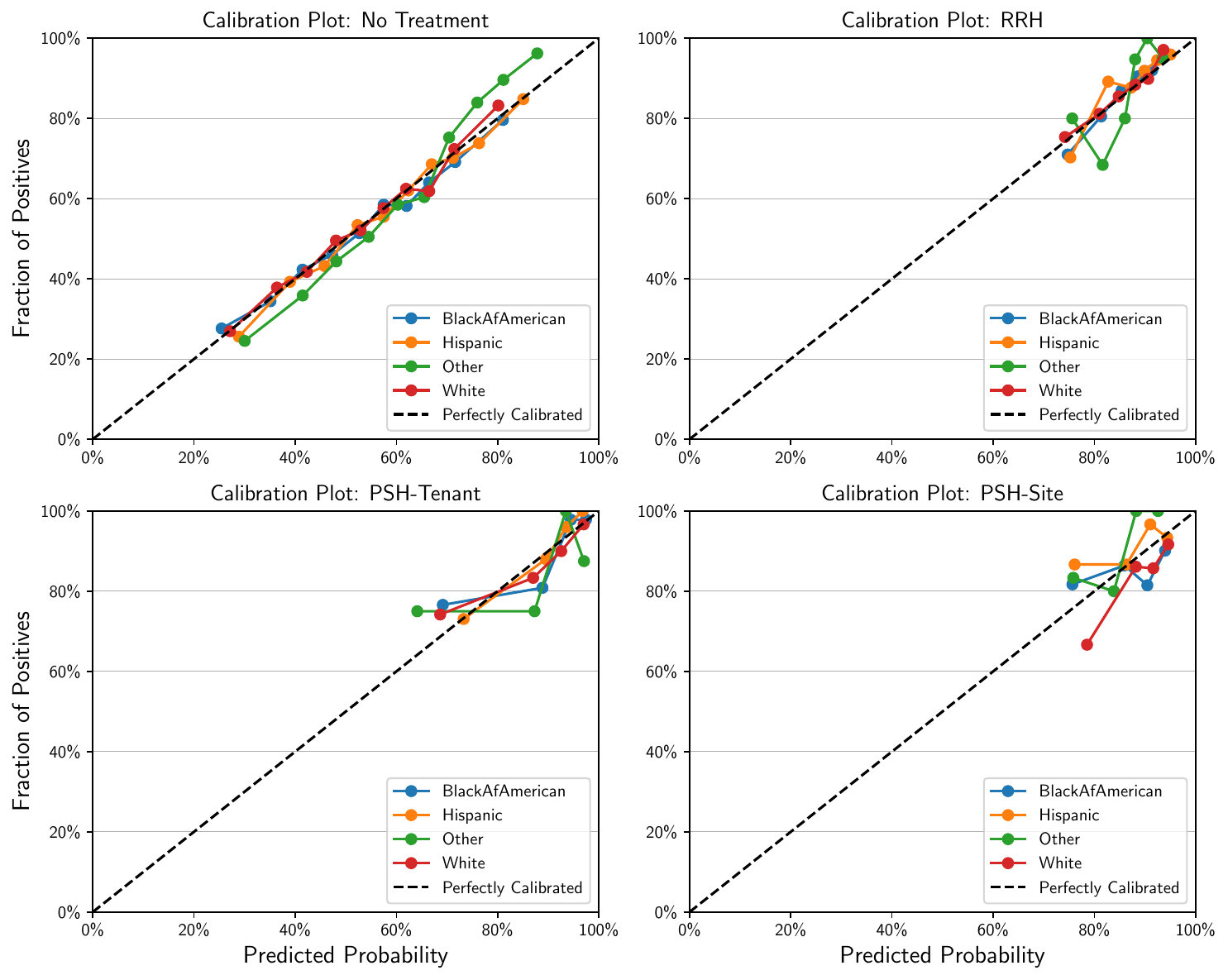}
\caption{Out-of-sample calibration plots for different racial groups and treatments from final counterfactual generating models chosen.}
\label{fig: calib_race}
\end{center}
\end{figure*}

An additional concern with our counterfactual models for evaluation was potential bias by race. Since we want to use our generated counterfactuals as reflections of the true generating process, any calibration bias would mean our evaluation results would be potentially biased by race. In Figure~\ref{fig: calib_race}, we plot and compare calibration curves for each racial group across treatments on the test-set for the final chosen models. For the most part, the calibration curves across racial groups are similar for each treatment, with the exception of the `Other' racial group. This is likely because the `Other' group is only a small portion of the sample, representing only~$\sim 7\%$ of the sample. Furthermore, the scarce PSH-Tenant and PSH-Site treatment groups also have worse calibration than RRH, again likely due to insufficient samples from this small treatment group.  
\subsubsection{Full Experimental Test Results}\label{app: full_experiment}
\begin{figure*}
\begin{center}
\includegraphics[width=\textwidth]{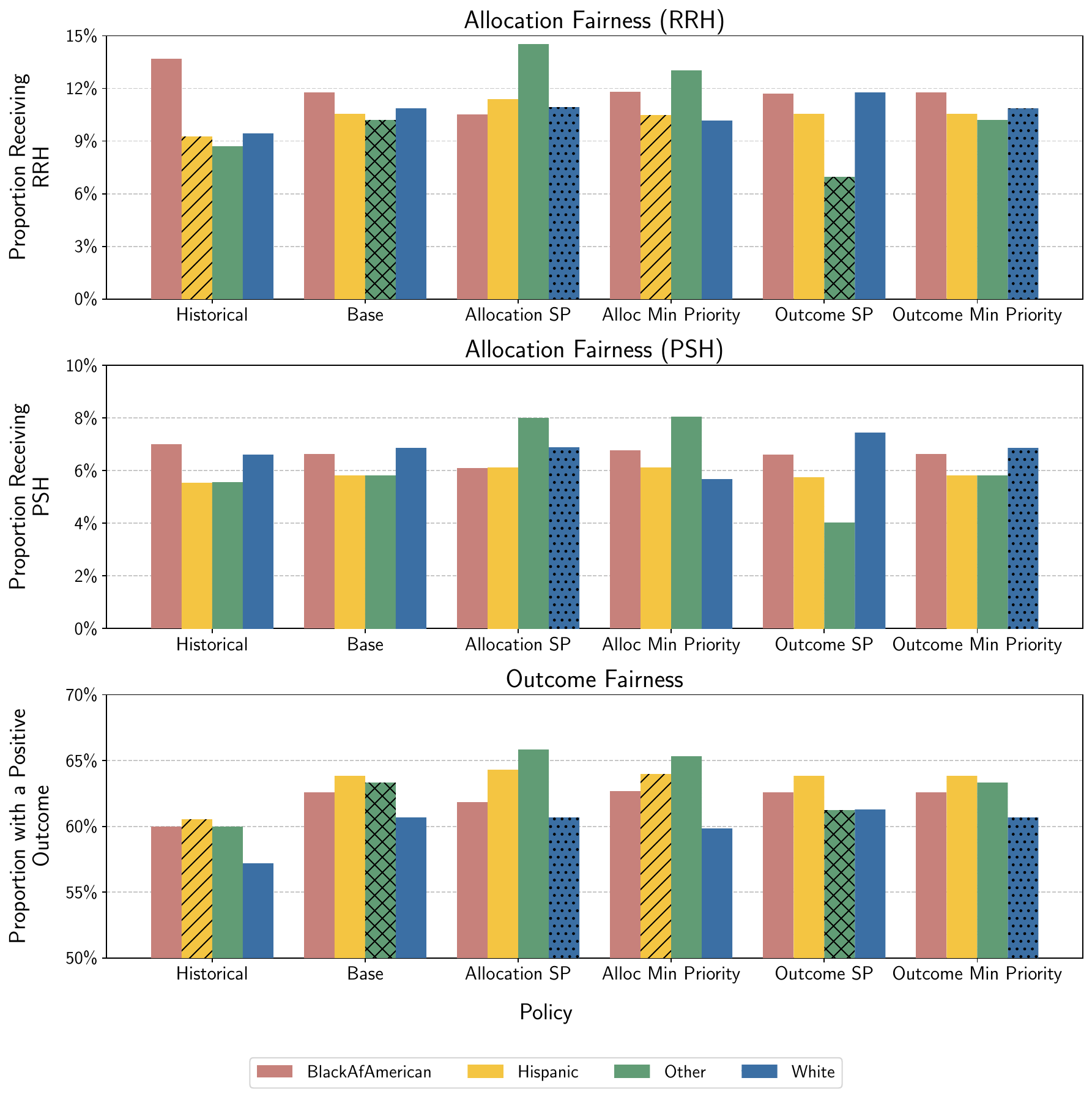}
\caption{Out-of-sample fairness results showing proportion of each racial group receiving \texttt{RRH}, \texttt{PSH}
(Allocation Fairness) and the proportion of positive outcomes for each racial group (Outcome Fairness) under different policies.}
\label{fig:alloc_fair_appendix}
\end{center}
\end{figure*}

In this subsection, we present further results on two additional fairness-constrained policies \texttt{Allocation SP} and \texttt{Outcome SP} policies, which seek to achieve statistical parity in allocation (constraint~\eqref{eq: sp alloc}) and statistical parity in outcomes (constraint~\eqref{eq: sp_outcome_constraint}) by race, respectively. The out-of-sample performance in terms of expected outcomes are summarized in Table~\ref{table:real_overall_outcomes_appendix} and allocation and outcome fairness properties are presented in Figure~\ref{fig:alloc_fair_appendix}. We use a~$\delta$ of~$0.01$ for the fairness constraints in the sample-based problem to learn the \texttt{Allocation SP} and \texttt{Outcome SP} policies.

\begin{table}[!htb]
\caption{Out-of-sample average expected outcome under each policy, and the relative percentage difference compared to the \texttt{Historical} allocation.}
\vspace*{2mm}
\centering
\begin{tabular}{ |s|c|c| }
\hline
 \rowcolor{tablecell} \textbf{Policy} & \textbf{\begin{tabular}{c}Proportion of \\ Positive Outcomes \end{tabular} } & \textbf{\begin{tabular}{c}Percentage Change\\ vs Historical \end{tabular} }  \\ 
\hline \hline
 {\texttt{No Treatment}} &~$54.25\%$ &~$-8.76\%$\\ 
 \hline
 {\texttt{Historical}} &~$59.46\%$ &~$0.00\%$\\ 
 \hline
 {\texttt{Base}} &~$62.53\%$ &~$5.16\%$\\  
 \hline
 {\texttt{Allocation SP}} &~$62.50\%$ &~$5.11\%$\\  
 \hline
 {\texttt{Alloc Min Priority}} &~$62.51\%$ &~$5.13\%$\\  
 \hline
 {\texttt{Outcome SP}} &~$62.52\%$ &~$5.15\%$\\  
 \hline
 {\texttt{Outcome Min Priority}} &~$62.53\%$ &~$5.16\%$\\  
 \hline
 {\texttt{Perfect Foresight}} &~$63.45\%$ &~$6.71\%$\\  
 \hline 
\end{tabular}
\label{table:real_overall_outcomes_appendix}
\end{table}

{\color{blue}For \texttt{RRH}, we see that \texttt{Allocation SP} has the highest minimum percentage receiving \texttt{RRH} across all racial groups compared to all other policies, but surprisingly does not result in the smallest differences in percentages between all racial groups out-of-sample. It also allocates fewer resources to Black individuals compared to \texttt{Historical} or any other policy, which may be undesirable when policymakers want to prioritize minority groups. In this instance, policymakers need to choose between statistical parity in allocation or prioritizing minority groups for allocation of \texttt{RRH} but our proposed policy can incorporate either notions into consideration or possibly some combination of the two. Similarly, for \texttt{PSH}, \texttt{Allocation SP} also has the highest minimum allocation rate across all racial groups, but again does not have the smallest differences between all racial groups and allocates to Black individuals at the lowest rate relative to all other policies.} Also, of note in the allocation fairness results is that the \texttt{Outcome SP} policy allocates \texttt{RRH} and \texttt{PSH} to the \texttt{Other} group at far lower rates compared to other racial groups. This is because in our predictions, individuals of the \texttt{Other} group had predicted outcomes under \texttt{no-treatment} that skewed more positive relative to other racial groups. Therefore, the \texttt{Outcome SP} policy allocated less to the \texttt{Other} group in an attempt to achieve outcome parity across all groups because their predicted outcomes under \texttt{no-treatment} were already higher compared to the other racial groups. In terms of statistical parity in outcomes, we see that \texttt{Outcome SP}, in addition to \texttt{Base} and \texttt{Outcome Min Priority}, resulted in the smallest differences in proportion of positive outcomes between all racial groups out-of-sample. 

Finally, we see in Table~\ref{table:real_overall_outcomes_appendix} that our fairness-constrained policies \texttt{Allocation SP} and \texttt{Outcome SP} suffer almost no performance gap in terms of expected outcomes compared to \texttt{Base} while potentially improving the fairness properties of the \texttt{Base} policy. Again, these further results suggest almost no `price of fairness'. 

{\color{blue}
\subsection{Qualitative Behavior of Dual-Price Queuing Policy}\label{subsec: qualitative_behavior}

\subsubsection{Allocations with Respect to Vulnerability}\label{subsec: alloc_vul}

One potential concern for our policy is that by optimizing for overall expected outcomes, it may overlook vulnerable individuals who are not responsive to treatments. Under our queue assignment mechanism for the \texttt{Base} policy (Definition~\ref{def: sample dual-price queuing policy}), each person is assigned to the queue $t$ that maximizes the difference $\hat{m}^{t}_{n}(x) - \hat{\mu}^{\star,t}_{n}$. Since the dual price for \texttt{No Treatment} is $0$ (it has no budget constraint), a resource queue $t$ is eligible for assignment only when $\hat{m}^{t}_{n}(x) - \hat{\mu}^{\star,t}_{n} -\hat{m}^{0}_{n}(x)$ is at least larger than $0$. In this case, optimizing for expected outcomes will risk deprioritizing vulnerable individuals with low treatment effectiveness. However, we find this potential concern does not occur in our experimental results due to a strong correlation ($93\%$) between an individual's vulnerability and the treatment effectiveness of the best resource for that individual. This means that highly vulnerable individuals also tend to benefit most from treatment and are prioritized by our policy. 

In Figure~\ref{fig:vul_to_matching}, we plot the probability of being matched to any treatment queue under the \texttt{Base} policy as a function of individual vulnerability, where vulnerability is defined as the probability of an adverse outcome under \texttt{No Treatment}. We observe a clear pattern where individuals with higher vulnerability are increasingly likely to receive a resource, because of their higher treatment effects. In particular, more vulnerable individuals are also more likely to receive the more supportive PSH resources. Therefore, in settings where individuals who are very vulnerable also tend to have large benefits from receiving treatments optimizing for expected outcomes does not necessarily come at the expense of vulnerable individuals. 

\begin{figure*} [!htb]
\begin{center}
\includegraphics[scale=0.7]{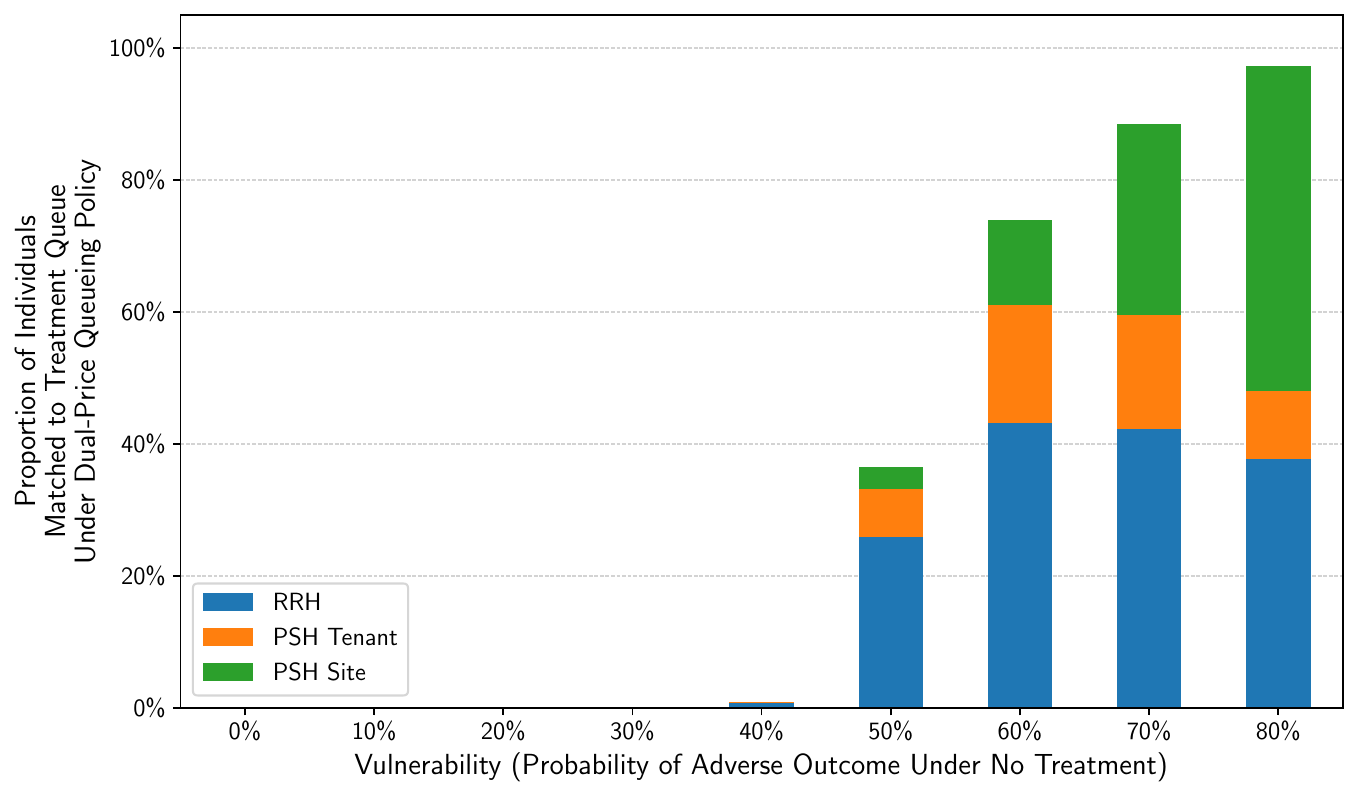}
\caption{\color{blue} Proportion of individuals matched to any treatment queue under a dual-price queueing policy, as a function of individual vulnerability. Vulnerability is defined as the probability of an adverse outcome under no treatment, where the x-axis corresponds to buckets with a width of $10\%$. For example, $0\%$ corresponds to the interval $[0\%, 10\%)$, $10\%$ to $[10\%, 20\%)$, and so on, where each bucket is left-inclusive and right-exclusive. The y-axis represents the average matching probability for individuals within each vulnerability bucket.}
\label{fig:vul_to_matching}
\end{center}
\end{figure*}

\subsubsection{Comparing Dual-Price Queueing Assignments to Historical Assignments Based on Characteristics} \label{subsec: policy_comparison}

In practice, policymakers will want to understand why our policy allocates resources differently from the existing \texttt{Historical} policy. Since our policy consistently makes resource assignments as a function of individual characteristics, we can examine individual covariates to understand cases where our \texttt{Base} policy makes different assignments than the \texttt{Historical} policy. We therefore examined every individual in our test data and isolated the cases where the two policies disagreed: either \texttt{Base} offered \emph{any} treatment while \texttt{Historical} offered none, or vice-versa.  (Observations on which both policies agreed were discarded.) 

Due to the high-dimensionality of our data (more than $100$ covariates), we used a simple feature-importance analysis approach to understand the disagreements. For individuals with policy disagreements, we assigned a label $1$ if \texttt{Base} assigned a treatment but \texttt{Historical} did not, and $0$ in the other case. Then we fit a logistic regression with one-hot encoded categorical covariates and normalized continuous covariates so that all covariates were in the same scale. Finally, we rank covariates by the absolute value of their coefficients as a heuristic for feature importance.

The results show that the covariates with the highest absolute value coefficients for predicting whether \texttt{Base} vs \texttt{Historical} assigned an treatment were covariates that did not come directly from the VI-SPDAT survey but from auxiliary HMIS data sources. These included characteristics like an individual's current employment type/status, whether they were looking for employment, or the specific category of living situation an individual was in. As an example, the \texttt{Historical} policy tended to assign treatments to individuals with some form of work whereas our policy was more likely to assign treatment to someone without known form of work. Since case managers primarily relied on the VI-SPDAT responses to make allocation decisions historically, our findings suggest that augmenting VI-SPDAT responses with additional HMIS data can improve allocation decisions. 
}

{\color{blue}
\section{Wait Time and Queue Length Statistics under the Dual-Price Queueing Policy} \label{sec: waittime simulation app}
Section \ref{app_sub: queue_discussion} discusses the theoretical properties of wait times and queue lengths under our proposed policy. Section \ref{app_sub: full_waittime_res} presents additional results on wait times observed in the numerical experiment described in Section~\ref{sec: policy_queue_waittimes}.
Section \ref{app_sub: waittime_causal_impact} investigates the causal effect of the wait times observed under our counterfactual policy on outcomes. Section \ref{app_sub: vi_spdat_changes} provides further empirical evidence from the VI-SPDAT that an individual's circumstances are unlikely to materially worsen under our policy's wait times.}

{\color{blue}
\subsection{Discussion of Arrival and Queue Statistics}\label{app_sub: queue_discussion}
In this section, we argue that the expected number of individuals assigned to the queue for any resource under the dual-price queuing policy matches the expected number of resource arrivals of that type, ensuring balance in the system over time. In constructing our policy optimization framework, we have not made any assumptions regarding the arrival processes of individuals and resources in the system. For the purposes of this discussion, we now assume that the arrival process of individuals is such that the expected number of arrivals in any interval of length $\ell$ is $\lambda \ell$, where $\lambda > 0$. Examples of such a process include cases where individuals arrive either with deterministic inter-arrival times of $1/\lambda$ or with i.i.d. exponential inter-arrival times with mean $1/\lambda$, i.e., a Poisson arrival process. Similarly, we assume that the arrival process of resources of type $t \in \mathcal{T} \setminus \{0\}$ is such that the expected number of arrivals in any interval of length $\ell$ is $\mu^t \ell$, where $\mu^t > 0$. Note that the system is assumed to be initially empty. To reflect the scarcity of resources, we assume that $\mu^t < \lambda$ for all $t \in \mathcal{T} \setminus \{0\}$. 

In our policy optimization framework, the parameter $b^t$, which represents the asymptotic capacity of treatment $t \in \mathcal{T} \setminus \{0\}$ per individual, can be related to the arrival processes through the equation $b^t = \mu^t/\lambda$, which denotes the expected number of arrivals of resources of type $t$ per unit time divided by the expected number of individual arrivals. The policy $\pi^\star$, defined in~\eqref{eq: policy from dual}, satisfies the capacity constraint $\mathbb{E}[\pi^{\star,t}(X)] \leq b^t$ for all $t \in \mathcal{T}$. Furthermore, under Assumption \ref{ass: meanoutcomes}, if resources of type $t \in \mathcal{T} \setminus \{0\}$ are very scarce and contribute to improving outcomes (which align with the housing application), these capacity constraints are binding for all $t \in \mathcal{T} \setminus \{0\}$.

The dual-price queueing policy introduced in Definition \ref{def: dual-price queuing policy} assigns each individual with covariates~$x$ to the queue for treatment $t$, for $t \in \mathcal{T}$, if $\pi^{\star,t}(x) = 1$. Given that the capacity constraints are binding for all $t \in \mathcal{T} \setminus \{0\}$, this implies that the expected number of individuals assigned to the queue for any treatment $t \in \mathcal{T} \setminus \{0\}$ at any time matches the expected number of arrivals of resources of type $t$, assuming that the covariates $X_k$ are i.i.d. across individuals. To see this, consider any time point $\ell$ and any treatment type $t \in \mathcal{T} \setminus \{0\}$. Under the dual-price queueing policy, the expected number of individuals assigned to the queue for resource $t$ at time $\ell$ is given by $b^t \lambda \ell$. This follows from the fact that the prior probability of assigning an arriving individual to the queue for treatment $t$ is $\mathbb{E}[\pi^{\star,t}(X)] = b^t$, and the expected number of individuals arriving in the system until time $\ell$ is $\lambda \ell$. Given the relation $b^t = \mu^t/\lambda$, this implies that the expected number of individuals assigned to the queue for resource $t$ by time $\ell$ matches the expected number of arriving resources of type $t$ up to the same time, i.e., $\mu^t \ell$. In this sense, the dual-price queueing policy ensures a well-balanced system. If this balance is not maintained (i.e., if the prior probability of assigning a treatment $t \in \mathcal{T} \setminus \{0\}$ is either low or high compared to that under $\pi^\star$), then, on expectation, the difference between the number of individuals waitlisted for treatment $t$ and the number of resources of type $t$ arriving in the system will continue to grow either toward infinity or negative infinity, indicating either longer waiting times for individuals or an increased number of unutilized resources—both of which are undesirable. So far, in the discussion, we have focused on the optimal dual-price queuing policy constructed under the assumption of knowing the distribution of the covariates and potential outcomes. Recall that, in reality, we only have access to data, meaning we can only compute and rely on the sample-based dual-price queuing policy, as defined in Definition~\ref{def: sample dual-price queuing policy}. Since we have shown this policy to be asymptotically optimal, it is reasonable to expect that the sample-based dual-price queuing policy will achieve a similar balance as the number of in-sample data points grows.

Studying the statistics of queue lengths for any treatment $t \in \mathcal{T} \setminus \{0\}$, or the waiting times of individuals and resources in the system, is significantly more complex and depends on further specifications about the arrival processes as well as the distribution of covariates. To illustrate this complexity, consider the simplest case where inter-arrival times are deterministic; that is, at any time $\ell$, $\lambda \ell$ (assume to be an integer) individuals arrive in the system, and $\mu^t \ell$ (assume to be an integer) resources of type $t \in \mathcal{T} \setminus \{0\}$ become available. As the dual-price queueing policy assigns each arriving individual with covariates $x$ to the queue for treatment $t \in \mathcal{T} \setminus \{0\}$ with probability $\pi^{\star,t}(x)$, the expected length of queue $t$ at time $\ell$ is given by
\begin{equation}\label{eq: expected queue length}
    \mathbb{E}\left[\left(\sum_{k=1}^{\lambda \ell} \pi^{\star,t}(X_k) - \mu^t \ell\right)^+\right],
\end{equation}
where $(a)^+ = \max\{a,0\}$. The term $\sum_{k=1}^{\lambda \ell} \pi^{\star,t}(X_k)$ represents the length of queue $t$ after $\lambda \ell$ arrivals of individuals with covariates $X_1, \dots, X_{\lambda \ell}$, while $\mu^t \ell$ is the number of resources of type $t$ that become available by time $\ell$. The expected queue length in equation~\eqref{eq: expected queue length} depends not only on the expected allocation of an individual to queue $t$ prior to knowing their covariates but also on the overall distribution of the covariates. Moreover, if we consider a stochastic arrival process, this quantity also depends on the specifics of that process (e.g., the distribution of inter-arrival times) and not merely on the expected number of arrivals, making the analysis significantly more complex. 

Additionally, in reality, we only have access to data and therefore lack knowledge of the true distribution of the covariates as well as the specifics of the arrival process. When assessing queue lengths, waiting times, and the number of idle resources under our proposed sample-based queueing policy, we will therefore rely on data to estimate. In Section \ref{app_sub: full_waittime_res}, we will indeed verify numerically that the allocation system under our sample-based queuing policy is well-balanced, in the sense that on average there is not a build up of individuals waiting in queues nor idle resources. 

}

\begin{figure*}[!htb]
\begin{center}
\includegraphics[scale=0.5]{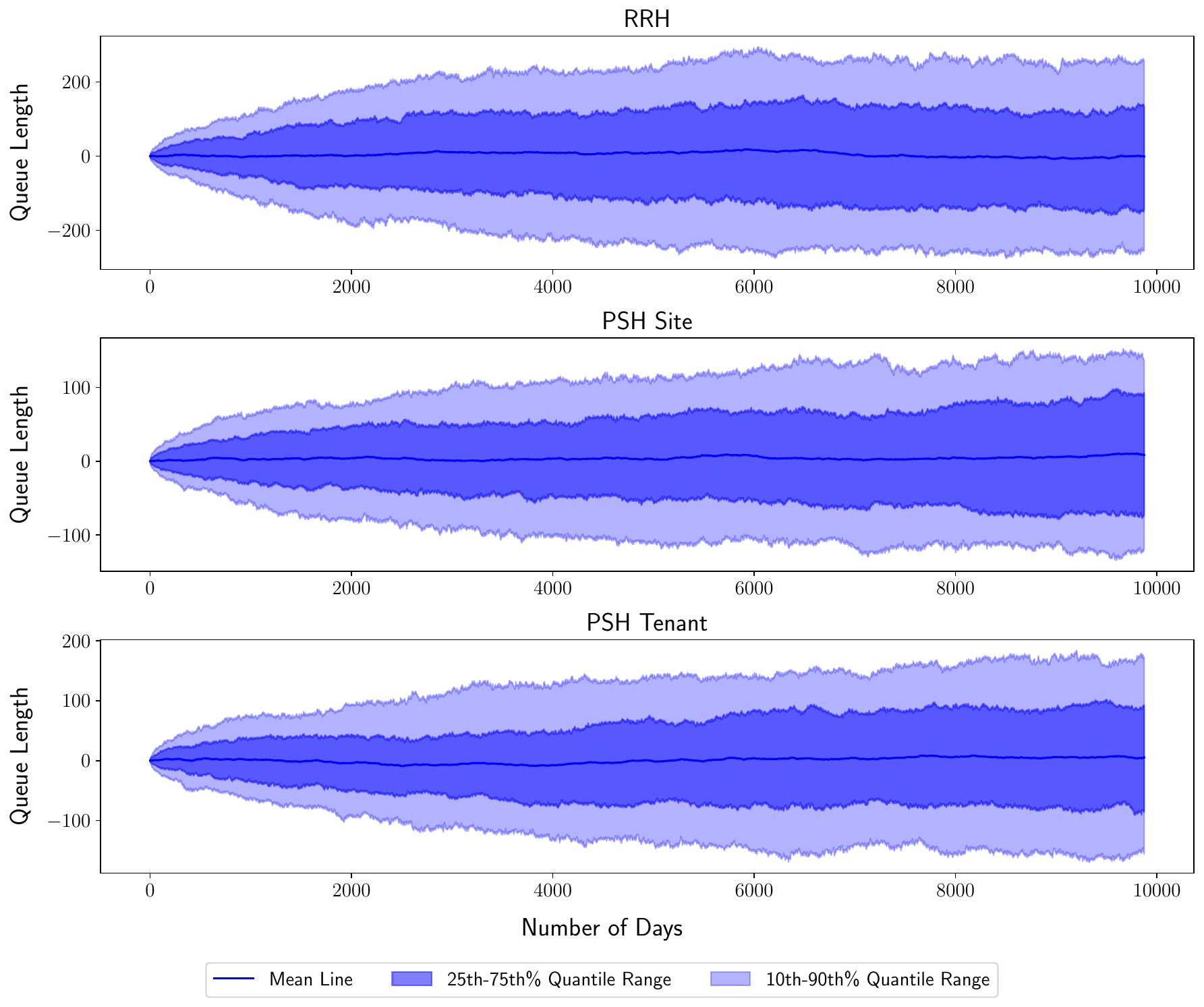}
\caption{Simulated treatment queue sizes across $250$ simulated sample paths, where positive queue size means individuals are waiting in queue while negative queue size means resources are waiting to be matched.}
\label{fig: queue_len_no_trunc}
\end{center}
\end{figure*}

{\color{blue}
\subsection{Full Wait Time Results for Section~\ref{sec: policy_queue_waittimes}}\label{app_sub: full_waittime_res}
{\color{blue}In this subsection, we present further numerical results on wait times (Section~\ref{app_sub: waittimes}) and queue lengths (Section~\ref{app_sub: queue_length}) under our policy, based on the experiment setup described in Section~\ref{sec: policy_queue_waittimes}.  

\subsubsection{Wait Times}\label{app_sub: waittimes}
Previously, we summarized our wait time results in Section~\ref{sec: policy_queue_waittimes} and Figure~\ref{fig:waittime_figure_maintext}, which showed that our policy results in relatively low wait times on average, across 250 trials.
One aspect that Figure~\ref{fig:waittime_figure_maintext} does not account for is idle housing ``waiting'' to be matched, which is also an undesirable situation based on our discussions with policy makers. Figure \ref{fig: waittime_no_trunc} presents the same statistics for the wait times of each $i^{\text{th}}$ arriving individual, adjusted by subtracting the idle time of their matched resource. This ``adjusted wait time'' can be computed as the difference between the individual’s arrival time and the arrival time of their matched resource, and provides a holistic view of the allocation system in terms of both individuals and resources. The ``adjusted wait time'' is positive when an individual is waiting for a resource and is negative when a resource is waiting to be matched to an individual, and we clarify that individuals and resources cannot both be waiting. Note that the quantities in~Figure \ref{fig:waittime_figure_maintext} are obtained by truncating negative adjusted wait times, which represent idle resource times, to zero. 
We see from Figure \ref{fig: waittime_no_trunc} that across each treatment queue, the mean adjusted wait times stabilize around $0$, indicating that the system is well-balanced on average, which aligns with the discussion in Section~\ref{app_sub: queue_discussion}. For applications in which a holistic view of the allocation system is important, we may desire such well-balanced behavior, where on average there is neither waiting individuals nor idle resources. We note that the $10^{\text{th}}$ quantiles (where adjusted wait times are negative and resources are idle) are $-57$, $-94$, and $-96$ days for \texttt{RRH}, \texttt{PSH Site}, or \texttt{PSH Tenant}, respectively. The $90^{\text{th}}$ quantiles (where adjusted wait times are positive and and individuals are waiting) are $57$, $99$, and $96$ days for \texttt{RRH}, \texttt{PSH Site}, or \texttt{PSH Tenant}, respectively.

}
\begin{figure*}[!htb]
\begin{center}
\includegraphics[scale=0.5]{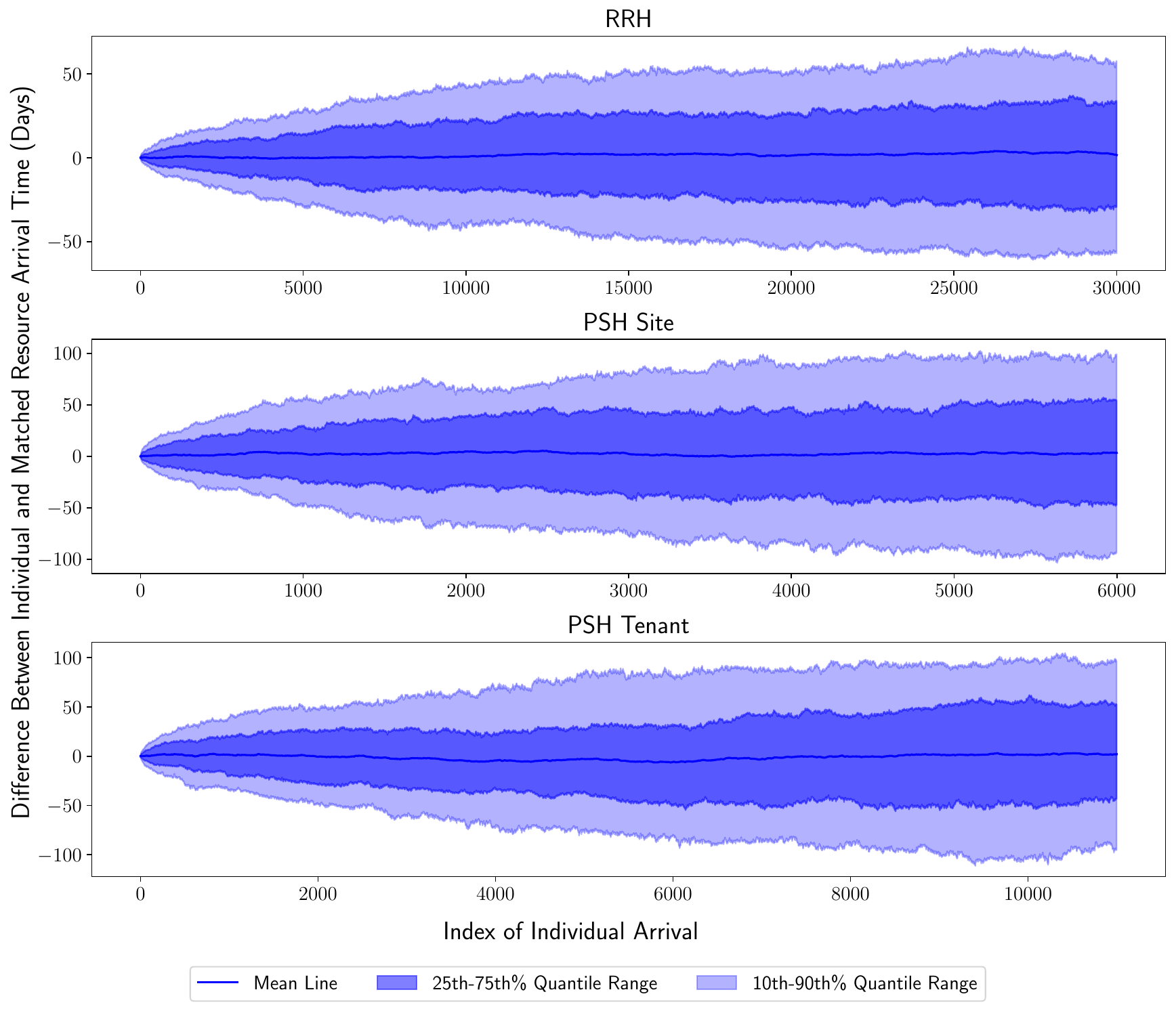}
\caption{Adjusted wait times for the different types of housing resources across $250$ sample paths, where, for each individual, we plot the difference between an individual's arrival time and the arrival time of their matched resource. A positive value means an individual arrives before their matched resource while a negative value means a resource arrived before its matched individual.}
\label{fig: waittime_no_trunc}
\end{center}
\end{figure*}


}
{\color{blue}
\subsubsection{Queue Lengths}\label{app_sub: queue_length}

To complement the results on wait times, we additionally investigate the queue lengths under our policy based on our experiment setup.
In Figure \ref{fig: queue_len_no_trunc}, we take a holistic system view by considering both individual and resource build-ups in treatment queues. We compute the queue size, displayed in this figure, as the difference between the number of individuals assigned and waiting for a particular housing type and the number of the respective resources waiting to be matched.
Therefore, at each time point, a positive queue size means there are individuals waiting for resource arrivals while a negative queue size means there are resources waiting to be matched to an individual. 
Similar to the wait time results in Figure \ref{fig: waittime_no_trunc}, we see that our policy, on average across simulation trials, results in a well-balanced system where neither individuals nor resources continuously build up within queues over time, but can result in longer queue lengths. The $10^{\text{th}}$ quantiles (where queue sizes are negative and there are resources waiting) are $-213$, $-118$, and $-156$ across the subplots, respectively, and the $90^{\text{th}}$ quantiles (where queue sizes are positive and individuals are waiting) are $256$, $137$, and $173$ across the subplots, respectively.
}

{\color{blue}\subsection{Impact of Wait Time on Outcomes} \label{app_sub: waittime_causal_impact}

One important consideration for whether the wait times observed under our policy are reasonable is the impact of these wait time levels on individual outcomes. Therefore, based on our historical data, we compute the negative impact of increasing wait times on individual outcomes for each group of individuals based on  disability characteristic and resources received. To study potential heterogeneity in how wait times may impact different individuals, we compute wait time conditional average treatment effects (CATE)---treatment referring to wait time in this case--- based on each individual's covariates using the Orthogonal/Double Machine Learning (DoubleML) framework from the causal inference literature \citep{chernozhukov2018double, nie2021quasi}. We chose this method because it allows for the incorporation of arbitrary machine learning techniques for estimation and naturally handles the continuous nature of wait times, in contrast to the binary or discrete treatment settings (recall that treatment refers to wait time in this analysis) commonly studied by other methods. For estimating the two nuisance functions required in DoubleML of outcomes and treatment assignment (in this case wait time) conditional on covariates and the final stage regression in DoubleML, we follow standard cross-validation and hyperparameter search procedures with the base estimators being linear/logistic regression (with and without regularization), random forests, gradient boosted trees, and XGBoost trees.

In Figure \ref{fig: waittime_impact_curves}, we show the mean, $25^{\text{th}}$ to $75^{\text{th}}$, and $10^{\text{th}}$ to $90^{\text{th}}$ quantile ranges of the decrease in probability of positive outcomes as wait times increase, where each subplot corresponds to a different combination of disability and resource received. In this case, a lower quantile (such as $10^{\text{th}}\%$) corresponds to larger decreases in positive outcome probability as wait times increases, while higher quantile corresponds to individuals where wait time has less effect on their positive outcome probability. For the sample path average \texttt{RRH} wait time of $18$ days, we see that the mean decrease in positive outcome probability is $-0.53$ p.p.and $-0.29$ p.p. for disabled and non-disabled individuals, respectively, while the $10^{\text{th}}$ percentile decreases by $-0.86$ p.p. and $-0.61$ p.p. for disabled and non-disabled individuals, respectively. For \texttt{PSH Tenant} and \texttt{PSH Site}, the sample path average wait time of $30$ days corresponds to mean decreases of $-0.64$ and $-0.36$ p.p. respectively, while the $10^{\text{th}}$ percentile individual's decrease in positive outcome probability are $-0.99$ and $-0.87$ p.p., respectively. Even at $90^{\text{th}}$ percentile wait times of $55$, $100$, and $100$ days of wait time for \texttt{RRH}, \texttt{PSH Tenant}, and \texttt{PSH Site}, the mean decrease in positive outcome probabilities are $-1.62$, $-2.12$, $-1.21$, and $-0.90$ p.p. for Disabled \texttt{RRH}, \texttt{PSH Tenant}, \texttt{PSH Site}, and Non-Disabled \texttt{RRH} individuals respectively. In addition, at $90^{\text{th}}$ percentile wait times, the $10^{\text{th}}$ percentile individual's decrease in positive outcome probability are $-2.62$, $-3.31$, $-2.89$, and $-1.85$ p.p. across the subplots, respectively. As expected, increasing wait times has a negative effect on positive outcome probabilities across individuals, but the impact is relatively minimal at the wait times seen in Figure~\ref{fig:waittime_figure_maintext}. 

\begin{figure*}[!htb]
\begin{center}
\includegraphics[scale=0.5]{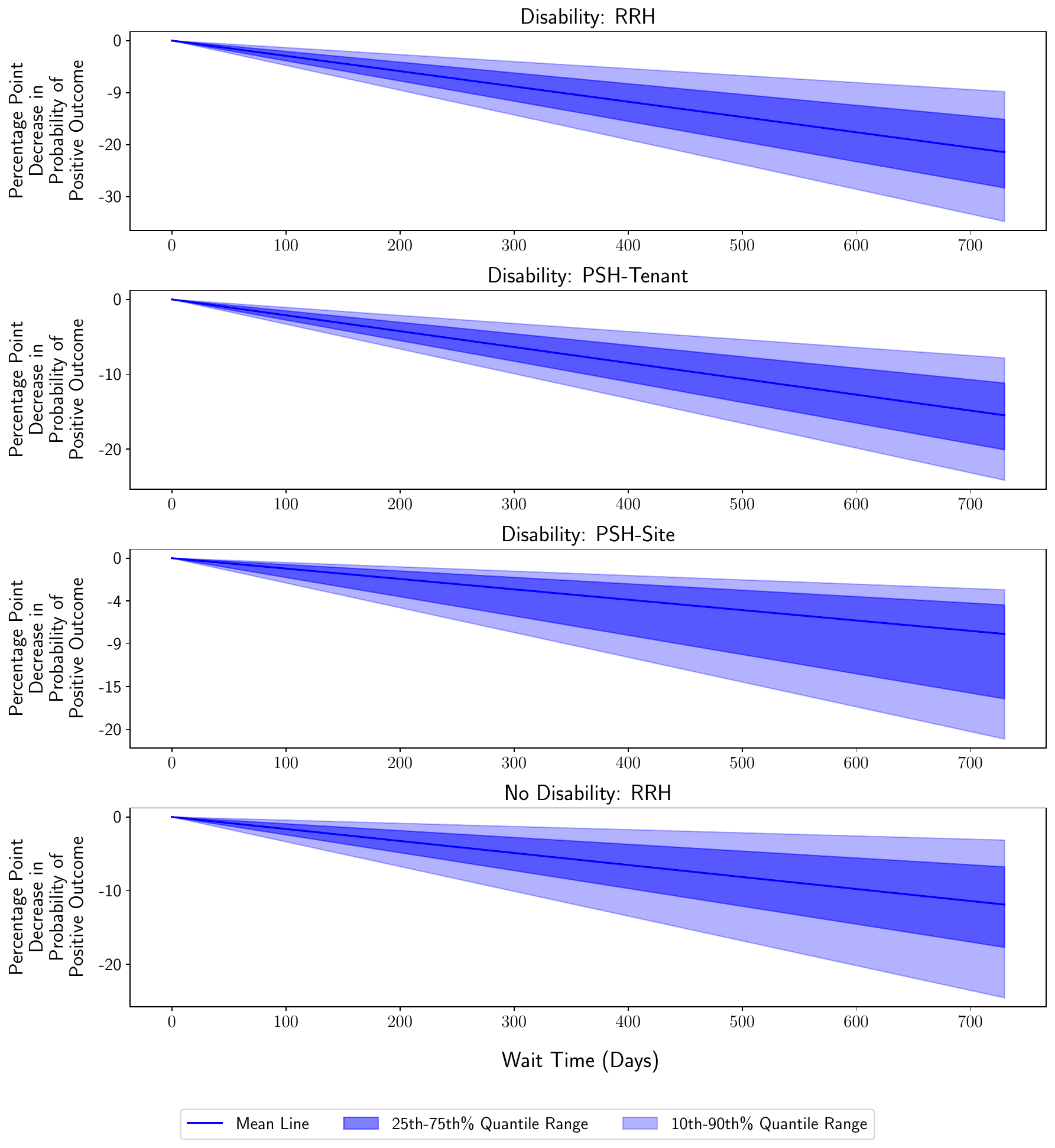}
\caption{DoubleML/R-Learner causal estimate of the heterogeneous negative impact of increasing wait times on probabilities of positive outcomes.}
\label{fig: waittime_impact_curves}
\end{center}
\end{figure*}
}

{\color{blue}
\subsection{Impact of Time on Vulnerability: Evidence from VI-SPDAT Reassessments} \label{app_sub: vi_spdat_changes}
Another way to gauge the potential impact of the wait times (Sections~\ref{sec: policy_queue_waittimes},~\ref{app_sub: full_waittime_res}, and~\ref{app_sub: waittime_causal_impact}) of our policy on outcomes is to examine changes in an individual's circumstances over time under no treatment, which we can observe from retakes of the VI-SPDAT. Because case managers typically administer a retake of the VI‑SPDAT only when they notice a visible shift in circumstances, the frequency of retakes, the time elapsed between them, and the magnitude of any score change together provide a proxy for how quickly vulnerability actually worsens. For example, if individuals rarely retake the assessment within $100$ days (the $90^{\rm th}$ quantile of wait time for PSH under our policy) and their scores do not increase significantly (indicating worsening vulnerability), then this suggests wait times at the levels of our policy do not result in material changes to circumstances. 

To analyze individual changes based on VI-SPDAT retakes, we subset our data on those who received no treatment historically since these individuals represent what happens to someone as they wait for their allocation. We see that only $6.7\%$ of all individuals had a reassessment within $100$ days. Among those who ever retook the assessment, the average and median time to retake were $385$ and $303$ days, respectively. These numbers indicate that case workers did not find most individuals' circumstances to change materially within the time frame of the wait times of our policy to administer a reassessment. Furthermore, amongst those who needed a reassessment within a year, only $51\%$ of the individuals had a score increase (their circumstances worsened) and the mean and median score changes both increased by $1$ on an $18$ point scale. Consequently, the wait times our policy introduces are short relative to the timescales on which VI‑SPDAT scores (a proxy for individuals' situations) tend to change, suggesting the practical acceptability of our policy’s queuing wait times.
}

\end{document}